\documentstyle{article}
\input amssymb.sty
\input epsf
\newtheorem{theorem}{Theorem}[section]
\newtheorem{lemma}[theorem]{Lemma}
\newtheorem{proposition}[theorem]{Proposition}
\newtheorem{corollary}[theorem]{Corollary}
\newtheorem{conjecture}[theorem]{Conjecture}

\begin{document}
\newcommand{\Z}{{\Bbb Z}}
\newcommand{\R}{{\Bbb R}}
\newcommand{\Q}{{\Bbb Q}}
\newcommand{\C}{{\Bbb C}}
\newcommand{\lra}{\longrightarrow}
\newcommand{\lms}{\longmapsto}

\begin{titlepage}
\title{  The dihedral Lie algebras and    Galois symmetries of 
$\pi^{(l)}_1({\Bbb P}^1 - (\{0,  \infty \} \cup \mu_N))$}
\author{A.B. Goncharov }

 \date{}
\end{titlepage}
\stepcounter{page}
\maketitle

\tableofcontents 

\section{An outline}

    This is the next in the series of papers [G1-3] 
devoted to 
study of higher cyclotomy, understood as motivic theory of 
multiple polylogarithms at roots of unity,  
and its relationship with modular varieties for $GL_m/_{\Q}$, for all 
$m \geq 1$.  
Let $\mu_N$ be the group of $N$-th roots of unity. 
Our main objective is a mysterious link between 
the {\it structure of the 
motivic fundamental group of } 
$$
X_N:= \quad {\Bbb P}^1 - (\{0, \infty\} \cup  \mu_N) \quad = \quad 
{\Bbb G}_m - 
\mu_N 
$$ 
and the 
 {\it geometry and topology of the following 
modular varieties for $GL_m/_{{\Q}}$}:   
\begin{equation} \label{7.19.00.1}
\Gamma_1(m;N) \backslash GL_m(\R)/ \R_+^*\cdot O_m \qquad \mbox{for  $m>1$} 
\end{equation}
where 
 $\Gamma_1(m;N) \subset GL_m(\Z)$ is the  subgroup 
stabilizing the vector $(0, ... , 0, 1)$ mod $N$. For $m=1$ it is 
$
S_N:=  {\rm Spec}\Bigl(\Z[\zeta_N][\frac{1}{N}]\Bigr)
$. 
Adelic approach provides a coherent description for all $m$, see s. 2.2.

In the present paper we turn to the Galois side of the story. However  
to keep the motivic perspective  let us recall the following. 
According to Deligne [D] the motivic fundamental group of a variety 
is not just a group, 
 but rather a Lie algebra object in the category of mixed motives.  
It can be viewed as a pronilpotent completion of the topological fundamental group 
of the corresponding complex variety equipped with lots of additional structures of 
analytic, geometric and arithmetic nature. 

Mixed motives are algebraic geometric objects. They can be seen 
through their realizations. The two most popular realizations  
which hypothetically capture all the information about the category 
 of the mixed motives are: 

{\it the Hodge realization} provided by analysis   
and  algebraic geometry, and

{\it the $l$-adic realization}
 provided by arithmetic and  algebraic geometry.

The $l$-adic realization of the motivic fundamental group is obtained from the 
action of the absolute Galois group ${\rm Gal}(\overline \Q/ \Q)$ on the pro-$l$ 
completion 
$\pi^{(l)}_1(X)$
of the fundamental group of the variety $X$. 
That is why  in this paper we study the action of 
the Galois group on 
$\pi^{(l)}_1(X_N)$. 

In the case $N=1$ 
this problem has been addressed by Grothendieck [Gr], Deligne [D], 
Ihara [Ih0-3] , Drinfeld [D] 
and others - see the wonderful survey of Ihara  [Ih0].  
For $N>1$ it has not been investigated. Our point of view  is that we should 
study the problem for all $N$, 
penetrating  the 
structures independent of $N$ 
(like modular complexes, see s. 2.5). 
%and then go to the limit using the natural maps 
%$X_{NM} \lra X_N$ given by $z \lms z^M$. 
%The last step will be discussed elsewhere. 

Briefly our approach is this. 
Let $G$ be a commutative group. In [G3] we constructed    a bigraded Lie 
algebra 
$ {D}_{\bullet \bullet}(G) $,  
called the dihedral Lie algebra of $G$, see section 3 below.   
We relate $ {D}_{\bullet \bullet}(\mu_N) $
to the Lie algebra of the image of ${\rm Gal}(\overline  \Q/\Q)$ in 
${\rm Aut}\pi^{(l)}_1(X_N)$. 
On the other hand in 
 [G2-3, 5]  
the structure of the Lie algebra $ {D}_{\bullet \bullet}(\mu_N)$ 
is related to  the 
 geometry of modular varieties (\ref{7.19.00.1}). 
Using these results    
we study the action of  
 ${\rm Gal}(\overline  \Q/\Q)$ on $\pi^{(l)}_1(X_N)$. 
In particular we obtain new results about the action 
on $\pi^{(l)}_1({\Bbb P}^1 - \{0, 1, \infty\})$. 

Here is a more detailed account. 
Let $X$ be a regular curve over $\overline \Q$. Let  
$\overline X$ be the corresponding projective curve and $v_x$ 
a nonzero tangent vector at a point 
$x \in \overline X$. 
Then according to Deligne [D] one can define the geometric 
profinite fundamental group $\widehat \pi_1(X, v_x)$ based at the vector $v_x$.  
If $X$, $x$ and $v_x$ are defined over a number field $F \subset \overline \Q$  
then the group $ {\rm Gal}(\overline  \Q/F) $ acts by automorphisms of $\widehat \pi^{(l)}_1(X, v_x)$. 
If $X = X_N$, there is a tangent vector $v_{\infty}$ at $\infty$ corresponding to 
the inverse $t^{-1}$ of the canonical coordinate $t$ on 
${\Bbb P}^{1} - (\{0,\infty \} \cup \mu_N)$. 

Since any finite $l$-group is nilpotent the pro-$l$ group 
$\pi^{(l)}_{1}(X_N, v_{\infty})$ is pronilpotent. 
Let ${\Bbb L}^{(l)}_N = {\Bbb L}^{(l)}(X_N, v_{\infty})$  be the $l$-adic 
pro-Lie algebra  
corresponding via the Maltsev theory to  
$\pi^{(l)}_{1}(X_N, v_{\infty})$ 
(see ch. 9 of [D]). It is a free pronilpotent Lie algebra with 
generators corresponding to the  loops around $0$ and all $N$-th 
roots of unity. 
We  call it  the $l$-adic fundamental Lie algebra of $X_N$. 
The Galois group acts by its automorphisms. 

Let $\zeta_n$ be a primitive $n$-th root of unity, and  $\Q(\zeta_{l^{\infty}N}):= \cup
\Q(\zeta_{l^{n}N})$. For the reasons explained in s 3.1 we restrict 
the action of the Galois group to the subgroup 
${\rm Gal} (\overline{\Bbb Q}/ {\Bbb Q}(\zeta_{l^{\infty}N}))$, picking up the homomorphism 
$$
{\rm Gal} \Bigl(\overline{\Bbb Q}/ {\Bbb Q}(\zeta_{l^{\infty}N})\Bigr) \lra {\rm Aut}{\Bbb L}^{(l)}_N
$$
Let ${\rm Der} {\Bbb L}^{(l)}_N$ be the Lie algebra of all derivations of the Lie algebra ${\Bbb L}^{(l)}_N$. 
Linearizing, as explained in s. 2.1 or 3.2,  the above  map
we get the Lie algebra
\begin{equation} \label{8.5.00.1}
{\cal G}_N^{(l)} \quad \hookrightarrow \quad {\rm Der} {\Bbb L}^{(l)}_N
\end{equation}

 The fundamental  Lie algebra ${\Bbb L}^{(l)}_N$ is equipped with two filtrations 
preserved by the Galois action. 
The weight filtration can be defined 
on the fundamental Lie algebra of any algebraic variety over $\Q$. In 
our case it coincides with the lower central series of the Lie algebra 
${\Bbb L}^{(l)}_N$. 
The depth filtration is more specific. It is given by the lower central series of 
the codimension one ideal
$$
{\cal I}_N:= \quad {\rm Ker}\Bigl( {\Bbb L}^{(l)}(X_N, v_{\infty})
 \lra {\Bbb L}^{(l)}({\Bbb G}_m, v_{\infty}) = \Q_l(1) \Bigr)
$$
where the map is provided by the natural inclusion $X_N \hookrightarrow {\Bbb G}_m$. 

These filtrations induce filtrations on 
${\rm Der} {\Bbb L}^{(l)}_N$, and hence, via (\ref{8.5.00.1}), on  ${\cal G}_{N}^{(l)}$. 
The associated graded for the weight and depth filtrations 
${\rm Gr}{\cal G}_{\bullet \bullet }^{(l)}(\mu_N)$ is a Lie algebra bigraded by negative 
integers $-w$ and $-m$. We call it the {\it level $N$ Galois Lie algebra}. 
  When $N=1$ it is 
denoted  by ${\rm Gr}{\cal G}_{\bullet \bullet }^{(l)}$. 

The weight filtration on ${\Bbb L}_N^{(l)}$ admits a splitting, i.e. is defined by a grading, compatible 
with the depth filtration. Such a weight splitting is provided by the 
eigenspaces of a Frobenius 
$F_p \in {\rm Gal} (\overline \Q/\Q), (p \not | N)$. Therefore taking ${\rm Gr}$ 
for the weight and depth filtrations we get an embedding 
$$
{\rm Gr}{\cal G}_{\bullet \bullet }^{(l)}(\mu_N) \quad \hookrightarrow \quad 
{\rm Gr}{\rm Der}{\Bbb L}^{(l)}_N
$$

The vector space ${\rm Gr}{\cal G}_{-w, -m}^{(l)}(\mu_N)$ is nonzero only if $w \geq m \geq 1$. As a ${\rm Gal} 
({\Bbb Q}(\zeta_{l^{\infty}})/ {\Bbb Q})$-module it is isomorphic to a 
direct sum of copies of 
$\Q_l(w)$. 

For every integer $m \geq 1$ we have the depth $m$ quotient 
$$
{\rm Gr}{\cal G}^{(l)}_{\bullet, \geq -m}(\mu_N) := \quad \frac{{\rm Gr}{\cal G}^{(l)}_{\bullet \bullet}(\mu_N)}{ 
{\rm Gr}{\cal G}^{(l)}_{\bullet, < -m}(\mu_N)}
$$
 It is  a  bigraded pro-nilpotent Lie algebra of the nilpotence class $\leq m$. 
In particular 
 ${\rm Gr}{\cal G}^{(l)}_{\bullet, \geq -1}(\mu_N)$ is abelian. 

The direct sum of the components of the Galois Lie algebra 
satisfying weight $=$ depth condition is a Lie subalgebra 
${\rm Gr}{\cal G}^{(l)}_{\bullet}(\mu_N)$ called 
the {\it diagonal Galois Lie algebra}.

\begin{center}
\hspace{4.0cm}
\epsffile{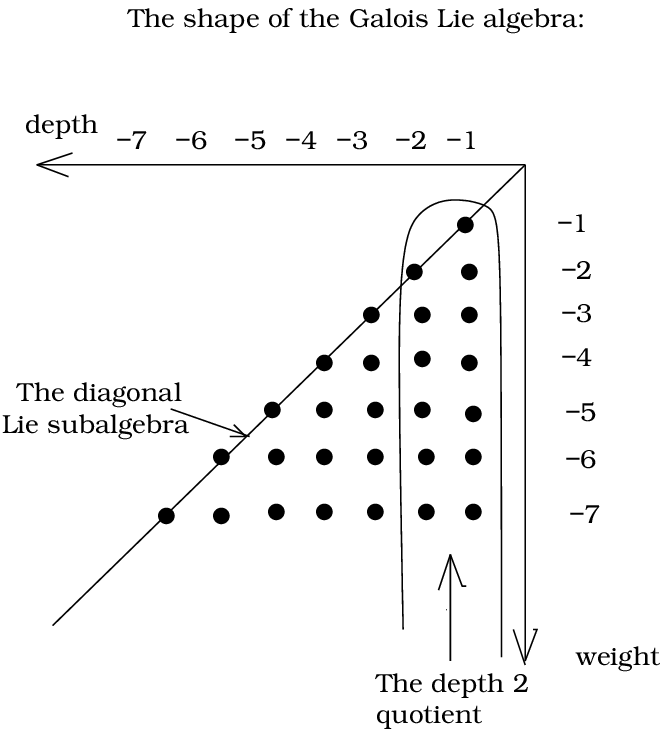}
\end{center}

Recall the standard cochain complex of a Lie algebra ${\cal G}$  
$$
{\cal G}^{\vee} \stackrel{\delta}{\lra} \Lambda^2 {\cal G}^{\vee} 
\stackrel{\delta}{\lra} \Lambda^3 {\cal G}^{\vee} \lra ... 
$$
where the first differential is dual to the commutator map $[ , ]: \Lambda^2 {\cal G} \lra {\cal G}$, 
and the others are obtained using 
the Leibniz rule. The condition $\delta^2 =0$ is equivalent to the Jacobi identity. 
If the Lie algebra ${\cal G}$ is graded its standard cochain complex inherits the grading. 

Let $V_m$ be the standard $m$-dimensional representation of $GL_m$. 
Our main goal is to show that 
\vskip 3mm
{\it the depth $m$, weight $w$ part of 
the standard cochain complex of the level $N$ Galois Lie algebra is related to 
the geometry of the local system with the fiber $S^{w-m}V_m$ over 
the level $N$ modular variety for the group ${GL_m}_{/{\Q}}$. 

In particular the structure of the depth $m$ quotient of the Galois Lie algebra is described by geometry of the modular variety for ${GL_m}_{/{\Q}}$}. 
\vskip 3mm
For $m =1$ 
this is deduced from the motivic theory of classical polylogarithms 
developed by Deligne and Beilinson 
([D], [Be], [BD], [HW]). 
%In the special case $m=1, N=1$ this has been known thanks 
%to the efforts of Soule [%So], Deligne [D] and Ihara [Ih2]. 
If $m=1$ we deal with the scheme $S_N$,  
 emphasizing that it is a modular variety for ${GL_1}_{/{\Q}}$. 

In this paper we prove general results about the Galois action on $\pi_1^{(l)}(X_N, v_{\infty})$, and  combining them 
with the results of [G2-3]   
 establish the relationship above for $m=2$ 
and essentially for $m=3$. The results [G5,7]  
indicate a similar story for $m=4$. 
 
For example in the case $N=1$ our main results describe the structure of the 
Lie algebra ${\rm Gr}{\cal G}^{(l)}_{\bullet, \geq -m}$  
very explicitly for $m=2,3$: for  $m=2$ it is given  
 in terms of  the classical modular triangulation of the hyperbolic plane 
(see the picture in s. 2.3) 
and for $m=3$ 
via a similar $GL_3(\Z)$-equivariant structure of the symmetric space 
${\Bbb H}_3 := GL_3(\R)/O(3) \cdot \R_+^*$ defined in [G3] using Voronoi's 
decomposition of ${\Bbb H}_3$.  
A detailed exposition of these results see in section 2. 

We expect a similar description for all $m$ 
in terms of the rank $m$ modular complexes, see 
conjecture \ref{dwhn1*} below. For a reformulation without modular complexes 
see conjecture 
\ref{ramier}. Theorem 1.2 in [G3] shows that 
they are equivalent.

Let me say few words about the methods. The Lie algebra 
${\cal G}^{(l)}_N$ enjoys the following 
properties:

i) ${\cal G}^{(l)}_N$ acts by derivations of the Lie algebra 
${\Bbb L}^{(l)}_N$. 

ii) There is a  subspace $\Q_l(1) \subset 
{\Bbb L}^{(l)}_N$ trivial as a Galois module 
(``canonical loop around infinity''). 

iii) The Galois action preserves conjugacy classes of 
loops around $0$ and $\mu_N$. 

iv) The group $\mu_N$ acts on 
${\Bbb L}^{(l)}_N$ 
commuting with the Galois action, see s. 3.1. 

v) The Galois action is compatible 
with the maps $X_{NM} \lra X_{M}$ 
given by $x \lms x$ and $x \lms x^N$, see s. 5.7. 

We construct explicitly 
in s. 5 the Lie subalgebra of ${\rm Gr}{\rm Der}{\Bbb L}^{(l)}_N$ 
respecting all these properties. One of our key points is that the Lie algebra 
${\rm Gr}{\cal G}^{(l)}_N$ should lie in a smaller Lie subalgebra which we single out by imposing 

(!) The ``power shuffle'' relations. 

To explain the origin of these relations 
let me recall the multiple polylogarithm functions [G9, 3]:
$$
Li_{n_1, ..., n_m}(x_1, ..., x_m):= \sum_{0 < k_1 < ... < k_m}\frac{x_1^{k_1} ... x_m^{k_m}}{k_1^{n_1} ... k_m^{n_m}}
$$
Their immediate  property is the 
shuffle product formula, which in the 
simplest case reads as follows:
$$
Li_{n}(x) Li_{m}(y) \quad = \quad Li_{n,m}(x,y) + Li_{n+m}(xy)  + Li_{m, n}(y,x)
$$
Indeed, 
$$
\sum_{k_1, k_2 >0}\frac{x_1^{k_1} x_2^{k_2}}{k_1^{n_1}  k_2^{n_2}} = 
\Bigl(\sum_{0 < k_1 < k_2} + \sum_{k_1= k_2 >0} + \sum_{0 < k_2 < k_1 }\Bigr)\frac{x_1^{k_1} x_2^{k_2}}{k_1^{n_1}  k_2^{n_2}}
$$
Similar arguments obviously lead to the general formula. 
Our most nontrivial constraint  (!) on the image of the Galois group 
is an $l$-adic/motivic version of these relations, considered modulo the 
lower depth terms, which allows to avoid complications related to 
regularization of the divergent relations. 
Although the proof of shuffle product relations for multiple 
polylogarithms is so simple, it is rather difficult to deal with 
their $l$-adic/motivic version. We proved in 
s. 5 that imposing relations (!) and  i)-v) we get a Lie 
subalgebra in ${\rm Gr}{\rm Der}{\Bbb L}^{(l)}_N$. Its precise description 
goes  as follows. 
We construct in s.  5(see theorem \ref{ga}c)) an embedding 
$$
\xi_{\mu_N}: D_{\bullet \bullet}(\mu_N) \quad \hookrightarrow \quad {\rm Gr}{\rm Der}{\Bbb L}_N
$$
The Lie subalgebra 
$
\xi_{\mu_N}(D_{\bullet \bullet}(\mu_N))
$  
coincides with the Lie algebra described 
by the six properties above.  If we assume relations imposed by 
ii)-iv)  there is an amazing duality 
between the relations (!) and  i),  which is naturally build in the definition of the dihedral Lie algebras, see Remark in s. 4.3.

\begin{conjecture} \label{ramierz} 
$
{\rm Gr}{\cal G}_{\bullet \bullet}^{(l)}(\mu_N) \subset \xi_{\mu_N}({D}_{\bullet \bullet}(\mu_N)) \otimes_{\Q} {\Q_l}
$. 
\end{conjecture}

If $G$ is a trivial group we set ${D}_{\bullet \bullet} : = {D}_{\bullet \bullet}(\{e\})$, 
 and $\xi:= \xi_{\{e\}}$. 
\begin{conjecture} \label{ramier}
One has $\xi({D}_{\bullet \bullet})\otimes \Q_l  = {\rm Gr} {\cal G}^{(l)}_{\bullet \bullet} $.
\end{conjecture}
Since the Lie algebra ${D}_{\bullet \bullet}$ is defined {\it very} explicitly, this  would give a 
precise description of the associated graded of the image of the Galois group.

The $(-w, -m)$-component of  
$ 
{D}_{\bullet \bullet}(G)
$ 
can be nonzero only 
if $w \geq m \geq 1$. 
The {\it diagonal} Lie algebra ${D}^{\Delta}_{\bullet}(G)$
is the subalgebra  of ${D}_{\bullet \bullet}(G)$ formed by the 
``depth $=$ weight'' components. It is graded by the weight. 

\begin{conjecture} \label{CD} Let $p$ be a prime number. Then 
$
{\rm Gr}{\cal G}^{(l)}_{\bullet}(\mu_p)
= \xi_{\mu_p}({D}^{\Delta}_{\bullet}(\mu_p)) 
$. 
\end{conjecture}

In this paper we prove conjecture \ref{ramierz} for the depth $2$ quotient 
in s. 7.3,  
and for the depth $3$ quotient in the $N=1$ case in s. 7.4. Its complete proof will appear elsewhere.
We prove conjectures \ref{ramier} and \ref{CD} 
in the depth 2 and 3. The relationship with geometry of modular varieties comes out of blue via the results of [G2-3, 5], and essentially used in the proofs. 
The  results of [G5,7] 
indicate a possibility to extend the proofs for the depth 4.  

We conjecture that 
the level $N$ Galois 
Lie algebra is ``very close'' to the image of the dihedral Lie algebra of $\mu_N$. 
For example the asymptotics  of dimensions of their $(w,m)$-components when 
$N \to \infty$ 
should coincide.

For $N=1$ there is another Lie subalgebra, the  
Grothendieck-Teichmuller Lie algebra ${\cal G}rt$ defined by Drinfeld [Dr],   
sitting  inside of ${\rm Der}{\Bbb L}({\Bbb P}^1 \backslash \{0, 1, \infty\}, v_{\infty})$. 

\begin{conjecture} The Lie subalgebra $\xi({D}_{\bullet \bullet})$ coincides 
with the associated graded for the depth filtration 
of the Grothendieck-Teichmuller Lie algebra ${\cal G}rt_1$.  
\end{conjecture}

We do not even know that one of these Lie algebras contains the other. 

We summarize the relationship between these Lie algebra in the following diagram, 
where the top left arrow assumes conjecture \ref{ramierz}, and all the other 
arrows are embeddings. 
$$
\begin{array}{ccc}
\xi({D}_{\bullet \bullet})\otimes \Q_l & & \\
\cup & \searrow &\\
{\cal G}^{(l)}& \hookrightarrow & {\rm Gr}^D{\rm Der}{\Bbb L}^{(l)}({\Bbb P}^1 \backslash \{0, 1, \infty\}, v_{\infty})\\
\cap &\nearrow &\\
{\rm Gr}^D {\cal G}rt_1 \otimes \Q_l&  &
\end{array}
$$

In the next section, which continues the introduction,  
we formulate our results in detail.

\section{Formulations of the main results}

{\bf 1. The Lie algebras of Galois symmetries of $\pi^{(l)}(X_N)$}. 
For any group $H$, denote by $H(m)$ its lower central series for  $H$:
$$
H(1):= H,  \quad H(m+1) := [H(m), H]
$$
Here $[,]$ stands for the closure of the commutator subgroup. 

The group $\pi^{(l)}_1(X_N):= 
\pi^{(l)}_1(X_N, v_{\infty})$ has two filtrations by normal subgroups, indexed by integers $n \leq 0$: 

i) {\it The weight filtration}  ${\cal F}^W_{\bullet}$ by the lower central series: 
$$
{\cal F}^W_{-w}\pi^{(l)}_1(X_N) := \quad \pi^{(l)}_1(X_N)(w)
$$

Let $I_N$ be the kernel of the map 
$
\pi^{(l)}_1(X_N) \lra 
\pi^{(l)}_1({\Bbb G}_m )
$
 induced by the embedding 
$X_N \hookrightarrow {\Bbb G}_m $. 

ii) {\it The depth filtration} ${\cal F}^D_{\bullet}$ is given 
by the lower central series $I_N(m)$ for $I_N$:
$$
{\cal F}^D_{0}\pi^{(l)}_1(X_N) := \pi^{(l)}_1(X_N), \qquad {\cal F}^D_{-m}\pi^{(l)}_1(X_N) := I_N(m)
$$ 
Notice that $I_N$ is not profinitely generated, however since we will always divide by
${\cal F}^W$, this causes no problem. 

These filtrations provide a projective system of quotients 
\begin{equation} \label{4-11} 
\pi^{(l)}_1(X_N)_{[w, m]}:= \quad 
\pi^{(l)}_1(X_N)/{\cal F}^W_{-w-1} \cdot {\cal F}^D_{-m-1}
 \end{equation}
of the 
group $\pi^{(l)}_1(X_N)$. 
They are unipotent  $l$-adic Lie groups. 
The Galois group acts by automorphisms of (\ref{4-11}), so we are picking up  homomorphisms
\begin{equation} \label{4-4} 
{\rm Gal} (\overline{\Bbb Q}/ {\Bbb Q}) \lra {\rm Aut} 
\Bigl( \pi^{(l)}_1(X_N)_{[w+1, m+1]}\Bigr)
 \end{equation}

The images of maps (\ref{4-4}) sit in the middle of short exact sequences
$$
0 \lra U^{(l)}_{N; [w, m]} \lra G^{(l)}_{N; [w, m]} \lra {\rm Gal} 
\Bigl({\Bbb Q}(\zeta_{l^{\infty}N})/ {\Bbb Q}\Bigr)\lra 0 
$$
In particular  the map  (\ref{4-4}) for $m=w=1$ is provided by the action of the Galois group on 
$H_{et}^1(X_N\otimes_{\Q}\overline \Q, \Z_l)(1)$.  So its image is ${\rm Gal} 
({\Bbb Q}(\zeta_{l^{\infty}N})/ {\Bbb Q})$, and $U^{(l)}_{N; [0, 0]}$ is the trivial group.  

 The weight and depth filtrations induce the two filtrations on the automorphism group of ${\Bbb L}^{(l)}_N$. Taking the quotients
of the image of the Galois group with respect to these filtrations we come to the groups 
$U^{(l)}_{N; [w, m]}$. 
The associated graded for the weight filtration is isomorphic to the free Lie algebra 
generated by $H_{et}^1(X_N\otimes_{\Q}\overline \Q, \Q_l)(1)$. Since 
${\rm Gal} (\overline{\Bbb Q}/ {\Bbb Q}(\zeta_{l^{\infty}N}))$ acts trivially on it $U^{(l)}_{N; [w, m]}$ is a unipotent $l$-adic Lie group  of the nilpotence class $\leq m$. 

There are natural projections
$$
 U^{(l)}_{N; [w, m]} \lra U^{(l)}_{N; [w', m']} \quad 
\mbox{if $w \geq w', m \geq m'$} $$
These groups  form  a projective system. 
Denote by ${\rm Lie} (*)$ the corresponding Lie algebras over $\Q_l$, and 
consider the pronilpotent Lie algebra 
$$
{\cal G}_N^{(l)}:= \lim_{\leftarrow} {\rm Lie}U^{(l)}_{N; [w, m]}
$$
It is, by construction, bifiltred by the {\it weight} $w$ and {\it depth} $m$. 
The associated graded 
${\rm Gr}{\cal G}_{\bullet \bullet }^{(l)}(\mu_N)$ is a Lie algebra bigraded by negative 
integers $-w$ and $-m$.

{\bf 2. The depth 1 case as part of the general picture}. Recall that the depth $1$ quotient 
 of the Galois 
 Lie algebra is abelian. 
Therefore the  result below, which follows from the motivic theory of 
classical polylogarithms developed by Deligne and Beilinson ([D], [BD], [HW]), 
settles the $m=1$ case of our story. 

\begin{theorem} \label{4-12.60090}
There is canonical isomorphism 
$$
{\rm Gr}{{\cal G}^{(l)}_{-w, -1}}(\mu_N) \quad = \quad
{\rm Hom}\Bigl(K_{2w-1}(S_N),  \Q_l\Bigr)
$$ 
\end{theorem}
The right hand side is calculated by  Borel's theorem [B].  
The result fits in our general framework as follows.  Then 
\begin{equation} \label{7.8.00.3} 
K_{2w-1}(S_N)\otimes \Q \quad = \quad
H^0(S_N(\C), {\cal L}_{S^{w-1}V_1})^+\qquad  w >  1
\end{equation}
where ${\cal L}_{S^{w-1}V_1}$ is the local 
system on  $S_N(\C)$ corresponding to the $GL_1$-module $S^{w-1}V_1$
 and $+$ means invariants of the  involution acting on $S_N(\C)$ as the complex conjugation $c$, and on the local system via  $v \lms -v$ on $V_1$. 
 
In particular for $N=1$ we get 
  \begin{equation} \label{demfed} 
{\rm dim}{\rm Gr}{\cal G}^{(l)}_{-w, -1}  = {\rm dim}K_{2w-1}(\Z) \otimes \Q = \quad   \left\{ \begin{array}{ll}
0 &  w \quad \mbox{even}, \quad \mbox{or $w=1$}  \\ 
 1 &    w>1 \quad \mbox{odd} \end{array} \right.
 \end{equation}
This has been known thanks to Soul\'e [So], Deligne [D] and Ihara [Ih2]. 

If we treat the pair 
 $\{S_n(\C), c\}$ as a stack and ${\cal L}_{S^{w-1}V_1}$ 
as a local system on this stack then 
$(\ref{7.8.00.3}) = H^0(\{S_n(\C), c\},  {\cal L}_{S^{w-1}V_1})$. 

Let ${\Bbb A}_{\Q}$ be the adels of $\Q$. 
In general the depth $m$, weight $w$ part of 
the standard cochain complex of the Lie algebra 
${\rm Gr}{\cal G}^{(l)}_{\bullet \bullet}(\mu_N)$ is related to the 
modular stack
\begin{equation} \label{7.8.00.2}
Y_1(m;N):= \quad GL_m(\Q) \backslash GL_m({\Bbb A}_{\Q})/ K_1(m;N) 
\cdot \R_+^* \cdot O_m
 \end{equation}
It is given by the quotient 
of 
\begin{equation} \label{7.8.00.1} 
GL_m(\Q) \backslash GL_m({\Bbb A}_{\Q})/ K_1(m;N) 
 \cdot \R_+^* \cdot SO_m
\end{equation}
by the group $\Z/2\Z = O(m)/SO(m)$ acting on the right. Here the subgroup 
$K_1(m;N) \subset \prod_p GL_m(\Z_p)$ is defined 
by imposing congruence conditions at 
the primes $p|N$. Namely, if $N = \prod p^{v_p(N)}$ then its $p$-component consists of the elements of $GL_m(\Z_p)$ whose last row is congruent to 
$(0, ..., 0, 1)$ mod $p^{v_p(N)}$. 

{\bf Example}. When $m=1$ the set 
(\ref{7.8.00.1}) is  isomorphic to $S_N(\C)$ and 
$O(1)$ acts as the complex involution $c$.  For $m>1$ the 
stack (\ref{7.8.00.2}) coincides with the modular variety (\ref{7.19.00.1}). 
When $m=2$ one has 
$GL_2(\R)/R_+^* \cdot SO(2) = \C - \R$ and $O(2)/SO(2)$ 
acts as the complex conjugation $z \lms \overline z$. 

Let us go 
beyond the depth $1$ case. We start from the $N=1$ case. 
The following result generalizes the first line of (\ref{demfed}):

\begin{theorem} \label{8.17.00.1}
${\rm dim}{\cal G}^{(l)}_{-w, -m} =0$ if $w+m$ is odd. 
\end{theorem}

{\bf 3. The Galois action on $\pi^{(l)}_{1}({\Bbb P}^{1} - \{0, 1,\infty \}, v_{\infty}  )$: a description of the depth $2$ quotient via 
 the modular triangulation of the hyperbolic plane}. 
The structure of the Lie algebra  
${\cal G}^{(l)}_{\bullet, \geq -2}$ is completely 
described by the Lie commutator map
\begin{equation} \label{2/04/00.3}
[,]:  \Lambda^2 {\cal G}^{(l)}_{\bullet, -1}
 \lra  
{\cal G}^{(l)}_{\bullet, -2}    
 \end{equation}
The left hand side is known to us by (\ref{demfed}). 
So to describe (\ref{2/04/00.3}) one needs 
to define the right hand side and the commutator map. 
Dualizing (\ref{2/04/00.3}) we get the  depth two piece of 
the standard cochain complex of the Lie algebra 
${\cal G}^{(l)}_{\bullet, \geq -2}$: 
\begin{equation} \label{dem} 
{{\cal G}^{(l)}_{\bullet, -2}}^{\vee} 
\stackrel{\delta}{\lra} {\Lambda^2 {\cal G}^{(l)}_{\bullet, -1}}^{\vee} 
\end{equation} 
To describe it we need to introduce the following two complexes:
\begin{equation} \label{4-12.1}
{M}_{(2)}^{\ast}:= {M}_{(2)}^{1} \lra {M}_{(2)}^{2}\quad\quad \mbox{and} \quad\quad  
{\Bbb M}_{(2)}^{\ast}:= {M}_{(2)}^{1} \lra {M}_{(2)}^{2} \lra 
{M}_{(2)}^{3}
 \end{equation} 
The first one is the chain complex  
 of the  
classical modular triangulation 
\begin{center}
\hspace{4.0cm}
\epsffile{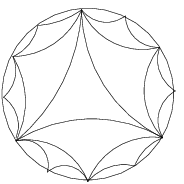}
\end{center} 
 of the hyperbolic plane ${\Bbb H}_2$ where
 the central ideal triangle has vertices at $0, 1$ and $\infty$. We place it in degrees $[1,2]$.  For example  
${M}_{(2)}^{1}$ is the group generated by the triangles. The second complex is the chain complex of the modular triangulation of the hyperbolic plane {\it extended by cusps}, i.e. by ${\Bbb P}^1(\Q)$. 
We place it in degrees $[1,3]$. 

The group  $GL_2(\R)$, acting on $\C - \R$ by  
$\left (\matrix{a&b\cr c&d\cr}\right )z = \frac{az+b}{cz+d}$, commutes with $z \lms \overline z$. We let $GL_2(\R)$ act on the upper half plane ${\Bbb H}_2$ by identifying ${\Bbb H}_2$ 
with the quotient of $\C - \R$ by complex conjugation. 
The action of the subgroup $GL_2(\Z)$ preserves the modular picture. So  
 ${M}_{(2)}^{\ast}$ and ${\Bbb M}_{(2)}^{\ast}$ are  complexes of $GL_2(\Z)$-modules.

{\it A digression on cohomology of a subgroup $\Gamma$ of $GL_2(\Z)$}. Let $V$ be a   
$GL_2$-module. For a torsion free subgroup $\Gamma$ of $GL_2(\Z)$ the group cohomology 
$H^*(\Gamma, V)$ are isomorphic to the cohomology of $\Gamma \backslash {\Bbb H}_2$ with coefficients in the local system ${\cal L}_V$ corresponding to $V$. Notice that 
$H^*(\Gamma \backslash {\Bbb H}_2, {\cal L}_V ) = H^*(\overline {\Gamma \backslash {\Bbb H}_2}, 
Rj_*{\cal L}_V )$ where 
$j: \Gamma \backslash {\Bbb H}_2 \hookrightarrow \overline {\Gamma \backslash {\Bbb H}_2}$.  

For a torsion free finite index subgroup 
$\Gamma$ of $GL_2(\Z)$
one defines the cuspidal cohomology $H^*_{{\rm cusp}}(\Gamma, V)$ as the cohomology 
of $\overline {\Gamma \backslash {\Bbb H}_2}$ with coefficients in a middle extension of the local system 
${\cal L}_V$. In our case the middle extension means 
the sheaf $j_*{\cal L}_V$.

For any finite index subgroup $\Gamma \hookrightarrow GL_2(\Z)$ there is a normal 
torsion free finite index subgroup $\widetilde \Gamma \hookrightarrow  \Gamma$.
So if $V$ is a $\Q$-rational $GL_2$-module one has, using the 
Hochshild-Serre spectral sequence,  
\begin{equation} \label{4-11.1d}
H^*(\Gamma, V) \quad = \quad H^*(\widetilde \Gamma, V)^{\Gamma/\widetilde \Gamma} 
\end{equation}
We define the cuspidal cohomology for arbitrary finite index subgroup $\Gamma$ 
by reducing it to the  torsion free case by formula similar to  (\ref{4-11.1d}).  

\begin{lemma} \label{modcoho}
Let $\Gamma$ be a finite index subgroup of $GL_2(\Z)$ and $V$ a $\Q$-rational  
$GL_2$-module. Then the complex 
 \begin{equation} \label{4-11.1}
{M}_{(2)}^{\ast} \otimes_{\Gamma} V [1]
\end{equation}
computes the cohomology $H^{*}(\Gamma,  V\otimes \varepsilon_2)$, and the  complex
\begin{equation} \label{4-11.1q}
{\Bbb M}_{(2)}^{\ast} \otimes_{\Gamma} V [1]
 \end{equation}
computes the cuspidal cohomology $H^*_{{\rm cusp}}(\Gamma,  V\otimes \varepsilon_2)$. 
\end{lemma}

{\bf Proof}. The complex (\ref{4-11.1}), up to a shift,  is the complex of chains with 
coefficients in the local system ${\cal L}_V$; it is relative to the modular triangulation of 
$\overline {\Gamma \backslash {\Bbb H}_2}$, which has finitely many cells, with ${\Gamma \backslash {\Bbb H}_2}$
a union of open cells, hence giving homology groups for locally finite chains 
(the Borel-Moore homology).

If we take the dual for the complex (\ref{4-11.1}), as well as for (\ref{4-11.1q}), 
we obtain a cochain complex computing 
the cohomology of $\overline {\Gamma \backslash X}$ with coefficients in a sheaf: 
the sheaf is 
$j_!{\cal L}_{V^{\vee}}$ for (\ref{4-11.1}) and the middle extension 
$j_*{\cal L}_{V^{\vee}}$ for (\ref{4-11.1q}) 
(notice that invariants by the stabilizer of the cusp dual to the coinvariant group). 

For the complex itself by Poincar\'e duality we get cohomology of 
${\Gamma \backslash {\Bbb H}_2}$ 
(respectively $\overline {\Gamma \backslash {\Bbb H}_2}$) with value in the local 
system ${\cal L}_{V}$ twisted by the orientation class, i.e. ${\cal L}_{V \otimes 
\varepsilon_2}$ (respectively its middle extension). The lemma is proved.

To clarify the homological algebra meaning of the modular complex consider complex 
 ${\cal M}^{(2)}_{-1} \stackrel{\partial}{\longleftarrow} {\cal M}^{(2)}_{-2}$,
  where ${\cal M}^{(2)}_{-*}:= {M}_{(2)}^{*}$ and $\partial$ is dual to 
the differential in the complex ${M}_{(2)}^{*}$. Then ${M}_{(2)}^{*}$ is a subcomplex of  
${\rm Hom}({\cal M}^{(2)}_{-*}, \Z)$. As a complex of $SL_2(\Z)$-modules ${\cal M}_{(2)}^{*}$ is 
the chain complex of the tree dual to the modular triangulation, shifted by $1$. 
As a complex of $GL_2(\Z)$-modules it is a resolution of $\varepsilon_2[1]$. 
The stabilizers of the action of $GL_2(\Z)$  on the sets of triangles and edges 
of the modular triangulation are finite. 
Thus for any $\Q$-rational $GL_2$-module $V$ and any subgroup $\Gamma \subset GL_2(\Z)$
$$
{\rm Hom}_{\Gamma}\Bigl({\cal M}^{(2)}_{-*}, V\Bigr)\quad \mbox{computes} \quad H^{*-1}(\Gamma, V)
$$
%The group $GL_2(\Z)$ acts transitively on the sets of triangles and edges 
%of the modular triangulation. So if 
%$\Gamma$ is a finite index subgroup of $GL_2(\Z)$ then the canonical map
%$
%{M}_{(2)}^{*} \otimes_{\Gamma}V \lra {\rm Hom}_{\Gamma}({\cal M}^{(2)}_{-*}, V)
%$ 
%is an isomorphism. 

Let us return to description of the image of the Galois group.

Let ${\cal G}_{\bullet \bullet}$ be any bigraded Lie algebra. Consider the bigraded Lie algebra
\begin{equation} \label{4-12.3}
\widehat {\cal G}_{\bullet \bullet} := \quad 
{\cal G}_{\bullet \bullet} \oplus \Q(-1,-1) 
\end{equation}
where $\Q(-1,-1)$ is a one dimensional Lie algebra 
 of the bidegree $(-1, -1)$. The standard cochain complex of 
$\widehat {\cal G}_{\bullet \bullet}$ admits a canonical decomposition 
\begin{equation} \label{4-16.10+}
\Lambda^* \widehat {\cal G}^{\vee}_{\bullet \bullet} \qquad =\qquad 
\Lambda^* {\cal G}^{\vee}_{\bullet \bullet}  \quad \oplus \quad
\Lambda^* {\cal G}^{\vee}_{\bullet \bullet} \otimes \Q(1,1)
\end{equation}

Strangely enough it is simpler to describe the structure of the Lie algebra 
${\widehat {\cal G}}^{(l)}_{\bullet, \geq -2}$ than ${\cal G}^{(l)}_{\bullet, \geq -2}$.   
The  Lie algebra structure of $\widehat {\cal G}^{(l)}_{\bullet, \geq -2}$ 
is completely described by the commutator map
\begin{equation} \label{7.7.00.2}
[ , ]: {\Lambda^2 {\widehat {\cal G}}^{(l)}_{\bullet, -1}} \lra {  {\cal G}^{(l)}_{\bullet, -2}}
\end{equation} 
\begin{theorem} \label{4-13.1} a) The weight $w$ part of the complex 
\begin{equation} \label{dem1}
{  {\cal G}^{(l)}_{\bullet, -2}}^{\vee}   
%={   {\widehat {\cal L}}^{(l)}_{\bullet, -2}}^{\vee}   
\quad \stackrel{\delta}{\lra} \quad 
{\Lambda^2 {\widehat {\cal G}}^{(l)}_{\bullet, -1}}^{\vee}
\end{equation} 
dual to  (\ref{7.7.00.2}), 
is isomorphic to the complex  
\begin{equation} \label{dem2}
\Bigl({M}_{(2)}^{\ast} \otimes_{GL_2(\Z)} S^{w -2}V_2 \Bigr) \otimes \Q_l
 \end{equation}
b)\begin{equation} \label{demfedd} 
 \qquad {\rm dim} {\rm Gr}{\cal G}^{(l)}_{-w, -2}  = \quad    \left\{ \begin{array}{ll}
0 &  w: \quad \mbox{odd}  \\ 
   \left[ \frac{w-2}{6} \right]  &    w:  \quad \mbox{even} \end{array} \right.
 \end{equation}
\end{theorem} 

Using the part a) and lemma \ref{modcoho} we compute the Euler characteristic of the complex 
(\ref{dem1}). Then we use (\ref{demfed}) to get formula (\ref{demfedd}). 
A Hodge-theoretic version of theorem \ref{4-13.1} appeared first 
as theorem 7.2 in [G1]. 
The estimate ${\rm dim} {\cal G}^{(l)}_{-w, -2} \geq [\frac{w-2}{6}]$  
has been independently obtained by 
Ihara and Takao ([Ih3]).

In particular there is a {\it canonical} isomorphism
\begin{equation} \label{4-12.5}
m_1^{(l)}: {{\cal G}^{(l)}_{-w, -2}}^{\vee}  \stackrel{=}{\lra} 
\Bigl({M}_{(2)}^{1} \otimes_{GL_2(\Z)} S^{w -2}V_2 \Bigr)\otimes \Q_l
\end{equation}
providing a description of the vector space ${{\cal G}^{(l)}_{-w, -2}}$. 
To describe the isomorphism of complexes from theorem \ref{4-13.1} we need to define the map 
\begin{equation} \label{odin}
m_2^{(l)}:  \Lambda_w^2\Bigl(\widehat{\cal G}^{(l)}_{\bullet, -1}\Bigr)^{\vee} \stackrel{=}{\lra}  
\Bigl({M}_{(2)}^{2} \otimes_{GL_2(\Z)} S^{w -2}V_2\Big) \otimes \Q_l
\end{equation}
where $\Lambda_w^2$ the weight $w$ part of $\Lambda^2$. We do it  
as follows. 
The stabilizer in $GL_2(\Z)$ of the geodesic from $0$ to $i\infty$ on the upper half plane 
is generated by 
$$
\left (\matrix{0&1\cr -1&0\cr}\right ) \quad \mbox{and} \quad 
\left (\matrix{-1&0\cr 0&1\cr}\right )$$ 
Therefore the right hand side of (\ref{odin}) is identified with the space  
of degree $w-2$ polynomials  $f(t_1, t_2)$, skewsymmetric in the variables  $t_1, t_2$ and of 
even degree in each of  them: $f(t_1, t_2) = -f(t_2, -t_1) = f(-t_1, t_2)$. 

Let $\zeta^{(l)}_{{\cal M}}(n)$ for $n>1$ be Soul\'e's generator 
of $\Bigl({\rm Gr}{\cal G}^{(l)}_{-2n+1, -1}\Bigr)^{\vee}$, see section 6. We define 
$\zeta^{(l)}_{{\cal M}}(1)$ as a generator of $\Q_l(-1,-1) = \Bigl({\rm Gr}{\cal G}^{(l)}_{-1, -1}\Bigr)^{\vee}$ and set
$$
m_2^{(l)}: \zeta^{(l)}_{{\cal M}}(2m+1) \wedge \zeta^{(l)}_{{\cal M}}(2n+1) \lms 
t_1^{2m}t_2^{2n} - t_1^{2n}t_2^{2m} 
$$
The  map $m_2^{(l)}$ identifies 
$\Lambda^2\Bigl({\cal G}^{(l)}_{\bullet, -1}\Bigr)^{\vee}$ with the subspace 
of polynomials  $f(t_1, t_2)$, 
skewsymmetric in $t_1, t_2$ and of 
{\it positive} even degree in each of them. Thus we get a precise description of the 
Lie algebra ${\cal G}^{(l)}_{\bullet, \geq -2}$ as well.

The decomposition (\ref{4-16.10+}) of complex (\ref{dem1})  
corresponds to decomposition of complex 
(\ref{dem2}) into a direct sum of the 
subcomplex computing the (truncated) cuspidal cohomology, and the subcomplex, 
consisting of a single group in the degree $2$, computing 
the Eisenstein part of the cohomology.

More precisely, consider the truncated complex
$$
\tau_{[1,2]}\Bigl({\Bbb M}_{(2)}^{\ast} \otimes_{GL_2(\Z)} V \Bigr) 
\quad :=  \quad 
$$
\begin{equation} \label{mc3f}
{M}_{(2)}^{1} \otimes_{GL_2(\Z)} V \lra 
{\rm Ker}\Bigl({M}_{(2)}^{2} \otimes_{GL_2(\Z)} V \lra {M}_{(2)}^{3} \otimes_{GL_2(\Z)} V \Bigr)
\end{equation}
 It is a subcomplex of the complex (\ref{4-11.1}). 

\begin{theorem} \label{mth1}
The complex  
\begin{equation} \label{dem3}
\tau_{[1,2]}\Bigl({\Bbb M}_{(2)}^{\ast} \otimes_{GL_2(\Z)} S^{w -2}V_2 \Bigr) \otimes \Q_l
\end{equation}
is canonically isomorphic to the  weight $w$ part of the complex (\ref{dem}).
\end{theorem}

{\bf 4. The Galois action on $\pi^{(l)}_{1}({\Bbb P}^{1} - \{0, 1,\infty \}, v_{\infty}  )$: the depth $3$ quotient}. To describe the structure of the Lie algebra 
${\cal G}^{(l)}_{\bullet, \geq -3}$ we need only to define the  commutator map 
\begin{equation} \label{3/04/00.1}
[,]: {\cal G}^{(l)}_{\bullet, -2} \otimes {\cal G}^{(l)}_{\bullet, -1} 
\lra {\cal G}^{(l)}_{\bullet, -3} 
\end{equation}
obeying the Jacoby identity. Indeed, the commutator $
[,]: \Lambda^2{\cal G}^{(l)}_{\bullet, -1} 
\lra {\cal G}^{(l)}_{\bullet, -2} 
$ has been described in the part a) of  theorem \ref{4-13.1}. 
Dualizing one sees that we need to describe the 
depth three part of the standard cochain complex of 
${\cal G}^{(l)}_{\bullet, \geq -3}$: 
\begin{equation} \label{demfeddu} 
{{\cal G}^{(l)}_{\bullet, -3}}^{\vee} \lra {{{\cal G}^{(l)}_{\bullet, -2}}^{\vee} \otimes 
{\cal G}^{(l)}_{\bullet, -1}}^{\vee} \lra 
{\Lambda^3{\cal G}^{(l)}_{\bullet, -1}}^{\vee} 
\end{equation}
The first map is dual to the map (\ref{3/04/00.1}). 
Complex (\ref{demfeddu}) is described in s. 1.7. 
Here is an important corollary:
 
\begin{theorem} \label{mth1t}  a) The complex (\ref{demfeddu}) 
computes the cuspidal cohomology 
$$
H^i_{{\rm cusp}}(GL_3(\Z), S^{w-3}V_3)\quad \mbox{at} \quad i = 1,2,3
$$ 

b) The complex (\ref{demfeddu}) is acyclic. 

c)
\begin{equation} \label{demfeddd} 
 \qquad {\rm dim} {\rm Gr}{\cal G}^{(l)}_{-w, -3} =  
\quad    \left\{ \begin{array}{ll}
0 &  w: \quad \mbox{even}  \\ 
   \left[ \frac{(w-3)^2-1}{48} \right]  &    w:  \quad \mbox{odd} \end{array} \right.
 \end{equation} 
\end{theorem}  

By b) the Euler characteristic of the complex (\ref{demfeddu}) is zero. So using (\ref{demfed}) 
and theorem \ref{4-13.1}b) we get formula (\ref{demfeddd}). 

{\bf Remark}. Compare theorems \ref{4-13.1}, \ref{mth1t} with theorems 1.4, 1.5 in [G3]
where the dimension of the $\Q$-space of reduced multiple $\zeta$-values of depth $2$ and $3$ is estimated from {\it above} by (\ref{demfedd}) and  (\ref{demfeddd}).  
According to some standard conjectures in arithmetic algebraic geometry 
one should have
$$
{\cal G}^{(l)}_{-w, -m} = (\mbox{the space of reduced multiple $\zeta$'s of weight $-w$, depth $-m$})\otimes_{\Q}\Q_l
$$
This supplies an additional  evidence for  the 
statement that the estimates given in theorems 1.4, 1.5 in [G3] are exact. 
Notice that computation of the dimension of the $\Q$-space of reduced multiple $\zeta$-values 
seems to be a transcendently  difficult problem 
(we can not prove that $\zeta(5) \not \in \Q$!),
while its more 
sophisticated $l$-adic analog is easier to approach. 

To describe results and conjectures about 
the structure of the Lie algebra ${\rm Gr}{\cal G}^{(l)}_{\bullet \bullet}$ we need 
to recall 
the definition of the modular complexes given in the section 5 of [G3].

{\bf 5`. The modular complexes}.  Let $L_m$ be a rank $m$ lattice. 
The modular complex $
{M}_{(m)}^{*}$ is a complex of left ${\rm Aut}(L_m) = GL_m(\Z)$-modules 
which  sits in the degrees $[1,m]$ and looks as follows:
$$
M_{(m)}^{1} \stackrel{\partial}{\lra} M_{(m)}^{2} \stackrel{\partial}{\lra}
... \stackrel{\partial}{\lra} M_{(m)}^{m}
$$
By definition $M_{(1)}^1$ is the trivial $GL_1(\Z)$-module. 
 When $m=2$ it is isomorphic to the complex on the left in (\ref{4-12.1}). 
In general its definition is purely combinatorial. We recall it now.

i) {\it Extended basis and homogeneous affine basis of a lattice}. 
We say that an  {\it extended basis} of a lattice $L_m$ is an  $(m+1)$-tuple of vectors 
$v_1, ..., v_{m+1}$ 
of the lattice such that  $v_1 + ... + v_{m+1} = 0$ and $v_1,...,v_m$ is a basis. 
Then omitting any of the vectors $v_1, ..., v_{m+1}$ we get a basis of the lattice. 

The group $GL_m(\Z)$ acts from the left on the set of basis of $L_m$, considered as columns
of vectors $(v_1, ..., v_m)$. This provides the set of extended basis with the 
structure of the left principal homogeneous space  for $GL_m(\Z)$.

Let $u_1, ..., u_{m+1}$ be elements of the lattice $L_m$ such that the set of elements 
$\{(u_i, 1)\}$ form a basis of  $L_m \oplus \Z$. The lattice $L_m$ acts on 
such sets by 
$
l: \{(u_i, 1)\} \lms \{(u_i+l, 1)\}
$.  
We call the coinvariants of this action 
{\it homogeneous affine basis} of $L_m$ and denote them by $\{u_1: ... : u_{m+1}\}$. 
Notice that 
 $\{u_1: ... : u_{m+1}\}$ is a homogeneous affine basis if and only if 
 $\{u_2-u_1, u_3-u_2, ..., u_1-u_{m+1}\}$ is an extended basis. 
So there is a canonical bijection 
\begin{equation} \label{5.26.1}
\mbox{homogeneous affine basis of $L_m$} \quad <-> \quad 
\mbox{extended basis of $L_m$}
\end{equation} 
$$
\{u_1: ... : u_{m+1}\} \quad <-> \quad  \{u_2-u_1, u_3-u_2, ..., u_1-u_{m+1}\}
$$

ii) {\it The group $M_{(m)}^{1}$}. 
The abelian group $M_{(m)}^{1}$ is generated  by the elements $<v_1, ... , v_{m+1}>$ corresponding to 
extended  basis $v_1, ... , v_{m+1}$ of $L_m$. 
To list the relations we need another set of the generators corresponding to the 
homogeneous affine basis of $L_m$ via (\ref{5.26.1}):
$$
<u_1: ... :u_{m+1}>:= \quad <u'_1, u'_2, ..., u'_{m+1}>, \quad u_i':= u_{i+1}-u_i
$$
Let $\Sigma_{p,q}$ is the set of all shuffles of the ordered sets 
$\{1, ..., p\}$ and $\{p+1, ..., 
p+q\}$.

{\bf Relations}.   For $m=1$ we have $<v_1, v_2> = <v_2, v_1>$. 

For any $1 \leq k \leq m$ one has:
\begin{equation} \label{sshh1}
\sum_{\sigma \in \Sigma_{k,m-k}} 
<v_{\sigma(1)}, ... ,v_{\sigma(m)}, v_{m+1}> \quad = \quad 0
\end{equation}
 \begin{equation} \label{sshh3}
\sum_{\sigma \in \Sigma_{k,m-k}} <u_{\sigma(1)}: ... :u_{\sigma(m)}: u_{m+1}> \quad = \quad 0
\end{equation} 

\begin{theorem} \label{DAHS} The double shuffle relations (\ref{sshh1}), (\ref{sshh3}) 
imply the following dihedral symmetry relations for $m \geq 2$:
$$
<v_{1}, ... ,v_{m}, v_{m+1}> = <v_{2}, ... ,v_{m+1}, v_{1}> =  
$$
$$
<-v_{1}, ... ,-v_{m}, -v_{m+1}> = (-1)^{m+1}<v_{m+1}, v_{m}, ... , v_{1}> 
$$
\end{theorem}

We picture both types of the generators on the circle as shown on the picture:
\begin{center}
\hspace{4.0cm}
\epsffile{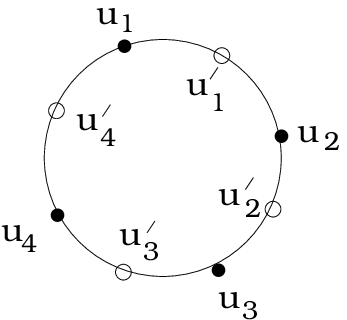}
\end{center}
Vectors of homogeneous affine basis are outside, 
and vectors of extended basis are inside of the circle (compare with s. 4.1 below). 

iii) {\it The group $M_{(m)}^{k}$}. For each decomposition of $L_m$ as a direct 
sum of $k$ non zero lattices $L^i$, we consider the tensor product of the 
$M^{1}(L^i)$. We consider 
 $M^{1}(L^i)$ as odd, and use the sign rule to identify $\otimes M^{1}(L^i)$ and 
$\otimes M^{1}(L^{i'})$ when the decomposition $L$ and $L'$ differ only in the 
ordering of the factors. The group $M_{(m)}^{k}$ is defined to be the sum over 
all such unordered 
decomposition $L_m=\oplus L^i$ of the corresponding 
$\otimes M^{1}(L^i)$.  
In other words it is generated by the elements $<A_1> \wedge ... \wedge <A_k>$ where 
$A_i$ is an extended  basis of the sublattice $L_i$ and $<A_i>$'s anticommute. 
Let us set 
$$
[v_1,...,v_k]:= <v_1,...,v_k, v_{k+1}>, \quad v_1 + ... + v_k + v_{k+1} = 0
$$
 We define a homomorphism  
$\partial: M_{(m)}^{1} \lra M_{(m)}^{2}
$  
by setting $\partial = 0$ if $m=1$ and 
$$
\partial:\quad <v_1,...,v_{m+1}> \quad \lms \quad -{\rm Cycle}_{m+1}\Bigl(\sum_{k=1}^{m-1} 
[v_1,...,v_k] \wedge [v_{k+1},...,v_m] \Bigr)
$$
where  the indices are modulo $m+1$. We get the differential in  $M_{(m)}^{\bullet}$ 
by extending   $\partial$ 
using  the Leibniz rule:
$$
\partial( [A_1] \wedge [A_2] \wedge ... ) := \quad 
\partial([A_1]) \wedge [A_2] \wedge ... \quad  -  \quad  [A_1] \wedge \partial([A_2]) \wedge ...  \quad + \ldots
$$

\begin{theorem} \label{6.15.00.1}
The map $\partial$ is  well defined and $\partial^2=0$, so we get a complex. 
\end{theorem}

{\bf Remark}. 
Both the definition of the modular complex and these proofs are very similar to the definitions 
and basic properties of the dihedral Lie algebras discussed in detail in the 
section 4 below. 
Yet I do not know a general framework unifying them. 

{\bf Remark}. Let $A$ be an 
arbitrary commutative ring. Then replacing 
$L_m$ by a free rank $m$ $A$-module one can construct modular complexes 
corresponding to  the ring $A$, recovering the 
construction above when $A = \Z$. 

{\bf 6. A description of the Lie algebra ${\rm Gr}\widehat {\cal G}^{(l)}_{\bullet \bullet}
$ via the modular complex}.
Set $V_m:= L_m\otimes \Q$. 

\begin{conjecture} \label{dwhn1*} 
There exists a canonical isomorphism between the complex 
\begin{equation} \label{dep5}
\Bigl(M^*_{(m)} \otimes_{GL_m(\Z)} S^{w-m}{\rm V}_m\Bigr) \otimes_{\Q}\Q_l
\end{equation}
and the depth $m$, weight $w$ part of the standard cochain complex of the Lie algebra ${\rm Gr}{\widehat {\cal G}}_{\bullet \bullet}^{(l)}$. 
\end{conjecture}

{\bf Example}. The rank $1$ modular complex $M^1_{(1)}$ is 
the trivial $GL_1(\Z)$-module $\Z$. Thus 
$$
M^1_{(1)}\otimes_{GL_1(\Z)} S^{w-1}V_1 = \quad   \left\{ \begin{array}{ll}
0 &  w \quad \mbox{even}  \\ 
 \Q &    w \quad \mbox{odd} \end{array} \right.
$$
Thus formula (\ref{dep5}) for $m=1$ is just equivalent to (\ref{demfed}) plus 
${\rm dim} {\rm Gr}\widehat {\cal G}^{(l)}_{-1, -1} = 1$. 
For $m=2$ we get the 
description of the depth two part given in theorem \ref{4-13.1}.

\begin{theorem} \label{dwhn1} 
There exists a canonical isomorphism between the complex 
\begin{equation} \label{dep55}
\Bigl(M^*_{(3)} \otimes_{GL_3(\Z)} S^{w-3}{\rm V}_3\Bigr) \otimes_{\Q}\Q_l
\end{equation}
and the depth $3$, weight $w$ part of the standard cochain complex of the Lie algebra 
$\widehat {\cal G}_{\bullet, \geq -3}^{(l)}$. 
\end{theorem}

\begin{theorem} \label{dwhn2} 
{\bf [G3]} There rank $3$ modular complex is quasiisomorphic to the chain complex of 
the Voronoi decomposition 
of the symmetric space ${\Bbb H}_3$ for $GL_3(\R)$, truncated in the degrees $[1,3]$ . 
\end{theorem}
This is theorem 6.2 in [G3]. To prove it we constructed
 a geometric realization of the rank $3$ modular complex as 
 a subcomplex of the truncated Voronoi complex.

On the geometric realization of the 
rank $m$ modular complex see in [G5-7]. Using it  
we proved that the rank $4$ modular complex, shifted by $[-2]$, 
 is quasiisomorphic to chain complex of the Voronoi decomposition  of the 
symmetric space ${\Bbb H}_4$, truncated in the degrees 
$[3,6]$.

{\bf 7. The diagonal Galois Lie algebras  
and  modular complexes}.  
The diagonal Lie algebra ${\rm Gr}{\cal G}^{(l)}_{\bullet}(\mu_N)$ is a Lie subalgebra of the Galois Lie algebra ${\rm Gr}{\cal G}^{(l)}_{\bullet \bullet}(\mu_N)$:
$$
{\rm Gr}{\cal G}^{(l)}_{\bullet}(\mu_N) := \quad \oplus_{w \geq 1}{\rm Gr}{\cal G}^{(l)}_{-w, -w}(\mu_N)\quad \hookrightarrow \quad 
{\rm Gr}{\cal G}^{(l)}_{\bullet \bullet}(\mu_N)
$$ 
It is graded by the weight. 

Equivalently, one can define the diagonal Galois Lie algebra 
as a Lie  subalgebra of ${\rm Gr}^W{\cal G}^{(l)}(\mu_N)$ by 
imposing the depth $\leq -w$ condition on each weight $-w$ subquotient:
$$
{\rm Gr}{\cal G}^{(l)}_{\bullet}(\mu_N) := \quad 
\oplus_{w \geq 1}{\cal F}^D_{-w}{\rm Gr}^W_{-w}{\cal G}^{(l)}(\mu_N)\quad \hookrightarrow \quad 
{\rm Gr}^W{\cal G}^{(l)}(\mu_N)
$$ 

Since the Lie algebra ${\rm Gr}^W{\cal G}^{(l)}(\mu_N) $ is non canonically isomorphic to the Lie algebra ${\cal G}^{(l)}(\mu_N) $, the diagonal Galois Lie algebra is isomorphic, although 
non canonically,  to a certain Lie subalgebra of ${\cal G}^{(l)}(\mu_N) $.
This is an important difference between the Galois Lie algebra and its diagonal part: the Galois Lie algebra in general is not isomorphic to ${\cal G}^{(l)}(\mu_N) $.

\begin{theorem} \label{4-12.6}
There is canonical isomorphism 
$$
{\rm Gr}{{\cal G}^{(l)}_{-1, -1}}(\mu_N) \quad = \quad
{\rm Hom}_{\Q}\Bigl(\mbox{the group of the cyclotomic units in $
\Z[\zeta_N][\frac{1}{N}]$}, \quad \Q_l\Bigr)
$$ 
\end{theorem}
This rather elementary result 
 is a particular case of theorem \ref{4-12.60090} for $w=1$.
It suggested the name ``higher cyclotomy'' for our story. 
Set
$$
E_m(N):= \Z[\Gamma_1(m;N)\backslash GL_m(\Z)] \quad \mbox{if $m>1$} \quad E_1(N):= \Z/N\Z
$$
\begin{conjecture} \label{CD*} Let $p$ be a prime number. Then there exists a canonical 
isomorphism between the complex 
\begin{equation} \label{dep5*}
\Bigl(M^*_{(m)} \otimes_{GL_m(\Z)} E_m(N)\Bigr)\otimes \Q_l
\end{equation}
and the depth $m$ part of the standard cochain complex of 
${\rm Gr}{\widehat {\cal G}}^{(l)}_{\bullet}(\mu_p)$. 
\end{conjecture}
Notice that if $m>1$ then 
$$
M^*_{(m)} \otimes_{\Gamma_1(m;N)} \Q = M^*_{(m)} \otimes_{GL_m(\Z)} E_m(N)
$$

 The depth $\geq -2$ quotient of ${\rm Gr}^{(l)}_{\bullet}\widehat {\cal G}(\mu_p) $ 
is described by the commutator map 
\begin{equation} \label{2/7/00/2}
[,]: \quad \Lambda^2{\rm Gr}\widehat {\cal G}^{(l)}_{-1}(\mu_p) \lra 
{\rm Gr}{\cal G}^{(l)}_{ -2}(\mu_p)
 \end{equation}
Let $\Gamma_1(p):= \Gamma_1(2;p) \cap SL_2(\Z)$. 
Consider the  modular curve $Y_1(p):= \Gamma_1(p)\backslash {\Bbb H}_2$. 
Projecting the modular triangulation 
of the hyperbolic plane onto $Y_1(p)$ we get  the modular triangulation of 
$Y_1(p)$. The complex involution acts on the modular curve preserving the 
triangulation. Consider the  following complex
\begin{equation} \label{MAN}
\Bigr(\mbox{{\rm the 
chain complex of the modular triangulation of $Y_1(p)$} }\Bigl)^+
\end{equation}
Here $+$ means the invariants of the action of the 
complex involution. 

\begin{theorem} \label{4-15.10} The dual to complex (\ref{2/7/00/2})
is naturally {\bf isomorphic} to complex (\ref{MAN}).
\end{theorem}

This is a depth two analog of 
theorem \ref{4-12.6}. 
In particular there is canonical isomorphism 
\begin{equation} \label{4-12.15}
\Q_l[\mbox{triangles of the modular triangulation of $Y_1(p)$}]^+ = 
\Bigl( {\rm Gr}{\cal G}^{(l)}_{ -2}(\mu_p)\Bigr)^{\vee}
\end{equation}

It turns out that the subcomplex 
$$
\tau_{[1,2]}\Bigl({\Bbb M}_{(2)}^* \otimes_{\Gamma_1(2;p)} \Q_l\Bigr) \quad  
\hookrightarrow \quad  {M}_{(2)}^* \otimes_{\Gamma_1(2;p)} \Q_l
$$
corresponds to the maximal quotient of the Galois group 
acting on $\pi^{(l)}(X_p)$, of weight=depth $\geq -2$, 
 which is {\it unramified at} $1-\zeta_p$, see section 7.9.

Now we turn to the depth three case. Let ${\rm Gr}{\cal L}^{(l)}_{\bullet}(\mu_N) \hookrightarrow 
{\rm Gr}{\cal G}^{(l)}_{\bullet}(\mu_N)$ be the Lie subalgebra  generated by 
 ${\rm Gr}{\cal L}^{(l)}_{-1}(\mu_N)$.  
\begin{theorem} \label{4-15.11} a) Conjecture {\rm \ref{CD*}} is valid for $m=3$ for the Lie algebra 
${\rm Gr}\widehat {\cal L}^{(l)}_{ \bullet}(\mu_p)$. So there exists canonical isomorphism between the complex (\ref{dep5*}) and the depth $3$ 
part of the standard cochain complex of this Lie algebra. 

b) Let us assume conjecture \ref{ramierz} below. Then 
${\rm Gr}{\cal L}^{(l)}_{-3}(\mu_p) = {\rm Gr}{\cal G}^{(l)}_{-3}(\mu_p)$, and therefore conjecture  {\rm \ref{CD*}} would be valid for $m=3$. 
\end{theorem}

There is a description of the vector space 
${\rm Gr}{\cal L}^{(l)}_{ -3}(\mu_p)$ similar to (\ref{4-12.15}) 
in terms of certain 
$4$-cells on the $5$-dimensional modular variety 
$\Gamma_1(3; p)\backslash {\Bbb H}_3$ 
which were constructed in [G3]. 
See [G5-7] for a construction of certain $2(m-1)$-cells in 
${\Bbb H}_m$ which should play a similar role in general.  

\begin{corollary} \label{4-19.10} If $p \geq 5$ then 
$$
{\rm dim}{\rm Gr}{\cal G}^{(l)}_{-2}(\mu_p) \quad = \quad 
\frac{(p-1)(p-5)}{12}   
$$
$$
{\rm dim}{\rm Gr}{\cal G}^{(l)}_{-3}(\mu_p) \geq 
{\rm dim}{\rm Gr}{\cal L}^{(l)}_{-3}(\mu_p)\quad = \quad     
  \frac{(p-5)(p^2-2p - 11)}{48}   
$$
\end{corollary}

{\bf 9. The structure of the paper}. In the section $3$ we explain how to 
linearize the action of the Galois group 
${\rm Gal}(\overline \Q/\Q(\zeta_{l^{\infty}N}))$ on 
$\pi^{(l)}_{1}(X_N, v_{\infty})$. 
We show that 
the Lie algebra of the image of the Galois group 
  acts by the so-called special equivariant derivations of 
${\Bbb L}^{(l)}(X_N, v_{\infty})$, preserving the two filtrations. 

In the section 4 we define the dihedral Lie coalgebra 
${\cal D}_{\bullet \bullet}(G)$ of a commutative group $G$ and derive its main  
properties. 
Completely similar arguments settle the basic properties 
of the modular complexes.

In the section $5$ we study the Lie algebra 
${\rm Der}^{SE}{L}(G)$ of 
special equivariant derivations of the free Lie algebra 
${L}(G)$ generated by the set $\{0\} \cup G$. 
We realize the dihedral Lie algebra of $G$ as a Lie subalgebra 
of ${\rm Gr}{\rm Der}^{SE}{L}(G)$.
We prove the distribution 
relations in s. 5.7.  

The depth $2$ and $3$ parts of the cohomology of the dihedral Lie algebra of $\mu_N$ 
were related in [G3] to the cohomology of the groups $\Gamma(m; N)$ with coefficients 
in $S^{w-m}V_m$ 
for $m=2,3$. In the section $6$ 
we compute the cohomology groups
$$
H^*(GL_3(\Z), S^{w-3}V_3) \quad \mbox{and} \quad H^*(\Gamma(3; p), \Q)
$$

Finally in the chapter $7$ we apply the previous results 
to prove the theorems from the introduction. 

\section{The setup}

 {\bf 1. The action of the Galois group}. 
Recall the   map
\begin{equation} \label{121212*}
\Phi^{(l)}_{N}: \quad {\rm Gal}(\overline{\Bbb Q}/ {\Bbb Q}) \longrightarrow 
{\rm Aut} \pi^{(l)}_{1}(X_N, v_{\infty} )
\end{equation} 
 and the projection
\begin{equation} \label{PP}
 \pi^{(l)}_{1}(X_N, v_{\infty} )\lra  
 \pi^{(l)}_{1}({\Bbb G}_m, v_{\infty} )   = \Z_l(1) 
\end{equation}
Let $X$ be a regular curve over $\overline \Q$,   
$\overline X$ the corresponding projective curve and $v$ a tangent vector at 
$x \in \overline X$. Then  there is a natural map of Galois modules ([D]):
$$
\Z_l(1) = \pi_1^{(l)}(T_x \overline X - 0, v) \to \pi_1^{(l)}(X, v)
$$ 
For 
$X = {\Bbb P}^1 - \{0, \infty\}$, $x = \infty, v = v_{\infty}$ 
it is an isomorphism. So for $X = X_N$ it 
provides a 
  splitting  of  (\ref{PP}): 
\begin{equation} \label{III}
I_{\infty}: \pi^{(l)}_{1}({\Bbb G}_m, v_{\infty} ) \hookrightarrow 
\pi^{(l)}_{1}(X_N, v_{\infty})
\end{equation}
The 
 subgroup ${\rm Gal}(\overline \Q/\Q(\zeta_{l^{\infty}}))$ preserves the elements of the subgroup 
\begin{equation} \label{II}
I_{\infty}(\Z_l(1)) \subset \pi_1^{(l)}(X_N, v_{\infty})
\end{equation}

Let $v$, $v'$ be tangent vectors at $x, x' \in \overline X$. Denote by 
$\pi^{(l)}_{1}(\overline X; v, v')$ the pro-$l$ completion of 
the torsor of path from $v$ to $v'$. 
Let $\eta \in \{0\} \cup \mu_N \subset {\Bbb P}^1$. 
Choose  a ${\rm Gal}(\overline \Q/\Q(\zeta_N))$-invariant tangent vector $w_{\eta}$ at the point $\eta$. 
Let $p$ be a pro-$l$ path from $v_{\infty}$ to $w_{\eta}$. The composition  
$$
p^{-1} \circ \pi_1^{(l)}({\Bbb P}^1 - \{0, \eta\}, w_{\eta}) \circ p
$$ 
provides a map
\begin{equation} \label{IIP}
I_{\eta}(p): \Z_l(1) \lra \pi_1^{(l)}(X_N, v_{\infty}), \quad 
\eta  \in \{0\} \cup \mu_N 
\end{equation}
For a different path $p'$ the map $I_{\eta}(p')$ is conjugated to $I_{\eta}(p)$ in $\pi_1^{(l)}$. So the  conjugacy class of this map is well defined, hence stable by ${\rm Gal}(\overline \Q/\Q(\zeta_N))$.

The restriction of the map $\Phi_N^{(l)}$ to 
${\rm Gal}(\overline \Q/\Q(\zeta_{l^{\infty}N}))$ satisfies additional constraints. 
The group $\mu_N$ acts 
on  $X_N$ by $z \lms \zeta_N z $. Moreover, 
there is a natural 
action of 
$\mu_N$ on $\pi^{(l)}_{1}(X_N, v_{\infty}) \otimes \Q_l$  
 (regarding   $\pi^{(l)}_{1}\otimes \Q_l $ see [D], ch. 9) commuting with the action of ${\rm Gal}(\overline \Q/\Q(\zeta_{l^{\infty}N}))$,  
coming as follows.  
Let $\xi, \zeta, \zeta' \in \mu_N$. The action of $\xi$ on ${\Bbb G}_m$ provides 
 a tangent vector  
$v_{\xi }:= \xi_* v_{\infty}$ at $\infty$ and a map  
$$
\xi_*: \pi^{(l)}_{1}({\Bbb G}_m; v_{\zeta },  v_{\zeta'}) \lra 
\pi^{(l)}_{1}({\Bbb G}_m; v_{\xi\zeta }, v_{\xi\zeta' }) 
$$
We have a canonical path 
$
p_{\zeta, \zeta'} $ in $\pi^{(l)}_{1}({\Bbb G}_m;  v_{\zeta}, v_{\zeta'}) \otimes \Q_l
$ 
such that
$$
p_{\zeta', \zeta''} \circ p_{\zeta, \zeta'} = p_{\zeta, \zeta''}, \quad \xi_*p_{\zeta, \zeta'} =  
p_{\xi\zeta, \xi\zeta'}
$$
Indeed, the composition of path provides an isomorphism of the $N$-th power of the torsor 
of path from $v_{\zeta }$ to $v_{\xi \zeta }$ on ${\Bbb G}_m$ with 
$\pi^{(l)}_{1}({\Bbb G}_m)$ which, being abelian, does not depend on the choice of 
the base point/vector. 
It is given by 
$$
com: p_1 \otimes ... \otimes p_N \lra \xi_*^{N-1} p_1 \circ ... \circ \xi_* p_{N-1} \circ p_N
$$
Then 
$$
p_{\zeta, \xi \zeta}:= p \circ  com (p ^{\otimes N})^{-1/N} \in \pi^{(l)}_{1}({\Bbb G}_m; 
v_{\zeta}, v_{\xi \zeta}) \otimes \Q_l
$$ 
Notice that if $(N,l)=1$ then $p_{\zeta, \zeta'} \in \pi^{(l)}_{1}({\Bbb G}_m; 
v_{\zeta }, v_{\zeta'})$. 
%Similarly to (\ref{III}) 
There is a natural map 
$$
\pi^{(l)}_{1}({\Bbb P}^{1} - \{0,  \infty \}; v_{\zeta}, v_{\zeta'} ) 
\quad \hookrightarrow \quad
\pi^{(l)}_{1}(X_N; v_{\zeta}, v_{\zeta'})
$$
commuting with the action of the Galois group. It provides canonical path 
$$
p_{\zeta, \zeta'} \in \pi^{(l)}_{1}(X_N;  
v_{\zeta}, v_{\zeta'}) \otimes \Q_l
$$
We define the action of an element $\xi \in \mu_N$ on 
$\alpha \in \pi^{(l)}_{1}(X_N;  v_{\infty}) \otimes \Q_l$ by 
$\xi(\alpha):= p_{1, \xi}^{-1} \circ \xi_*(\alpha) \circ p_{1, \xi}$.

 So we see that ${\rm Gal}(\overline \Q/\Q(\zeta_{ \zeta_{l^{\infty}}N}))$ 
preserves the elements of $I_{\infty}(\Z_l(1)) $, the conjugacy classes of the 
``loops around $0$ and $\zeta^a_N$'' provided by  (\ref{IIP}), and commutes with the action of 
the group $\mu_N$. Moreover, it is compatible in a natural sense with the maps 
$X_{MN} \lra X_{M}$ given by $z \lms z$, $z \lms z^N$,  
see s. 5.7. 

{\bf 2. Passing to the Lie algebras}. Let ${\Bbb L}(X_N, v_{\infty} )$ 
be the pronilpotent Lie algebra over $\Q$ 
corresponding by the Maltsev theory to the pronilpotent completion of the fundamental group 
$\pi_{1}(X_N(\C), v_{\infty})$ (see [D] ch. 9). We use a shorthand ${\Bbb L}_N$ for it.  
${\Bbb L}_N$ is isomorphic, not canonically,  to the 
 pronilpotent completion  of a free Lie algebra with  $N+1$ generators corresponding to the 
loops around $0$ and all $N$-th roots of unity.

   Let   ${\Bbb L}_{N}^{(l)}$   be the pronilpotent Lie algebra over $\Q_l$ 
corresponding to the pro-$l$ group 
$ \pi^{(l)}_{1}(X_N, v_{\infty} )$.  
Namely, set $\pi^{(l)}_{1}:= \pi^{(l)}_{1}(X_N, v_{\infty} )$. Recall that 
 $\pi^{(l)}_{1}(k)$ is the lower central series for the group $\pi^{(l)}_{1}$. Then $\pi^{(l)}_{1}/
\pi^{(l)}_{1}(k)$ is an $l$-adic Lie group, and 
$$
{\Bbb L}_{N}^{(l)}:= \lim_{\longleftarrow} {Lie} \left( \pi^{(l)}_{1}/
\pi^{(l)}_{1}(k)\right) = {\Bbb L}_N  \widehat \otimes_{\Q}\Q_l
$$

Similar to (\ref{PP}) and (\ref{III}) there is a  canonical projection 
\begin{equation} \label{PP1}
 p: {\Bbb L}({\Bbb P}^{1} - (\{0, \infty\} \cup  \mu_N), v_{\infty} )\lra  
 {\Bbb L}({\Bbb P}^{1} - \{0, \infty\} ), v_{\infty} ) 
\end{equation}
and its canonical splitting: 
\begin{equation} \label{splittin}
i_{\infty}: {\Bbb L}({\Bbb P}^{1} - \{0,  \infty \}, v_{\infty}) = {\Bbb L}(T_0{\Bbb P}^{1}, v_{\infty}) =  \Q(1) \hookrightarrow 
{\Bbb L}({\Bbb P}^{1} - (\{0,  \infty \} \cup \mu_N), v_{\infty} )
\end{equation} 
Just as in s. 1.2, there are well defined 
conjugacy classes of ``loops 
around $0$ or $\zeta \in \mu_N$'' based at  $v_{\infty}$.  Set $X_{\infty}:= i_{\infty}$. 

\begin{lemma}  \label{splitting}
There exist maps $X_0, X_{\zeta}: \Q(1) \lra {\Bbb L}_N$ which belong to the conjugacy classes  of the 
``loops around $0, \zeta$'' such that 
$
X_0 + \sum_{\zeta \in \mu_N} X_{\zeta} + X_{\infty} =0
$ and the action of $\mu_N$ permutes $X_{\zeta}$'s (i.e. $\xi_*X_{\zeta} = X_{\xi \zeta}$) 
and fixes $X_0$, $X_{\infty}$. 
\end{lemma}

We will use the following notations. Let $G$ be a commutative group written multiplicatively.
Let  ${L}(G)$ be the free Lie algebra with the generators $X_i$ where $i \in \{0\} \cup G$ (we assume  $0 \not \in G$). Set 
$X_{\infty}:= -X_0 -\sum_{g \in G}X_g$.   

Denote by ${\Bbb L}(\mu_N)$ the pronilpotent completion of the Lie algebra ${L}(\mu_N)$. It 
is isomorphic (non canonically) to ${\Bbb L}_N$. 

{\bf Proof}. We follow the argument sketched  by the referee. 
Let ${\Bbb H}$ be the proalgebraic group over $\Q$ consisting of all automorphisms 
of the Lie algebra ${\Bbb L}(\mu_N)$ which commute with the action of $\mu_N$ and preserve 
the conjugacy classes of  the generators $X_0, X_{\zeta}, X_{\infty}$. Let ${\Bbb T}$ be the ${\Bbb H}$-torsor of all 
maps 
$X_0, X_{\zeta}, X_{\infty}:\Q(1) \lra {\Bbb L}_N$ satisfying all the conditions of the lemma. Then
 ${\Bbb T}(\C)$ is nonempty. Indeed, in the De Rham realization of ${\Bbb L}_N$ (see [D]) 
$X_{\zeta}$, $X_{\infty}$ are dual to the basis $d \log (z-\zeta)$, $d \log (z)$ in 
$H^1_{DR}({\Bbb P}^{1} - (\{0,  \infty \} \cup \mu_N), \C)$ and 
$X_0:= - \sum_{\zeta} X_{\zeta} - X_{\infty}$. 
Applying the comparison theorem between the Betti and De Rham realizations we get an element 
in ${\Bbb T}(\C)$. Since ${\Bbb H}$ is a prounipotent algebraic group one has 
$H^1({\rm Gal} (\overline \Q/ \Q), {\Bbb H}(\overline \Q)) =0$, so ${\Bbb T}$ has a $\Q$-point.  
The lemma is proved. 

There is a homomorphism
$$
\varphi_N^{(l)}: {\rm Gal}(\overline \Q/\Q) \lra {\rm Aut}{\Bbb L}_N^{(l)}
$$
  Let ${\cal Z}^{\bullet}$ be the lower central series for ${\Bbb L}_N^{(l)}$. The quotient 
${\Bbb L}_N^{(l)}/{\cal Z}^{k}$ is a finite dimensional Lie algebra over $\Q_l$. So 
the image of the subgroup ${\rm Gal}(\overline \Q/\Q(\zeta_{l^{\infty}N}))$ in ${\rm Aut}({\Bbb L}_N^{(l)}/{\cal Z}^{k})$ 
is an $l$-adic Lie group. 
Denote by ${\cal G}_N^{(l)}$ 
the projective limit (over $k$) of the Lie algebras of these Lie groups. It is a pro-Lie algebra over $\Q_l$. It acts by derivations of the Lie algebra ${\Bbb L}_N^{(l)}$. 
We will describe the constraints on the derivations we get using a more general set up 
presented  below.

{\bf 3. Special equivariant derivations}. 
A derivation ${D}$ of the Lie algebra 
${L}(G)$ is called special if there are elements  $S_i \in {L}(G)$ such that 
\begin{equation} \label{g2}
{D}(X_{i}) = [S_{i}, X_{i}] \quad \mbox{for any $i \in \{0\} \cup G$},   
\quad \mbox{and} \quad {D}(X_{\infty}) = 0
\end{equation}
The special derivations of $L(G)$ form a Lie algebra,  denoted ${\rm Der}^{S}{L}(G)$. Indeed, 
if $D(X_i) = [S_i,X_i], \quad D'(X_i) = [S'_i,X_i]  $, then 
\begin{equation} \label{DDD}
[D,D'](X_i) = [S''_i, X_i], \quad \mbox{where} \quad S''_i := 
D(S'_i) - D'(S_i) +[S'_i, S_i]
\end{equation}
The group $G$ acts on the generators by $h: X_0 \lms X_0, X_g \lms X_{hg}$. So it acts by 
 automorphisms of the Lie algebra ${L}(G)$. A derivation ${D}$ of 
${L}(G)$ is called equivariant if it  commutes with 
the action of $G$. 
 Let ${\rm Der}^{SE}{L}(G)$ be the Lie algebra of all special equivariant derivations 
of the Lie algebra ${L}(G)$.

{\bf 4. The weight and depth filtration on ${\Bbb L}_N$}. There are two increasing 
filtrations by ideals on the 
Lie algebra ${\Bbb L}_N$, indexed by integers $n \leq 0$. 

{\it The  weight filtration ${\cal F}^{W}_{\bullet}$}. 
It coincides with the lower central series for  ${\Bbb L}_N$: 
$$
{\Bbb L}_N = {\cal F}^{W}_{  -1}{\Bbb L}_N; \quad {\cal F}^{W}_{ -n -1}{\Bbb L}_N := 
[{\cal F}^{W}_{-n}{\Bbb L}_N, {\Bbb L}_N]
$$

{\it The   depth filtration ${\cal F}^{D}_{\bullet}$}. Let   ${\cal I}_N$ 
be the kernel of projection (\ref{PP1}). 
Its powers   give  
 the  depth filtration:   
$$
{\cal F}^{D}_{ 0}{\Bbb L}_N = {\Bbb L}_N, \quad {\cal F}^{D}_{ -1}{\Bbb L}_N = {\cal I}_N, 
\quad {\cal F}^{D}_{ -n-1}{\Bbb L}_N = [{\cal I}_N, {\cal F}^{D}_{ -n}{\Bbb L}_N] 
$$

The weight filtration can be defined on the Lie algebra corresponding to the 
pronilpotent completion of 
the fundamental group of an arbitrary algebraic variety. The weights on 
${\Bbb L}_N$ are obtained by 
dividing by $2$ the usual weights. The weight filtration admits a splitting, i.e. it 
is defined by a grading. 

These filtrations induce  two filtrations on the Lie algebra ${\rm Der}^{SE}{\Bbb L}_N$. 
Taking the associated graded   with respect to these filtrations we get a 
Lie algebra ${\rm Gr}{\rm Der}^{SE}_{\bullet \bullet}{\Bbb L}_N$ bigraded 
by the weight $-w$ and depth $-m$.

 {\bf 5. The weight grading and depth filtration  on ${\rm Der}^{S}L(G)$.} 
The Lie algebra $L(G)$ is bigraded 
by the weight and depth. Namely, the free generators $X_0$, $X_{g}$ are bihomogeneous: 
they are of weight $-1$, $X_0$ is of depth $0$ and the $X_g$'s are of depth $-1$.

Each of the gradings induces a filtration of $L(G)$. 
The weight filtration is given by the lower central series. It goes from 
$-\infty$ to $-1$. 
Let ${\cal I}$ be the kernel of the natural  projection 
$ {L}(G) \lra \Q$ given by $X_g \lms 0, X_0 \lms 1$. Its powers provide the depth filtration on
${L}(G)$. It goes from 
$-\infty$ to $0$. 

The Lie algebra ${\rm Der}{L}(G)$ is bigraded by the weight and depth. Its Lie subalgebras 
 ${\rm Der}^S{L}(G)$ and ${\rm Der}^{SE}{L}(G)$ are compatible with the weight grading. However they are {\it not} compatible with the depth grading. Therefore they are graded by the weight, and 
{\it filtered} by the depth. A
 derivation (\ref{g2}) is of depth $-m$ if 
 each $S_j$ mod $X_j$ is of depth $-m$, i.e. there are at least $m$ $X_i$'s different 
from $X_0$ in $S_j$ 
mod $X_j$. The depth filtration is compatible with the weight grading. 
Let ${\rm Gr}{\rm Der}^{SE}_{\bullet \bullet}{L}(G)$ be the 
associated graded for the depth filtration.

One of our 
key tools is the following theorem proved in s. 5.

\begin{theorem} \label{ramie1} 
Let $G$ be a finite commutative group. Then there exists an injective morphism of bigraded 
Lie algebras  
$$
\xi_{G}: {D}_{\bullet \bullet}(G) \hookrightarrow  
{\rm Gr}{\rm Der}^{SE}_{\bullet \bullet}{L}(G)
$$
\end{theorem}

{\bf 6. The problem of describing of the map $\varphi_N^{(l)}$.} There is no canonical 
choice of the generators $X_0, X_{\zeta}$ of the Lie algebra ${\Bbb L}_N$ 
satisfying all the 
conditions of lemma 
\ref{splitting}. 
However their projections in ${\Bbb L}_N/
[{\Bbb L}_N, {\Bbb L}_N]$ are independent of the choice involved. 
So we have a {\it canonical} 
isomorphism 
\begin{equation} \label{MYF*}
i: {L}(\mu_N) \stackrel{=}{\lra} {\rm Gr}^W_{\bullet}{\Bbb L}_N
\end{equation} 
It preserves the depth filtration. 

{\bf Remark}. Let ${\Bbb L}(X, x)$ be the fundamental Lie algebra of $X$ 
based at a point/tangent vector $x$. If $y$ is another base point/vector then 
there is a torsor of  isomorphisms 
${\Bbb L}(X_N, x) \lra {\Bbb L}(X_N, y)$ defined up to a conjugation. 
However the isomorphism of the  
associated graded for the lower series filtration $L_{*}$ 
is a {\it canonical} isomorphism
$$
i_{x,y}: {\rm Gr}^{L}{\Bbb L}(X, x) \lra {\rm Gr}^{L}{\Bbb L}(X, y)
$$ 
If $X=X_N$ then the lower series filtration coincides with 
the weight filtration. Thus the graded 
Lie algebra ${\rm Gr}^W{\Bbb L}_N:= {\rm Gr}^W{\Bbb L}(X_N, v_{\infty})$ does 
not depend on the choice of the base point or 
tangent vector $v_{\infty}$. So (\ref{MYF*}) is indeed a canonical isomorphism.

It follows from s. 3.2 that ${\cal G}_N^{(l)}$ acts by special equivariant derivations of the Lie algebra 
${\Bbb L}_N^{(l)}$, i.e. 
\begin{equation} \label{MYF}
{\cal G}_N^{(l)} \subset {\rm Der}^{SE}{\Bbb L}^{(l)}_N 
\end{equation} 
The map $\varphi_N^{(l)}$ obviously respects the two filtrations. The filtrations on ${\rm Der}^{SE}{\Bbb L}_N^{(l)}$ induce two filtrations
on ${\cal G}_N^{(l)}$. Let ${\rm Gr}{\cal G}^{(l)}_{\bullet \bullet}(\mu_N) $ be the associated graded   of the Lie algebra ${\cal G}_N^{(l)}$. There is  an inclusion 
$
{\rm Gr}{\cal G}^{(l)}_{\bullet \bullet}(\mu_N) \subset 
{\rm Gr}{\rm Der}^{SE}_{\bullet \bullet}{\Bbb L}_N^{(l)}
$ 
 which, as was stressed by the referee, is provided by the fact that 
the weight filtration admits a splitting, i.e. is defined by a grading, 
and both the depth filtration and the subspace (\ref{MYF}) are compatible with this grading. 
Such a weight splitting is provided by the eigenspaces of a Frobenius 
$F_p$, $p \not | N$. Notice that $F_p$ normalizes 
the normal subgroup ${\rm Gal} (\overline \Q/\Q(\zeta_{l^{\infty}N}))$.

Let $W_{\bullet}L$ be a filtration on a vector space $L$. A splitting $\varphi: 
{\rm Gr}^W L \lra L$ of the filtration leads to an isomorphism 
$\varphi^*: {\rm End}(L) \lra {\rm End}({\rm Gr}^WL)$. The space ${\rm End}(L)$ inherits a natural filtration, while ${\rm End}({\rm Gr}^WL)$ is graded. The map $\varphi^*$ respects the corresponding filtrations.
The map
$$
{\rm Gr}\varphi^*: {\rm Gr}^W({\rm End}L)  \lra   
{\rm End}({\rm Gr}^WL)
$$ 
does not depend on the choice of the splitting. 
Applying this to the case $L = {\Bbb L}_N$ we get a {\it canonical} isomorphism 
\begin{equation} \label{PPPQ}
{\rm Gr}^W({\rm Der}^{SE}{\Bbb L}_N) \quad \stackrel{\sim}{=} \quad {\rm Der}^{SE}{L}(\mu_N)
\end{equation}
respecting the weight grading. 
Thus there is a {\it canonical} 
injective morphism
\begin{equation} \label{myLI}
{\rm Gr}{\cal G}^{(l)}_{\bullet \bullet}(\mu_N) \quad \stackrel{}{\hookrightarrow} \quad 
{\rm Gr}{\rm Der}^{SE}_{\bullet \bullet}{\Bbb L}^{(l)}_N 
\quad \stackrel{\sim}{=} \quad {\rm Gr}{\rm Der}^{SE}_{\bullet \bullet}{L}(\mu_N) \otimes \Q_l 
\end{equation}
The goal of this paper is to study  Lie subalgebra (\ref{myLI}).  
It is canonically isomorphic to the Galois Lie algebra introduced in section 1. 
Indeed, the Lie algebras of images if the following two maps coincide:
$$
{\rm Gal}(\overline \Q/\Q(\zeta_{l^{\infty}N})) 
\lra {\rm Aut}(\pi_1^{(l)}(X_N, v_{\infty})_{[-w,-m]}
$$
$$
{\rm Gal}(\overline \Q/\Q(\zeta_{l^{\infty}N})) 
\lra {\rm Aut}{\Bbb L}_N^{(l)}/({\cal F}^W_{-w-1} + {\cal F}^D_{-m-1})
$$

Why did we take the associated graded of ${\cal G}^{(l)}_N$ 
for the weight and depth filtrations?
  The Lie algebra ${\cal G}^{(l)}_N$ is isomorphic to 
${\rm Gr}_{\bullet}^W{\cal G}^{(l)}_N$, 
but this isomorphism is not canonical. ${\rm Gr}_{\bullet}^W{\cal G}^{(l)}_N$ is a Lie 
subalgebra of ${\rm Gr}^W({\rm Der}^{SE}{\Bbb L}_N)\otimes \Q_l$. 
Via canonical isomorphism (\ref{PPPQ}) it became a Lie subalgebra of 
${\rm Der}_{}^{SE}{L}(\mu_N) \otimes \Q_l$, which has natural generators 
provided by the canonical generators of ${L}(\mu_N)$ (see s. 5.2). This gives canonical 
``coordinates'' for description of elements of ${\rm Gr}_{\bullet}^W{\cal G}^{(l)}_N$. Another reason to consider the Lie algebra 
${\rm Gr}_{\bullet}^W{\cal G}^{(l)}_N$ is provided by its motivic interpretation, see s. 3.7 below.  
The  
surprising benefit of taking its associated graded for the {\it depth} 
filtration is an 
unexpected relation with the geometry of modular varieties for $GL_m$, 
where $m$ is the depth. 

{\bf Remark}. Usually the associated graded for the depth filtration 
is not isomorphic to a subalgebra of a 
Lie algebra of any quotient of the Galois group. The situation is different 
in the following two  cases when 
the associated graded for the depth filtration is isomorphic to the 
original Lie algebra:

 i) {\it In the depth two case}  
${\rm Gr}{\cal G}_{\bullet, \geq -2}^{(l)}(\mu_N)$  
is isomorphic to  
${\rm Gr}^W{\cal G}^{(l)}_N/ {\cal F}_{-3}^D{\rm Gr}^W{\cal G}^{(l)}_N$  because there is no room for the difference. 

ii) {\it In the diagonal case}, as we noticed in s. 2.7, the diagonal 
Lie algebra defined as a Lie subalgebra of 
${\rm Gr}^W{\cal G}^{(l)}_{\bullet \bullet}(\mu_N)$, is non canonically 
isomorphic to a Lie subalgebra of ${\rm Gr}^W{\cal G}^{(l)}_N$.

{\bf 7.  Motivic interpretation of  the Lie algebra 
${\rm Gr}_{\bullet}^W{\cal G}^{(l)}_N$}. Let me recall the 
mixed Tate category of lisse $l$-adic sheaves 
defined by   Beilinson and Deligne, following  closely to the 
first 4 pages of 
chapter 1 of [BD]. 
For a connected coherent scheme $S$ over $\Z[1/l]$ such that 
$\mu_{l^{\infty}} \not \subset {\cal O}^*(S)$ denote 
by ${\cal F}_{\Q_l}(S)$ the Tannakian category of 
lisse $\Q_l$-sheaves on $S$.   
There are the  Tate sheaves $\Q_l(m):= \Q_l(m)_S := \Q_l(1)^{\otimes m}_S$. 
Thanks to our condition 
they 
 are mutually nonisomorphic. 
Call an object of  ${\cal F}_{\Q_l}(S)$ a mixed Tate object if 
it admits a finite increasing filtration $W$, indexed by $\Z$, such 
that ${\rm Gr}^W_{k} $ is a direct sum of copies of $\Q_l(-k)$. Let 
${\cal T}{\cal F}_{\Q_l}(S)$ be the full subcategory of mixed Tate 
objects in ${\cal F}_{\Q_l}(S)$. Then it is a Tannakian $\Q_l$-category, 
and  obviously one has
$$
{\rm Ext}^1_{{\cal T}{\cal F}_{\Q_l}(S)}(\Q_l(0), \Q_l(m)) =0 \quad 
\mbox{for $m\leq 0$}
$$
(This may not be so for the ${\rm Ext}^1$ in the category ${\cal F}_{\Q_l}(S)$).  
So ${\cal T}{\cal F}_{\Q_l}(S)$ is a mixed Tate category in the 
terminology of [BD]. 
It is easy to deduce from this 
that any its object admits a unique filtration $W$ such that 
${\rm Gr}^W_{k} $ is a direct sum of copies of $\Q_l(-k)$, 
and any morphism is strictly compatible with $W$. 
There is canonical 
functor $\omega$ to the category of graded $\Q_l$-vector spaces: 
$$
\omega: X \lms \oplus_m {\rm Hom}(\Q_l(m), {\rm Gr}^W_{-m}X)
$$
It is an exact functor commuting with $\otimes$-product. It is called canonical fiber functor. Let $L_{{\cal T}\Q_l}(S)$ be  
the Lie algebra of all $\otimes$-derivations of  $\omega$. 
By definition a degree $i$ element $\alpha \in L_{{\cal T}\Q_l}(S)_i$ 
is a collection of natural transformation $\alpha_j: \omega_j \lra \omega_{j+i}$
such that $\alpha_{M \otimes N} = \alpha_{M }\otimes {\rm Id}_{ \omega(N)} + 
{\rm Id}_{ \omega(M)} \otimes \alpha_{N }$ and $\alpha_{jM(1) }  = 
\alpha_{(j+1)M }$. Then 
$L_{{\cal T}\Q_l}(S)$ is 
graded pro-Lie algebra  over $\Q_l$. It is called the fundamental Lie algebra 
of the mixed Tate category 
${\cal T}{\cal F}_{\Q_l}(S)$. The fiber functor provides an equivalence between 
the category ${\cal T}{\cal F}_{\Q_l}(S)$ and the category of finite dimensional 
graded $L_{{\cal T}\Q_l}(S)$-modules. 

The standard Tannakian formalism in this case works as follows.  
Forgetting the graded structure on $\omega(X)$ we get a 
fiber functor $\widetilde \omega$.  
Its $\otimes$-automorphisms provide a group scheme over $\Q_l$ 
which is a direct product of ${\Bbb G}_m$ and a prounipotent group. 
Then $L_{{\cal T}\Q_l}(S)$ is the Lie algebra of the prounipotent part. 
The action of  ${\Bbb  G}_m$ provides a grading on it.

From now on let $S$ be the spectrum of the  ring of integers of a number field $F$ 
punctured in a finite set ${\cal S}$ containing all primes above $l$. Then  
${\cal F}_{\Q_l}(S)$ is identified with the category 
of finite dimensional $l$-adic representations of ${\rm Gal}(\overline F/F)$ 
unramified outside of ${\cal S}$. The underlying vector space of a Galois representation provides a fiber functor on ${\cal F}_{\Q_l}(S)$, and thus 
another fiber functor on the category ${\cal T}{\cal F}_{\Q_l}(S)$. 
It follows that  for any Galois representation 
$V$ from the category ${\cal T}{\cal F}_{\Q_l}(S)$ 
the Lie algebra of Zariski closure 
of the image of the Galois group in  ${\rm Aut}V$ is isomorphic to 
the image of the semidirect product of ${\rm Lie}{\Bbb G}_m$ and 
$L_{{\cal T}\Q_l}(S)$ acting 
on $\widetilde \omega(V)$. Moreover, 
let ${\cal G}_V$  be Zariski closure of the image of 
${\rm Gal}(\overline F/ F(\zeta_{l^{\infty}}))$ in ${\rm Aut}V$.  
Then the Lie algebra ${\cal G}_V$  
is isomorphic to the image of $L_{{\cal T}\Q_l}(S)$ in ${\rm Der}\omega(V)$. 
This isomorphism is canonical for 
 ${\rm Gr}^W{\cal G}_V$. 

As noticed in  s. 1.2.3 of [BD], it follows from a 
theorem of Soule [So] that 
the $l$-adic regulator map provides canonical isomorphism 
\begin{equation} \label{8.13.00.1}
H^i_{\cal M}(S, \Q(m))\otimes \Q_l := 
gr^{\gamma}_mK_{2m-i}(S) \otimes \Q_l 
\stackrel{=}{\lra} 
{\rm Ext}_{{\cal T}{\cal F}_{\Q_l}(S)}^i(\Q_l(0), \Q_l(m)) 
\end{equation}
Indeed, let $G_{F_{\cal S}}$ be the Galois group of of the maximal extension of $F$ unramified ourside ${\cal S}$. Soule proved that 
the $l$-adic regulator map provides isomorphism 
$$
H^i_{\cal M}(S, \Q(m))\otimes \Q_l \quad \stackrel{=}{\lra} \quad 
H^i(G_{F_{\cal S}}, \Q_l(m))
$$
for $i=1, m \geq 1$, and also for  $i=2, m \geq 2$ where both groups are zero: for the left one this follows from Borel's theorem;  
the $l=2$ case see s. B4 of [HW].
 
By the very definitions for $m \geq 0$
$$
{\rm Ext}_{{\cal T}{\cal F}_{\Q_l}(S)}^1(\Q_l(0), \Q_l(m)) = 
{\rm Ext}_{{\cal F}_{\Q_l}(S)}^1(\Q_l(0), \Q_l(m)) = H^1(G_{F_{\cal S}}, \Q_l(m)),
$$
so (\ref{8.13.00.1}) follows for $i=1$ case. Further, for $m \geq 2$
\begin{equation} \label{9.11.00.1}
{\rm Ext}_{{\cal T}{\cal F}_{\Q_l}(S)}^2(\Q_l(0), \Q_l(m)) \subset 
{\rm Ext}_{{\cal F}_{\Q_l}(S)}^2(\Q_l(0), \Q_l(m)) = H^2(G_{F_{\cal S}}, \Q_l(m)) = 0
\end{equation}
To check the 
 first inclusion notice the following.   
In any abelian category every element in the Yoneda 
${\rm Ext}^2$ is a product of ${\rm Ext}^1$'s. The  
Yoneda product of extension classes $\alpha \in {\rm Ext}^1(A,B)$ and 
$\beta \in {\rm Ext}^1(B,C)$ is zero in ${\rm Ext}^2(A,C)$ if and 
only if there exists an object $M$ with a three step filtration $F$ 
such that 
$F_0M = A, F_1M = \alpha$ and $M/F_0M = \beta$. 
In our case 
we may assume that a  class in 
${\rm Ext}^2(\Q_l(0), \Q_l(m))$ is a product the classes in 
 ${\rm Ext}^1(\Q_l(0), \Q_l(a))$  and ${\rm Ext}^1(\Q_l(a), \Q_l(m))$. 
Their product is zero in ${\rm Ext}_{{\cal F}_{\Q_l}(S)}^2$, 
and filtration $F$ on the corresponding object $M$ can be used as  a 
 weight filtration, so $M$ belongs to the category ${\cal T}{\cal F}_{\Q_l}(S)$. 
Finally, 
 ${\rm Ext}_{{\cal T}{\cal F}_{\Q_l}(S)}^2(\Q_l(0), \Q_l(m)) = 0$ 
for $m<2$ since, as before,  we can decompose it in a product of $Ext^1(\Q_l(0), \Q_l(k)$'s 
one of whom must have $k \leq 0$ and thus be zero. 
 Therefore 
$$
{\rm Ext}_{{\cal T}{\cal F}_{\Q_l}(S)}^2(\Q_l(0), \Q_l(m)) =0 \qquad 
\mbox{for all $m$}
$$
   This just means that  
the fundamental 
Lie algebra $L_{{\cal T}\Q_l}(S)$ is a free negatively 
graded 
Lie algebra generated in the degree $-m$ by the dual 
to $\Q_l$-vector space $K_{2m-1}(S) \otimes \Q_l$. 

Now let ${\cal S}':= {\cal S} - \{l\}$ and $S'$ be the spectrum of the ring of 
${\cal S}'$-integers, then the $l$-adic regulator map 
\begin{equation} \label{8.13.00.111}
H^i_{\cal M}(S', \Q(m))\otimes \Q_l \quad  
\stackrel{}{\lra} \quad
{\rm Ext}_{{\cal T}{\cal F}_{\Q_l}(S)}^1(\Q_l(0), \Q_l(m)) 
\end{equation}
is still isomorphism $i=1, 2$ and $m \geq 2$ since each of the 
groups does not change when we delete 
a closed point from the spectrum. 
If $i=1, m=1$ map (\ref{8.13.00.111}) is injective but not an isomorphism since 
 $$
H^1_{\cal M}(S', \Q(1))\otimes \Q_l = {\cal O}_{S'}^* \otimes \Q_l 
\hookrightarrow {\rm Ext}_{{\cal T}{\cal F}_{\Q_l}(S)}^1(\Q_l(0), \Q_l(1)) = 
E_S \otimes \Q_l
$$ where $E_S$ is the group of $S$-units. Thus 
$L_{{\cal T}\Q_l}(S')$ is a quotient of $L_{{\cal T}\Q_l}(S)$, and the space of 
generators of
 these
 Lie algebras differ only in the degree $-1$. 

Since ${\rm G}_m - \mu_N$ has a good reduction outside of $N$, 
in our case $S = S_N \cup \{l\}`$ and the Lie algebra 
${\Bbb L}_N$ is a pro-object in 
${\cal T}{\cal F}_{\Q_l}(S_N \cup \{l\})$. So the Lie algebra 
${\rm Gr}_{\bullet}^W{\cal G}^{(l)}_N$ is canonically isomorphic to 
the image of the Lie algebra $L_{{\cal T}\Q_l}(S_N \cup \{l\})$ acting by derivations of 
the pro-object 
$\Psi({\Bbb L}_N)$. 
However it follows from the distribution relations (see s. 5.7) that the degree
 $-1$ component of 
Lie algebra $L_{{\cal T}\Q_l}(S_N \cup \{l\})$ always act
 through its motivic quotient dual to the subspace ${\cal O}_{S'}^* \otimes \Q_l$. 

Recently some of the arguments similar to the one 
of Beilinson and Deligne described above were  rediscovered 
by Hain and Matsumoto [HM].   

Moreover let $L_{{\cal T}{\cal M}}(S)$ 
be the fundamental Lie algebra of the mixed Tate 
category of mixed Tate motives over $S$ ([G4]). Then 
$L_{{\cal T}\Q_l}(S) = L_{{\cal T}{\cal M}}(S)\otimes \Q_l$. 
We will not use this in the present paper. 
See also s. 3 of [G8] for the formalism of motivic Lie algebras 
in the  general setup. 

Deligne's arguments  show that  
(assuming the motivic formalism) the Lie algebra 
${\rm Gr}_{\bullet}^W{\cal G}^{(l)}_N$ is free for 
$N=2$ ([D2]) and $N=3,4$ ([D3]). Corollary \ref{4-19.10} 
implies that it is not free for sufficiently big $N$, for instance for 
prime $p>5$: the modular forms on $\Gamma_1(p)$ provide an  
obstruction to freeness.

 \section{The dihedral Lie coalgebra } 

{\bf 1. Definitions}.  Let $G$ and $H$ be two commutative groups, or more generally 
two commutative group schemes. Then, generalizing a construction given in [G3], 
one can define a graded Lie coalgebra ${\cal D}_{\bullet}(G|H)$, called the dihedral Lie algebra of 
 $G$ and $H$ (see [G4]). 
In the special case when $H = {\rm Spec}\Q[[t]]$ is the additive group of the 
formal line it essentially coincides with 
a bigraded Lie coalgebra  
${\cal D}_{\bullet \bullet}(G)$  defined in [G3] and called the dihedral Lie algebra of 
 $G$. (The second grading is related to the natural 
filtration on $\Q[[t]]$). The construction of  ${\cal D}_{\bullet}(G|H)$ is left as an 
 easy exercise, see the end of s. 4.5.

Below we give  a version of the definition of  the bigraded Lie coalgebra
$$
{\cal D}_{\bullet \bullet}(G ) = \oplus_{w \geq m \geq 1}
 {\cal D}_{w, m}(G )
$$
   given in [G3]. We will use it only when $G = \mu_N$. 

The $\Q$-vector space  $ {\cal D}_{w,m}(G)$ is generated by the symbols 
\begin{equation} \label{ginf}
I_{n_1,...,n_m}(g_1: ... : g_{m+1}),    \qquad w =
n_1+...+n_m, \quad n_i \geq 1
\end{equation} 
To define the relations we introduce the generating series 
\begin{equation} \label{ginf11}
\{g_1: ... : g_{m+1}|t_1:....:t_{m+1}\}:=  
\end{equation}
$$
\sum_{n_i >0} I_{n_1,...,n_m}(g_1: ... : g_{m+1})
(t_1-t_{m+1})^{n_1-1}...(t_m -t_{m+1})^{n_m-1}
$$
Thus
\begin{equation} \label{inversio}
\{g_1:  ... : g_{m+1}|t_1: ...: t_{m+1}\}  = \{g_1:  ... : g_{m+1}|t+t_1: ...: t+t_{m+1}\} 
\end{equation}
   We think about the generating series  (\ref{ginf11}) as a function of   $m+1$ pairs of variables   $(g_1,t_1), ... , (g_{m+1},t_{m+1})$ located cyclically on an oriented circle (and called a dihedral word) as follows. 
The oriented circle has slots, where the $g$'s sit, and in between the consecutive slots, 
dual slots, where $t$'s sit, see the picture below. The slots are marked by black points, and the 
dual slots by little circles.

\begin{center}
\hspace{4.0cm}
\epsffile{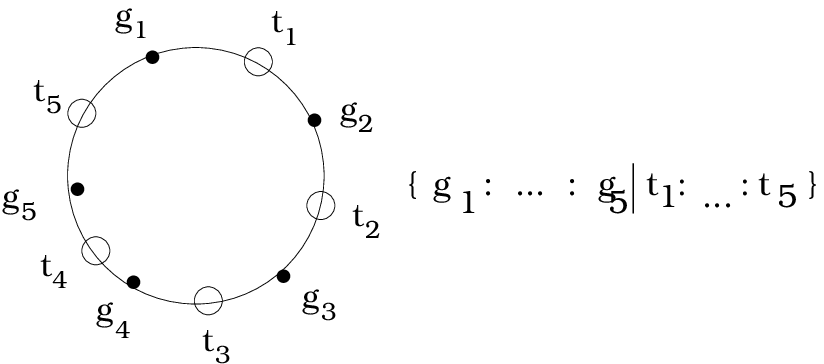}
\end{center}

We will also need two other generating series: 
 \begin{equation} \label{ginf112}
\{g_1: ... :  g_{m+1}|t_1,...,t_{m+1}\} := 
\end{equation} 
$$
\{g_1:   ... : g_{m+1}   
|t_1:t_1+t_2:...:t_1+...+t_{m}:0 \} 
$$
 where $t_1+ ...+ t_{m+1} =0$, and 
\begin{equation} \label{::}
\{g_1,   ... ,  g_{m+1}|t_1:...:t_{m+1}\} := \{1: g_1: g_1  g_2 : ... :  g_1   ...  g_{m }|
  t_1: ... :  t_{m+1} \}
\end{equation}
   where $g_1 \cdot  ...  \cdot  g_{m+1} = 1$. 

To make these definitions more transparent set  
$$
g_i' := g_i^{-1}g_{i+1}, \qquad t_i':= -t_{i-1} + t_i
$$
and put them on the circle together with $g$'s and $t$'s as follows

\begin{center}
\hspace{4.0cm}
\epsffile{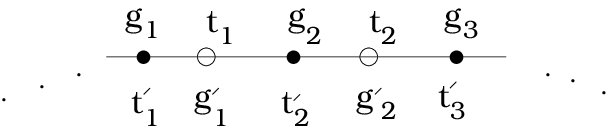}
\end{center}

Then it is easy to check that 
\begin{equation} \label{9.15.99.1}
\{g_1: ... :  g_{m+1}|t_1: ... : t_{m+1}\}  = 
\{g_1: ... :  g_{m+1}|t_1', ... , t_{m+1}'\} =  
\end{equation}
$$
\{g'_1, ... ,  g'_{m+1}|t_1: ... : t_{m+1}\}
$$
So we have 
\begin{equation} \label{:::}
\{g_1:  ... :  g_{m+1}|t_1,...,t_{m+1}\} = \{g'_{1},  ... ,  
g'_{m+1}|  t_1: t_1+t_2:...:t_1+...+t_{m}:0 \} 
\end{equation}
\begin{equation} \label{TRF}
\{g_1,  ... ,  g_{m+1}|t_1:...:t_{m+1}\} = \{1: g_1: g_1  g_2 : ... :  g_1   ...  g_{m }| 
t_1', t_2', ..., t'_{m+1}\}
\end{equation}
  
We picture the three generating series (\ref{9.15.99.1})  on  an oriented circle. 
Namely, for each of them we leave on the circle only the two sets of variables among 
$g$'s, $g'$'s, $t$'s, $t'$'s which appear in this generating series, see the picture. 
\begin{center}
\hspace{4.0cm}
\epsffile{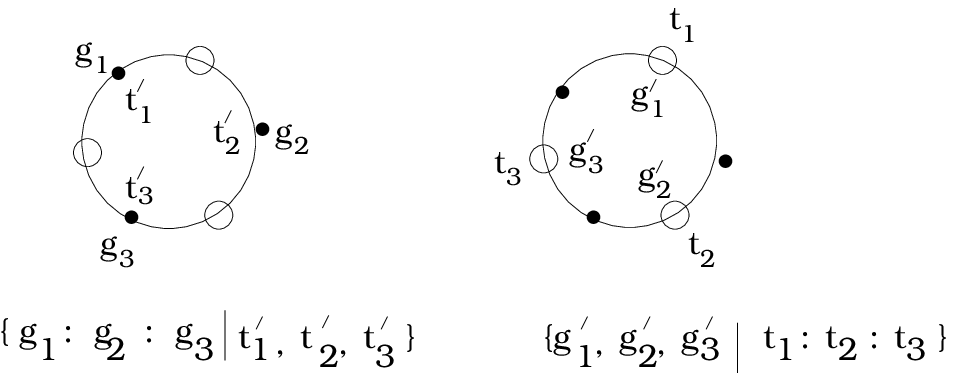}
\end{center}
The  ''$\{:\}$''-variables are outside, and the ``$\{,\}$''-variables are inside of the circle. (One can have the fourth type of the generating series, but it plays no role in our story). 
For the $\{g:|t,\}$-generating series  the variables sit only at slots, and for 
the $\{g,|t:\}$-generating series  only at dual slots.

{\bf Relations}. i) {\it Homogeneity}. For any $g \in G$ one has
\begin{equation} \label{gw0}
\{g\cdot g_1:   ... : g\cdot g_{m+1}|t_1: ... :t_{m+1}\} = 
\{g_1:  ... : g_{m+1}|t_1: ...: t_{m+1}\}  
\end{equation}
(A similar relation for $t$'s, when  $t_i \lms t_i +t$, is true by the definition (\ref{ginf11})). 

ii) {\it The double shuffle relations}  $(p+q = m, p\geq 1, q \geq 1)$. 
$$
s_1^{p,q}(g_1, ..., g_m, g_{m+1}|t_1: ... : t_m:t_{m+1}):=
$$
\begin{equation} \label{gw4}
 \sum_{\sigma \in \Sigma_{p,q}}
\{  g_{\sigma(1)}, ... , g_{\sigma(m)}, g_{m+1} | t_{\sigma(1)}: ...:  t_{\sigma(m)}: t_{m+1} \} =0
\end{equation}
$$
s_2^{p,q}(g_1: ... : g_m: g_{m+1}|t_1, ... , t_m, t_{m+1}):=
$$
\begin{equation} \label{gw3}
 \sum_{\sigma \in \Sigma_{p,q}}\{ g_{\sigma(1)} : ... :  g_{\sigma(m)}: g_{m+1}|   t_{\sigma(1)}, ...,  t_{\sigma(m)}, t_{m+1}\} =0
\end{equation}

iii) {\it The distribution relations}. Let $l \in \Z$. Suppose that the $l$-torsion subgroup $G_l$ of  $G$ is finite 
and its order is divisible by $l$. Then if $x_1, ..., x_m$ are $l$-powers 
$$
\{x_1:  ... : x_{ m+1}| t_1:  ...:  t_{ m+1 } \} -  
\frac{1}{|G_l|}\sum_{y_i^l = x_i}\{y_1: ... :y_{m+1}| l\cdot t_1: ...:  l\cdot t_{ m+1} \} =0
$$
  except that a constant (in $t$) is allowed when $m=1$ and $x_1 = x_2$, so that 
 $I_1(e:e) - \sum_{y^l=e} I_1(y:e)$ is not necessarily zero.

iv) $I_1(e:e) = 0$. (In [G3] we did not impose this condition). 

Denote by $\widehat {\cal D}_{\bullet \bullet}(G)$ the bigraded vector space defined 
by the conditions i)-iii) only. Then $\widehat {\cal D}_{\bullet \bullet}(G) = 
{\cal D}_{\bullet \bullet}(G) \oplus \Q_{(1,1)}$ where $\Q_{(1,1)}$ is generated by $I_1(e:e)$. 

{\bf Remark}. If $m \geq 2$, the relation (\ref{gw3}) for $g_1, ..., g_{m+1} = e$ implies that 
$I_{1, ..., 1}(e: ... :e) =0$, and iv) requires this vanishing to hold for $m=1$ as well.

{\bf Remark}. To define a map from ${\cal D}_{\bullet \bullet }(G)$ to a vector space $V$ amounts 
to define generating series $\{g_1: ... : g_{m+1}|t_1: ... : t_{m+1}\}$ with coefficients in $V$, 
that is, in the $V[[t_1, ..., t_{m+1}]]$, obeying (\ref{inversio}) and i) - v). 

{\bf Remark}. Altering one or both of the 
 conditions iii) and iv) we can still get important Lie coalgebras.

 {\bf 2.  The dihedral symmetry relations}. By definition 
they consist of  the following list of relations: 

{\it The cyclic symmetry relations}:
\begin{equation} \label{gw1}
\{g_1: g_2: ... :  g_{m+1}|t_1:t_2:...:t_{m+1}\} = \{g_2: ... : g_{m+1}: g_1|t_2:...:t_{m+1}: t_1\}  
\end{equation}

{\it The reflection relation}: 
\begin{equation} \label{gw2}
 \{g_1: ... :g_{m+1}|t_1:...:t_{m+1}\} = 
(-1)^{m+1}\{g_{m+1}:  ... : g_1|-t_{m}:...:-t_1: -t_{m+1}\}
\end{equation}

{\it The inversion relations}:
\begin{equation} \label{inversion}
\{x_1^{-1}: ... : x_{ m+1}^{-1}| t_1: t_{2}: ...:  t_{m+1} \} =
 \{x_1: ... :  x_{m+1} | - t_1: ...:  - t_{ m+1} \}
 \end{equation}
The inversion relations are 
precisely the distribution relations for $l=-1$.

One can check that the reflection relations (\ref{gw2})  for the $\{a,|b:\}$-generating 
series  look as follows:
\begin{equation} \label{reflexs}
\{g_1, ... , g_{m+1}| t_1: ...: t_{m+1}\} = (-1)^{m+1}\{g_{m+1}, ... , g_1| t_{m+1}: ...: t_1\}
\end{equation}
Using (\ref{inversion}), one checks the same identity for the $\{a:|b,\}$-generating series.  
This tells the effect of changing the orientation of the circle. 

The dihedral symmetry for $m=1$ reduces to the inversion relation. 

\begin{theorem} \label{zhdsrel} 
The double shuffle relations in the case $m\geq 2$ imply the 
 dihedral symmetry relations. 
\end{theorem}

{\bf Proof}. 
One has 
$$
\sum_{k=0}^m (-1)^k s_2^{k,m-k} (g_k: g_{k-1}:   ...   : g_1 : g_{k+1}:   ...  : g_m: g_{m+1}| 
$$
$$
t_k, t_{k-1},   ...   , t_1 , t_{k+1},   ...  , t_m, t_{m+1})
 = 0 
$$
Indeed, let us decompose the set of all shuffles of the ordered sets $\{k, k-1, ..., 1\}$ and 
$\{k+1, k+2, ..., m\}$ into a union of two subsets, denoted $S_k'$ and $S_k^{''}$. 
The subset $S_k^{''}$ consists of those shuffles  where $k+1$ appears on the left of $k$. 
After the summation with signs the terms corresponding to $S_k^{''}$ cancel the ones corresponding to 
$S_{k+1}'$. This means that
 \begin{equation} \label{sshh|}
\{g_1: g_{2}:   ...   : g_m: g_{m+1}| t_1, t_{2},   ...   t_m, t_{m+1}\} 
\quad = 
\end{equation} 
$$
(-1)^{m+1}\{g_m: g_{m-1}:   ...   : g_1: g_{m+1}| t_m, t_{m-1},   ...   t_1, t_{m+1}\} 
$$
Similarly 
\begin{equation} \label{sshh|1}
\{g_1, g_{2},  ...   ,g_m, g_{m+1}| t_1: t_{2}:   ... :  t_m: t_{m+1}\} = 
\quad 
\end{equation} 
$$
(-1)^{m+1} \{g_m, g_{m-1},  ...   , g_1,g_{m+1}| t_m: t_{m-1}:  ... :  t_1: t_{m+1}\} 
$$ 

To simplify the formulas below we use notation $h_i:= g_i' = g_i^{-1}g_{i+1}$. 
Using the identities we just get  we have 
$$
\{h_1,  ...   , h_{m+1}| t_1:   ... :  t_{m+1}\} \quad \stackrel {(\ref{9.15.99.1})}{=}\quad 
\{g_1:  ... :   g_{m+1}| t'_1,  ...,  t'_{m+1}\}  
\quad \stackrel {(\ref{sshh|})}{=}
$$
$$
(-1)^{m+1}\{g_{m}:  g_{m-1}: ... : g_1: g_{m+1} |
t'_m, t'_{m-1}, ..., ,t'_{1} , t'_{m+1} \} 
$$
$$
\stackrel {(\ref{:::})}{=} \quad (-1)^{m+1}\{h_{m-1}^{-1}, h_{m-2}^{-1}, ... , 
h_1^{-1}, h_{m+1}^{-1}, h_{m}^{-1}| 
$$
$$t_m - t_{m-1}: t_m - t_{m-2}: ... : t_m - t_{1}: t_m - t_{m+1}: 0 \} 
\quad \stackrel {(\ref{sshh|1})}{=}
$$
$$
\{h_{m+1}^{-1}, h_{1}^{-1}, h_2^{-1}, ... ,  h_{m-1}^{-1},  h_m^{-1} |
-t_{m+1}: -t_{1}: -t_{2}: ... :  -t_{m-1}: -t_{m}\} 
$$
Therefore we get 
\begin{equation} \label{fs}
\{g_1, g_2, ..., g_{m+1}|t_1: t_2: ... : t_{m+1}\} = 
\end{equation}
$$
\{g_{m+1}^{-1}, g_{1}^{-1}, g_2^{-1}, ... ,    g_m^{-1}|
-t_{m+1}: -t_{1}: -t_{2} : ... : -t_{m} \} 
$$

It is not difficult to check that  
\begin{equation} \label{s.shh+-} 
\{ g_1, ...,  g_{m+1}|t_1: ... : t_{m+1}\} + 
\end{equation}
$$
s^{m-1,1}_1(g_1, ..., g_{m-1}, g_{m+1}, g_m|t_1: ...: t_{m-1}: t_{m+1}: t_m ) - 
$$
\begin{equation} \label{5.30.00.1}
\{g_{m+1}, g_1, ...,  g_{m}|t_{m+1}: t_1: ... : t_{m}\} = 
\end{equation}
\begin{equation} \label{sshh+-}
\sum \{g_1, g_{\sigma(2)}, ... , g_{\sigma(m+1)}|
t_1: t_{\sigma(2)}: ... : t_{\sigma(m+1)}\}
\end{equation}
where the sum is over all shuffles of the sets $\{2,3,...,m\}$ and $\{m+1\}$. 
Expression (\ref{sshh+-}) does  not have the form of 
a shuffle relation. However 
applying (\ref{fs}) we see that it is equal to 
a shuffle relation. Therefore relation (\ref{s.shh+-}) = (\ref{5.30.00.1}) provides 
the cyclic symmetry for the generators (\ref{::}). 
Applying (\ref{fs}) we get inversion relations 
(\ref{inversion})  
for these generators. 
Finally, using (\ref{sshh|1}) and the cyclic relations 
we get reflection relation (\ref{reflexs}) 
for them. 
So we proved the dihedral symmetry  for  the 
 generators (\ref{::}), and thus for the generators (\ref{ginf112}) and (\ref{ginf11}). 
Theorem \ref{zhdsrel} is proved.

{\bf Proof of theorem \ref{DAHS}}. It is identical  to the proof of theorem \ref{zhdsrel} 
after we suppress $t'$s, and rename $g$'s by $v$'s, and use the additive notations instead of the multiplicative.

\begin{corollary} \label{8.26.00.100}
If $G$ is a trivial group then ${\cal D}_{w,m}=0$ if $w+m$ is odd. 
\end{corollary}

{\bf Proof}. One has $\{e: ... :e| t_1, ..., t_m,0\}= 
\{e: ... :e| -t_1, ..., -t_m,0\}$ by the inversion relation. This immediately leads to the statement.

{\bf 3. Motivation: relation with multiple polylogarithms when $G = \mu_N$}. There are several natural 
sets of the generators of the vector spaces ${\cal D}_{w, m}(G )$ which are reminiscent of different 
definitions of 
 multiple polylogarithms. It is useful to keep them in mind in order to trace the origins of the 
definitions given above. 

i) Let $x_1 ...x_{m+1} = 1$. We define the generators $L_{n_1, ..., n_m}(x_1, ..., x_m)$ as the 
coefficients of the 
generating  series (\ref{::}) when $t_{m+1} := 0$: 
$$
\{ x_1, ..., x_{m+1}| t_1: ... : t_m:0\} = : \sum_{n_i >0} L_{n_1, ..., n_m}(x_1, ..., x_m)
t_1^{n_1-1} ... t_m^{n_m-1}
$$
   If $G = \mu_N$ they are related to multiple polylogarithms
\begin{equation} \label{pas}
Li_{n_1, ..., n_m}(x_1, ..., x_m):= \quad 
\sum_{0 < k_{1} < k_{2} < ... < k_{m} } \frac{x_{1}^{k_{1}}x_{2}^{k_{2}}
... x_{m}^{k_{m}}}{k_{1}^{n_{1}}k_{2}^{n_{2}}...k_{m}^{n_{m}}}
\end{equation}
considered modulo the products of similar power series  and   modulo the 
lower depth power series. The shuffle relations (\ref{gw4}) correspond precisely 
to the shuffle product formula ([G3]) for this power series. 

The distribution relations are easy to check for the power series (\ref{pas}). 
%The constant $c_l$ in (\ref{EDR}) plays the role of $\log l$ in the relation 
%$$
%\log l = \sum_{a=1}^{l-1}\log (1 - \zeta_l^a) 
%$$
%if we 
%use the  dictionary $I_1(e:y) <-> -\log (1-y)$. So 
%relations (\ref{EDR}) reflect  $\log (lk) =\log (l) + \log (k)$.

ii) The generators $I_{n_1, ..., n_m}(a_1: a_2:  ... : a_{m+1})$ are the coefficients of the 
generating  series (\ref{ginf11}) when $t_{m+1} =0$: 
$$
\{  a_1: ... : a_{m+1}|t_1: ... : t_m:0\} = : \sum_{n_i >0}I_{n_1, ..., n_m}(  a_1: ... :   a_{m+1}) 
t_1^{n_1-1} ... t_m^{n_m-1} 
$$
When  $G = \mu_N$ they are related to the iterated integrals 
\begin{equation} \label{g21gr}
\int_0^{a_{m+1}}\underbrace{\frac{dt}{a_1-t} \circ \frac{dt}{ t} \circ ...\circ \frac{dt}{ t} }_{n_1   } \circ ... \circ \underbrace{\frac{dt}{a_m-t} \circ \frac{dt}{ t} \circ ...\circ \frac{dt}{ t} }_{n_m   }
\end{equation}
considered modulo the products of similar integrals  and   modulo the lower depth integrals. The shuffle relations (\ref{gw3}) correspond to the shuffle product formula for the iterated integrals (\ref{g21gr}), see theorem 2.2 in [G3] and 
more details in [G4]. Formula (\ref{::}) reflects  theorem 2.1 in [G3] relating 
power series (\ref{pas}) with iterated integrals (\ref{g21gr}). Formula 
(\ref{ginf112}) reflects the definition of the generating series $I^*$ which appear in 
theorem 2.2 in [G3]. 

{\bf Remark}. There is a striking symmetry between $g$'s and $t$'s in the generating series (\ref{ginf11}) or (\ref{ginf112}) and (\ref{::}). It is completely unexpected even when $G = \mu_N$: in this case $g$'s play role of variables of, say, iterated integrals (\ref{g21gr}), while $t$'s are formal parameters 
used to make the generating series. This symmetry is amplified by formulae for the coproduct given below: compare (\ref{ccc1}) 
and (\ref{ccc1q}). 

iii) Define 
$$
I_{  n_1, ..., n_m, n_{m+1}  }(a_1:  ... : a_m: a_{m+1}) \in {\cal D}_{w,m}(G) \quad \mbox{where} \quad 
w := n_1+...+n_{m+1}-1
$$ 
as the coefficients of the 
generating  series (\ref{ginf11}): 
\begin{equation} \label{frfr}
\{  a_1: ... : a_m: a_{m+1}|t_1: ... : t_m: t_{m+1}\} : =  
\end{equation}
$$
\sum_{n_i >0}I_{ n_1, ..., n_m, n_{m+1}}( a_1: ... : a_m: a_{m+1}) 
t_1^{n_1-1} ... t_m^{n_m-1} t_{m+1}^{n_{m+1}-1}
$$
More explicitly,
\begin{equation} \label{gauss4}
I_{  n_1 ,  ..., n_m,   n_{m+1}  }(a_1: ... : a_{m+1}) = 
 \end{equation}
$$
(-1)^{n_{m+1}-1} \cdot \sum_{i_1+...+i_m=n_{m+1}-1}{n_1+i_1 \choose  i_1}... {n_m+i_m \choose  i_m} I_{n_1+i_1, ... , n_m+i_m}(a_1: ... : a_{m+1})  
$$
Here $i_{k} \geq 0$. We recover (\ref{ginf}) for $n_{m+1} =1$. 

We obviously have the cyclic symmetry relations
\begin{equation} \label{ssa6}
I_{n_1,...,n_{m+1}}(a_1:... :a_{m+1}) = 
I_{n_2,...,n_{m+1}, n_1}(a_2:... :a_{m+1}:a_1)
\end{equation}

When $G = \mu_N$ the properties of the generators   $I_{n_1,...,n_m, n_{m+1}}(a_1: ... : a_{m+1})$ reflect  the properties of the iterated integral
\begin{equation} \label{g21}
\int_0^{a_{m+1}}\underbrace{ \frac{dt}{ t} \circ ...\circ \frac{dt}{ t} }_{n_{m+1}-1   }\circ \underbrace{\frac{dt}{a_1-t} \circ \frac{dt}{ t} \circ ...\circ \frac{dt}{ t} }_{n_1   } \circ ... \circ \underbrace{\frac{dt}{a_m-t} \circ \frac{dt}{ t} \circ ...\circ \frac{dt}{ t} }_{n_m   }
\end{equation}
This integral is divergent if $n_{m+1}>1$, and so has to be regularized (see theorem 7.1 in [G3] or [G1], or [G4]). 
     Its   regularized value   is given by formula (\ref{gauss4}).

The interpretation of the $\{g:|t:\}$-generating series on the circle is reminiscent of the structure 
of the iterated integral (\ref{g21}): the $\frac{dt}{t}$ differentials located between 
$\frac{dt}{t- a_i}$ and $\frac{dt}{t - a_{i+1}}$ correspond to the $t_i$-variable of the 
generating series, and thus sit between $a_i$ and $a_{i+1}$ on the circle.

{\bf 4. The cobracket 
$ 
\delta:     {\cal D}_{\bullet \bullet}(G)  \longrightarrow
  \Lambda^2{\cal D}_{\bullet \bullet}(G)$}. 
 It will be defined by
\begin{equation} \label{ccc3}
\delta \{g_1:  ... : g_{m+1}| t_1: ... :t_{m+1}\} = 
\end{equation}
$$
-\sum_{k=2}^{m} {\rm Cycle}_{m+1}\Bigl(\{g_{1}:... :g_{k-1}:g_k| t_{1}: ... : t_{k-1}: t_{m+1}    \} 
\wedge \{g_{k}  :   ... : g_{m+1}| t_k: ... : t_{m+1}  \} \Bigr)
$$ 
where the indices are modulo $m+1$ and  
$$
 {\rm Cycle}_{m+1} f(x_1,...,x_{m+1}) := \sum_{i = 1}^{m+1} f(x_i,...,x_{m+i})
$$  

{\bf Remark}. 
Let us add  formally to the generators the coefficients of the 
generating series $\{g_1|t_1\}$ 
and put them equal to zero. 
Then the sum in (\ref{ccc3}) will be from $1$ to $m+1$. 

Each term of the formula corresponds to the 
following procedure: choose a slot and 
a dual slot on the circle. 
Cut the circle 
 at the chosen  slot and dual slot  and make two   oriented circles      
with a dihedral words on each of them   out of the initial data. It is useful to 
think about the slots and dual slots as of little arcs, not points, so cutting one of them we 
get the 
arcs on each of the two new circles marked by the corresponding letters. 
The formula reads as follows: $$
\delta((\ref{ccc3})) = - \sum_{{\rm cuts}} \mbox{(start at the dual slot)} \quad  \wedge\quad 
\mbox{(start at the slot)}
$$ 
The only asymmetry between $g$'s and $t$'s is the order of factors.

\begin{center}
\hspace{4.0cm}
\epsffile{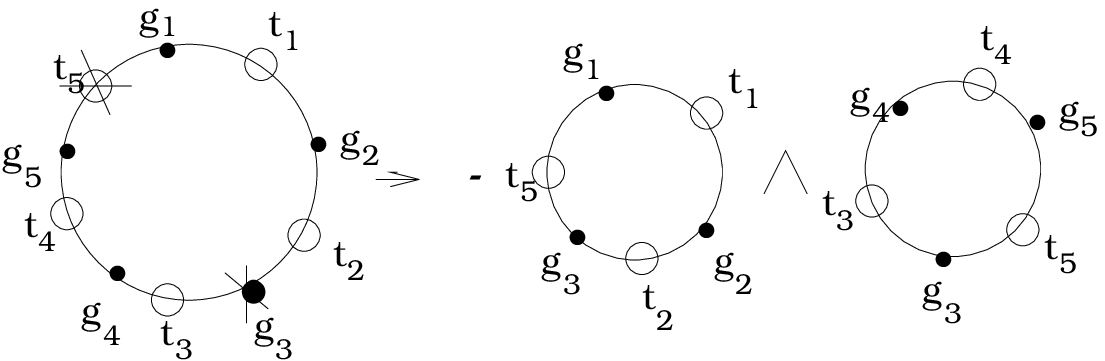}
\end{center}

The cobracket on the generators looks as follows. 
\begin{equation} \label{gauss6}
\delta   I_{n_1,...,n_{m+1}}(g_1:...:g_{m+1})= 
 \end{equation}
$$
 - {\rm Cycle}_{m+1}\Bigl(\sum_{k=2}^{m}\sum_{n'+n^{''}=n_{m+1}+1}  I_{n_1,...,n_{k-1},n'}(g_1:... :g_k) \wedge  
  I_{n_k,  ... , n_{m}, n''}(g_k:...:g_{m+1})\Bigr)
$$
 Notice that the second summation is over positive integers $n', n^{''}$ such that 
$
(n'-1) +( n^{''} - 1) = n_{m+1}  - 1
$.

  \begin{theorem} \label{9.99.18} a) There exists unique map $\delta: {\cal D}_{\bullet \bullet}(G) \lra 
\Lambda^2 {\cal D}_{\bullet \bullet}(G)$ for which (\ref{ccc3}) holds, providing 
a bigraded Lie coalgebra structure    on    
  ${\cal D}_{\bullet \bullet}(G)$.

b) A similar result is true for $\widehat {\cal D}_{\bullet \bullet}(G)$. Moreover there 
is an isomorphism of bigraded Lie algebras
$
\widehat {\cal D}_{\bullet \bullet}(G) =  {\cal D}_{\bullet \bullet}(G)\oplus \Q(1,1)
$ where $\Q(1,1)$ is a one dimensional Lie coalgebra of bidegree $(1,1)$. 
\end{theorem}

{\bf Remark}. The data of $f$ in $V[[t_1, ... , t_{m+1}]]$ is the same thing as the data, for any nilpotent ring $A$, and any nilpotent elements $t_1, ...,, t_{m+1}$ in $A$ of
$f(t_1, ..., t_{m+1}) \in A\otimes V $ functorial in $A$. Therefore to define $\delta: {\cal D}_{\bullet \bullet}(G) \lra 
\Lambda^2 {\cal D}_{\bullet \bullet}(G)$ it suffices to show that the {\it function} of  
$g, t$ with values in $\Lambda^2 {\cal D}_{\bullet \bullet}(G) \otimes A$ 
given by proposed $\delta\{g|t\}$ obeys the relations (i)-(v).

{\bf Proof}. Let ${\cal E}$ be the $\Q$-vector space generated by the coefficients 
of the $\{g:|t:\}$-generating series submitted only to the cyclic invariance (\ref{gw1}).  
To define a map ${\cal E} \lra V$ amounts to give generating series $\{g:|t:\}$ in $V[[t_1, ...]]$ obeying (\ref{gw1}), 
and the fact that (\ref{ccc3}) defines a map ${\cal E} \lra \Lambda^2{\cal E}$ is clear. 

{\it The Jacobi identity $\delta \delta =0$}. It holds in ${\cal E}$. One has to prove that \linebreak 
$\delta \delta \{g_1:  ... : g_{m+1}| t_1: ... :t_{m+1}\} =0$ in $\Lambda^3{\cal E}[[t_1, ... , t_{m+1} ]]$. 

Let $(1) \wedge (2)$ be a single term in  (\ref{ccc3}) corresponding to cuts at given slot and dual slot on the circle. 
Let us show that $\delta (1) \wedge (2) - (1) \wedge \delta (2)  = 0$: 
this implies the Jacobi identity.  
The terms in the expression 
$\delta (1) \wedge (2)$ could be of two different types, corresponding to 
the two situations shown on the left part of the picture. Here we marked by $I$ the initial cuts, so 
for $\delta (1) \wedge (2)$ the new cuts must be after the marked by $I$ dual slot. The four cuts define 
 four arcs on the circle, and computing $\delta (1) \wedge (2)$ we make three little circles out of them, wedged in a 
certain order. The numbers 
$1,2,3$ on the arcs indicate this order. 

  The terms in  $(1) \wedge \delta (2)$ correspond to the two drawings on the right of the picture. 
It is easy to see that the terms $N1$ and $N4$, and the terms 
$N2$ and $N3$, cancel each other (provided the cuts on the corresponding circles 
are made in the same slots and dual slots).

\begin{center}
\hspace{4.0cm}
\epsffile{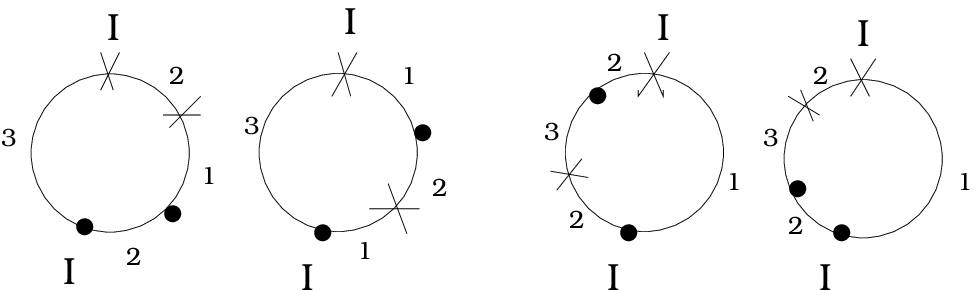}
\end{center}

It is easy to see that adding relations (\ref{gw0}) we still get a Lie coalgebra. 

{\it The shuffle relations}. We will show that the cyclic symmetry relations (\ref{gw1}) together with  each of the 
shuffle relations generate a coideal.  
 
  We start from the shuffle relations (\ref{gw3}). Let us impose (\ref{ginf112}) and (\ref{gw3}) 
to define the quotient 
 ${\cal E}_1$ of ${\cal E}$. To define a map ${\cal E}_1 \lra V$ amounts to define $\{...\}$ in 
$V[[t_1, ...]]$ obeying 
defining relations of ${\cal E}$ plus (\ref{ginf112}) and (\ref{gw3}).

We will 
need an explicit description of the cobracket for the generating series 
$\{g_1:| t_1,\}$ and $\{g_1, | t_1: \}$.

\begin{lemma} \label{ccc2} One has 
\begin{equation} \label{ccc1}
a) \qquad \qquad \qquad \delta \{g_1:  ... : g_{m+1}| t_1, ... ,t_{m+1}\}:= 
\end{equation} 
$$
-\sum_{k=2}^{m} {\rm Cycle}_{m+1}\Bigl(\{g_{1}:... :g_{k-1}:g_k| t_{1}, ... , t_{k-1}, x_k    \} 
\wedge 
$$
$$
\{g_{k}  : g_{k+1}: ... : g_{m+1}| y_k, t_{k+1}, ... , t_{m+1}  \} \Bigr)
$$
 where $t_1+...+t_{k-1} + x_k =0$, $y_k + t_{k+1}+ ...  +t_{m+1} =0$.  
\begin{equation} \label{ccc1q}
b) \qquad \qquad \qquad \delta \{g_1,  ... , g_{m+1}| t_1: ... :t_{m+1}\} = 
\end{equation} 
$$
\sum_{k=2}^{m} {\rm Cycle}_{m+1}\Bigl(\{g_{1},... ,g_{k-1}, a_k| t_{1}: ... :t_{k-1}: t_k    \} 
\wedge 
$$
$$
\{b_{k}, g_{k+1}, ... , g_{m+1}| t_k: t_{k+1}: ... : t_{m+1}  \} \Bigr)
$$
where $a_k$ and $b_k$ satisfy the conditions 
$g_{1} ...  g_{k-1} a_k = 1$ and $b_{k}g_{k+1} ... g_{m+1} = 1$. 
\end{lemma}

Notice the sign difference between the otherwise similar (\ref{ccc1}) and (\ref{ccc1q}). 
It might be explained by the different role the slots and dual slots play in the 
definition of $\delta$, as well as in the $\{a, | b:\}$ and $\{a: | b,\}$-generating series.  

{\bf Proof}. The proofs of these formulas are more or less identical, 
so we will present only the proof of the first one. 
We have 
$$
\delta \{g_1:  ... : g_{m+1}| t_1',  ... , t_{m+1}'\} =  
\delta \{g_1:  ... : g_{m+1}| t_1: ... : t_{m+1}\} =  
$$
\begin{equation} \label{ccc3q}
-\sum_{k=2}^{m} {\rm Cycle}_{m+1}\Bigl(\{g_{1}:... :g_{k-1}:g_k| t_1: ... :
t_{k-1}: t_{m+1}   \} \wedge 
\end{equation}
$$
\{g_{k}  : g_{k+1}:   ... : g_{m+1}| t_k: t_{k+1}: ...  :t_{m+1}  \} \Bigr)=
$$
$$
-\sum_{k=2}^{m} {\rm Cycle}_{m+1}\Bigl(\{g_{1}: ... : g_{k-1} : g_k| t_1', t'_2, ... , 
t'_{k-1}, -t_{k-1} + t_{m+1}  \} \wedge 
$$
$$
\{g_{k}  : g_{k+1} : ... : g_{m+1}| t_k- t_{m+1}, t_{k+1}', ... , t'_{m+1} \} \Bigr) 
$$
This is just what we wanted. The lemma is proved.

Recall that $(g_i|t_i)$ sits at the slot $i$.
A shuffle relation (\ref{gw3}) is determined 
 by choosing a slot, say $m+1$, and integer $1 \leq p \leq m-1$. Namely, we shuffle 
the slots numbered by 
\begin{equation} \label{ded}
\{1, ... , p\} \quad \mbox{and} \quad \{p+1 , ... , m\}
\end{equation}  
together with the symbols $(g_i|t_i)$ attached to them, and make generating series 
of type 
(\ref{ginf112})  out of each of the shuffles
 by putting it 
on the oriented circle 
before the slot $m+1$. We use a shorthand $s_1(1,...,p|p+1,...,m)$ for this shuffle relation.

Let us mark a slot $i$ and the dual slot $i+k$, which is between the slots  $i+k$ and $i+k+1$,  
see the picture. The marks determine a single term (the indices are modulo $m+1$)
\begin{equation} \label{dedqa}
(g_i: g_{i+1}:  ... : g_{i+k}| y_i, t_{i+1}, ... , t_{i+k}) 
\wedge (g_{i+k+1}: ...: g_{i-1}: g_i| t_{i+k+1}, ..., t_{i-1}, x_i) 
\end{equation}
in formula (\ref{ccc1}) for the coproduct, denoted $\partial_{[i,k]}(\{g_0: g_{1}: ... : g_{m}|t_0, t_{1}, ... , t_{m} \} )$. We will focus our attention on this term. 
\begin{center}
\hspace{4.0cm}
\epsffile{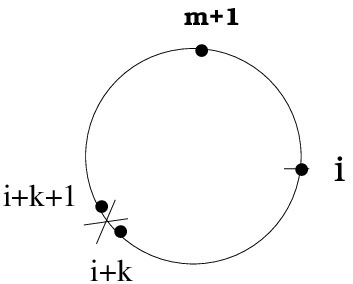}
\end{center}
Assume first that $m+1 \not \in \{i+1, ..., i+k\}$. 
Consider 
all the shuffles $\sigma$ of the sets (\ref{ded}) such that
$$
\{\sigma (1), ...,  \sigma (p)\}  \cap \{i+1, ... , i+k \} \quad \mbox{is a given set} \quad 
\{a, a+1, ..., b\} 
$$
Then the complement to $\{a, a+1, ..., b\}$ in $\{i, ..., i+k\}$ 
is also a given set,  denoted 
$\{\alpha,...,\beta\}$.  We assume that the slots $\{a,  ..., b\}$  as well as $\{\alpha,...,\beta\}$ are located on the circle in the order prescribed by the orientation of the circle. 
Then the terms in $\partial_{[i,k]}( s_1(1,...,p|p+1,...,m))$ corresponding to such $\sigma$'s are in the subspace 
$
s_1(a,...,b|\alpha,...,\beta) \wedge \{\mbox{all terms}\}
$. If $m+1 \in \{i+1, ..., i+k\}$ then  working with the 
second term in (\ref{dedqa}) we may argue 
 just as  before.

Using formula (\ref{ccc1q}) and formal symmetry $g <-> t$ we conclude that the cyclic symmetry together with shuffle relations (\ref{gw4}) also generate a coideal. 

{\it The distribution relations}. Assume first that there are no two consecutive equal elements 
 among $x_i$. Then computing 
$$
\delta \frac{1}{|G_l|} \sum_{y_i^l = x_i} \{y_1: ... : y_{m+1}|l t_1: ... : l t_{m+1} \}
$$
via formula (\ref{ccc3}), applying the distribution relations to the first factors, and after that to the second factors of $\delta(...)$ we get $\delta \{x_1: ... : x_{m+1}|t_1: ... : t_{m+1} \}$. 
However since the distribution relations for $\{e:e|t_1:t_2\}$ hold only up to a constant we need additional arguments  to
show that the total contribution of these constants is  zero. 

Suppose $x_1 = x_2 = e$. Let us show that if $m \not = 2$ or $ x_3 \not = e$ then 
\begin{equation} \label{9.20.99.1}
\{e:e|t_1:t_{m+1}\} \wedge \{ e: x_3: ... x_{m+1}|t_2: ... : t_{m+1}\} +
\end{equation}
\begin{equation} \label{9.20.99.2}
\{x_3: ... : x_{m+1}: e| t_3: ... : t_{m+1}: t_2\} \wedge \{e:e|t_1:t_{2}\} = 
\end{equation}
\begin{equation} \label{9.20.99.3}
\frac{1}{(|G_l|)^2} \sum_{y_i^l = x_i}\Bigl(\{y_1:y_2'| lt_1: lt_{m+1}\} \wedge 
\{ y_2: y_3: ... y_{m+1}|lt_2: ... : lt_{m+1}\} +
\end{equation}
\begin{equation} \label{9.20.99.4}
\{y_3: ... : y_{m+1}: y_1'| lt_3: ... : lt_{m+1}: lt_2\} \wedge \{y_1: y_2|lt_1: lt_{2}\} \Bigl)= 
\end{equation}
Here    (\ref{9.20.99.1}) $+$ (\ref{9.20.99.2}) is a sum of two appropriate terms in the formula for $\delta(...)$. 
We assume the summation over the arbitrary $l$-torsion elements $y_1'$ and $y_2'$. 
Applying the distribution relations we write the last two lines as
$$
\frac{1}{|G_l|} \sum_{y_i^l = x_i}\Bigl( \{e:e| t_1: t_{m+1}\} \wedge \{ y_2: y_3: ... y_{m+1}|lt_2: ... : lt_{m+1}\}+ 
$$
$$
\{y_3: ... : y_{m+1}: y_1'| lt_3: ... : lt_{m+1}: lt_2\} \wedge \{e: e|t_1: t_{2}\}\Bigl)+ 
$$
$$
+ \frac{1}{|G_l|} \Bigl(\Bigl(\sum_{y_1^l = e} I_1(y_1:e) - I_1(e:e) \Bigl)\wedge \sum_{y_i^l = x_i}
\{ y_2: y_3: ... y_{m+1}|lt_2: ... : lt_{m+1}\}+
$$
$$
\sum_{y_i^l = x_i} \{y_3: ... : y_{m+1}: y_1'| lt_3: ... : lt_{m+1}: lt_2\} \wedge 
(\sum_{y_1^l = e} I_1(y_1:e) - I_1(e:e)) \Bigl)
$$
The last two lines cancel each other thanks to the skewsymmetry relation in $\Lambda^2$, the cyclic symmetry
$$
\{ y_2: y_3: ... y_{m+1}|lt_2: ... : lt_{m+1}\} = \{y_3: ... : y_{m+1}: y_2| lt_3: ... : lt_{m+1}: lt_2\}
$$
and the observation that since $x_1=x_2=e$ we can replace summation over $y_1'$ by the summation over $y_2$. 
If $m \not = 2$ or $ x_3 \not = e$ we are done: applying the distribution relations to the first two lines 
 we get (\ref{9.20.99.1}) + 
(\ref{9.20.99.2}). 
Suppose that $m=2$ and $x_3=e$. Then in the last step we get  (\ref{9.20.99.1}) + 
(\ref{9.20.99.2}) 
plus the reminder terms
$$
\{e:e| t_1: t_{3}\} \wedge (\sum_{y_2^l = e} I_1(y_2:e) - I_1(e:e)) + 
$$
$$
(\sum_{y_3^l = e} I_1(y_3:e) - I_1(e:e)) \wedge \{e:e| t_1: t_{2}\} 
$$
In this case to get the total contribution of the reminder terms one has to apply the cyclic summation 
${\rm Cycle}_3$ along the $(y_1|t_1), ..., (y_3|t_3)$ variables. Therefore thanks to the skewsymmetry 
the total reminder term is zero. 
The distribution relations, and hence the part a) of the theorem are proved.

{\it Proof of the part b)}. The coproduct on $\widehat {\cal D}_{\bullet \bullet}(G)$ 
is well defined thanks to the part a) of the theorem. 
The decomposition into a direct sum 
of bigraded vector spaces is true by the very definitions. 
It remains to check that 
$$
\delta(\widehat {\cal D}_{\bullet \bullet}(G)) \subset \Lambda^2{\cal D}_{\bullet \bullet}(G)
$$
i.e. $I_{1, 1}(g: g)$ appears in $\delta I_{n_1,...,n_{m+1}}(g_1:...:g_{m+1})$ 
with the factor equal to zero. To calculate this factor we need only the following 
component of $\delta I_{n_1,...,n_{m+1}}(g_1:...:g_{m+1})$:
$$
 - {\rm Cycle}_{m+1}\Bigl(I_{n_1,1}(g_1: g_2) \wedge  
  I_{n_2,  ... , n_{m}, n_{m+1}}(g_2:...: g_{m+1})  + 
$$
$$
I_{n_1,  ... , n_{m-1}, n_{m+1}}(g_1: ... :g_{m-1}: g_m)\wedge I_{n_m,1}(g_m: g_{m+1}) \Bigr)
$$
 We rewrite its second term as 
$$
 - {\rm Cycle}_{m+1}\Bigl(
I_{n_3,  ... , n_{m+1}, n_{2}}(g_3: ... :g_{m+1}: g_1)\wedge I_{n_1, 1}(g_1: g_{2}) \Bigr)
$$
So if $g_1=g_2$ the factor appearing at $I_{1, 1}(g_1: g_{2})$ is zero 
thanks to the cyclic symmetry. The part b) is proved.

{\bf The dihedral Lie coalgebra ${\cal D}_{\bullet}(G|H)$}.  
Let $G$ and $H$ be two (finite) commutative groups. We define a  
Lie coalgebra ${\cal D}_{\bullet}(G|H)$, graded by the integers $m \geq 1$, as follows. 
The $\Q$-vector space 
${\cal D}_{m}(G|H)$ is generated by the symbols $\{g_1: ... :g_{m+1}|h_1: ... : h_{m+1}\}$ 
obeying the relations i), ii) and iii). We assume the homogeneity in both $g$'s and $h$'s. 
The distribution relations look as follows:
$$
\sum_{s_i^l = h_i}\{x_1: ... : x_m: 1| s_1: ... : s_m: 1\} = 
\sum_{y_i^l = x_i}\{y_1: ... : y_m: 1| h_1: ... : h_m: 1\} 
$$
This definition has an obvious generalization when $G$ or $H$ are commutative group schemes. 
(The role of constant in $t$ may play a constant function in $\Q[H]$).

{\bf Proof of the theorem \ref{6.15.00.1}}. It is completely similar to the proof of theorem 
\ref{9.99.18}a): suppress $t$'s from the notations, and denote $g$'s by, say, $v$'s understood as 
appropriate vectors of the lattice $L_m$, and use the additive notations instead of the multiplicative. The the proof goes literally the same.

{\bf 5. Two related bigraded Lie coalgebras: 
$\widetilde {\cal D}_{\bullet \bullet}(G)$ and 
${\cal D}'_{\bullet \bullet}(G)$}. 

i) The generators of $\widetilde {\cal D}_{w, m}(G)$  are
symbols
$$
   \widetilde I_{n_0,...,n_{m}}(g_0:...:g_{m})\quad  
\sum n_i = w-1, \quad n_i > 0, m\geq 1
$$
We package them into the generating series
\begin{equation} \label{12.23.99}
{\{g_0: ... : g_{m}| t_0: ... : t_{m}\} }\widetilde \quad:= 
\sum_{n_i>0} \widetilde I_{n_0,...,n_{m}}(g_0:...:g_{m}) t_0^{n_0} ... t_{m}^{n_{m}}
\end{equation}
The relations are the following: {\it the $G$-homogeneity}: for any $h \in G$
\begin{equation} \label{gauss10}
\{hg_0: ... : hg_{m}| t_0: ... : t_{m}\}\widetilde \quad = 
\{g_0: ... : g_{m}| t_0: ... : t_{m}\}\widetilde\quad  
\end{equation}
{\it the cyclic symmetry } (\ref{ssa6}), and in addition the constant term of (\ref{12.23.99}) 
is zero when the $g's$ are all equal, i.e. 
\begin{equation} \label{*!*!*}
\widetilde I_{1, ...,1}(e:...:e) =0
\end{equation}
The first part of the proof of theorem \ref{9.99.18} shows that formula (\ref{gauss6}) provides a bigraded 
Lie coalgebra structure on $\widetilde {\cal D}_{\bullet \bullet}(G)$. 

For $A$ the 
associative  algebra freely generated by a finite set $S$, the vector space  
$ 
{\cal C}(A):= A/[A, A]
$ has as basis the cyclic words in $S$. We denote by ${\cal C}: A \lra {\cal C}(A)$ the natural projection. 

The algebra $A$ is the universal enveloping algebra of the free Lie algebra generated by $S$, and as such has a commutative coproduct $\Delta$.   If we identify $A$ with its graded dual, using the basis afforded by words in $S$, this coproduct dualizes into the shuffle product
\begin{equation} \label{SH}
(s_1 ... s_k) \circ_{Sh}(s_{k+1} ... s_{k+l}) := \sum_{\Sigma_{k,l}}s_{\sigma(1)} ... 
s_{\sigma(k+l)}
\end{equation}

Let ${A}(G) $ be the free associative algebra with the generators $Y$ and $X_g$, $g \in G$. 
The group $G$ acts on the generators of  ${A}(G) $ as in s. 1.4. So it acts  
by automorphisms of ${A}(G) $. Let $\widetilde {\cal C}({A}(G))$ 
be the quotient of ${\cal C}({A}(G))$ by the subspace generated by $Y^n$, $X_g^n$, ${n \geq 0}$. 

\begin{lemma} \label{ill}
There is a canonical isomorphism of $\Q$-vector spaces 
\begin{equation} \label{ma}
\eta:  \widetilde {\cal D}_{\bullet \bullet}(G) \lra \widetilde {\cal C}({A}(G))_G
\end{equation}
 $$
\widetilde I_{n_0,...,n_{m}}(g_0:...:g_{m}) \lms {\cal C}\left( X_{g_0} Y^{n_0-1}... X_{g_{m-1}}
 Y^{n_{m-1}-1} X_{g_{m}} Y^{n_{m}-1}\right)
$$
\end{lemma}

{\bf Proof}. Cyclic invariance of $\widetilde I$ corresponds to cyclic words being considered, homogeneity in $g$ for $\widetilde I$ corresponds to taking coinvariants and the relation 
(\ref{*!*!*}) corresponds to the relation $X_g^n =0 (n \geq 1)$. We divided as well by 
$Y^n (n \geq 0)$ which  otherwise would have been left out of the image.

ii) 
Define a {\it shuffle relation} in $\widetilde {\cal D}_{\bullet \bullet}(G)$ as the image under the isomorphism $\eta^{-1}$ of 
\begin{equation} \label{ill1}
{\cal C}\Bigl(X_e \cdot \left\{ (Y^{n_0-1} X_{g_1} Y^{n_1-1}... X_{g_k}
 Y^{n_k-1}) \circ_{Sh}  (Y^{m_0-1} X_{h_1} Y^{m_1-1}... X_{h_l}
 Y^{m_l-1}) \right\}\Bigr)
\end{equation}  
where the expressions in each of the parentheses $(\cdot)$ are nonempty. 
We define 
${\cal D}'_{\bullet \bullet}(G)$ as the quotient of $\widetilde {\cal D}_{\bullet \bullet}(G)$ 
by  the subspace generated by  these  shuffle relations.

iii) Let ${\cal D}^{''}_{w, m}(G)$ be the $\Q$-vector space generated by the symbols 
(\ref{ginf}) 
subject to the relations (\ref{gw0}),(\ref{gw1}), (\ref{gw3}) and v). 
So it has the same generators as ${\cal D}_{w, m}(G)$, but the relations are relaxed.

\begin{proposition} \label{mfor1}
There is an isomorphism 
$i: {\cal D}^{''}_{w, m}(G) \lra {\cal D}'_{w, m}(G)$ given by 
\begin{equation}\label{mfor2}
i: I_{n_1,...,n_m}(g_1:...:g_{m+1}) \lms \widetilde I_{n_1,...,n_m,1}
(g_1:...:g_{m+1})
 \end{equation}
\end{proposition}

{\bf Proof}. {\it Comparison of the spaces ${\cal D}^{'}$ and ${\cal D}^{''}$}. 
Defining both in  terms of the $\{g:|t:\}$-generating series one requires:

for both: $G$-homogeneity, cyclic invariance, nullity of the constant terms of the 
$\{e: e: ... :e| \}$.  (This holds for ${\cal D}^{''}$ for $m=1$ by the definition, and 
for $m>1$ by the $\{g:|t,\}$-shuffle relation);

for ${\cal D}^{''}$: $t$-homogeneity, $\{g:|t,\}$-shuffle;

for ${\cal D}^{'}$: shuffle  (\ref{ill1}). 

{\it The map $i$ is well defined}. Observe that  (\ref{ill1}) 
for $X_e \cdot \{Y \circ_{sh}(...)\}$ gives in ${\cal D}^{'}$ that 
\begin{equation} \label{ikar6*}
\partial_t\{g_0: ... : g_{m}| t_0+t: ... : t_{m} +t\} =0 \quad \mbox{at} \quad t=0
\end{equation}
hence the $t$-homogeneity. 

Finally, one can show that  the $\{g:|t,\}$-shuffle relations (\ref{gw3}) go 
under the map $i$ precisely 
to the space 
\begin{equation} \label{ikar6}
\eta^{-1} (\mbox{the subspace of the shuffle relations (\ref{ill1}) with $n_0 = m_0 =1$ }  )
\end{equation}
(In fact in [G3] relations (\ref{gw3}) appeared as a way to write 
the shuffle relations (\ref{ill1}) with $n_0 = m_0 =1$). 

{\it Surjectivity of the map $i$}. 
For $m=0$ the $t$-homogeneity (\ref{ikar6*}) implies $\{g|t\}=0$, i.e.  
 $\widetilde I_n(g_0) =0$ in ${\cal D}^{'}$.  

The $t$-homogeneity (\ref{ikar6*}) in ${\cal D}'$ allows to express $\widetilde I_{n_0,n_1,...,n_m}(g_0:g_1:...:g_m)$ via $\widetilde I_{1,n_1,...,n_m}(g_0:g_1:...:g_m)$ for $m \geq 1$. More precisely, let $S$ be the subspace of $A(G)_G$ 
generated by elements $X_e \{Y\circ_{Sh}{\cal A}\}$. 
We write $a \stackrel{S}{= }b$ if $a-b \in S$. Then 
$t$-homogeneity is equivalent to the formula
\begin{equation} \label{ikar1} 
X_{e}Y^{n_0-1}X_{g_1} Y^{n_1-1}\cdot ...  \cdot X_{g_m} Y^{n_m-1} \stackrel{S}{= }
\end{equation}
$$
(-1)^{n_0-1}X_{e}X_{g_1} \left\{ (Y^{n_0-1})\circ_{Sh}  (Y^{n_1-1} X_{g_2}Y^{n_2-1} \cdot ...  \cdot 
X^{g_m}
 Y^{n_m-1}) \right\} 
$$
So the map $i$ is surjective. 

{\bf Remark}. Applying $\eta^{-1}$ to formula (\ref{ikar1}) we get relations (\ref{gauss4}).

{\it Injectivity of the map $i$}. 
We need to show that general shuffle relations (\ref{ill1}),
$$
X_{e}\left\{ (Y^{n_0-1}  X_{g_1} \cdot ...  ) 
\circ_{Sh} (Y^{m_0-1}  X_{h_1} \cdot ...  ) \right\}, 
$$
belong to the subspace of ${A}(G)_G$ generated by 
$S$ and the shuffle relations (\ref{ill1}) with $n_0 = m_0 =1$. 
Indeed, by (\ref{ikar1})
$$
X_{e}\left\{ Y^{n_0-1}  X_{g_1} \cdot ...  \right\}  = 
X_{e}\left\{ (Y \circ_{Sh}{\cal A}) \right\} + 
X_{e}\left\{ X_{g_1} \cdot ...  \right\}
$$
Since obviously 
$
X_{e}\left\{ (Y \circ_{Sh}{\cal A}_1)  \circ_{Sh}{\cal A}_2\right\}  = 
X_{e}\left\{ Y \circ_{Sh}  ({\cal A}_1  \circ_{Sh}{\cal A}_2)\right\} 
$
we see that for ${\cal B}:= Y^{m_0-1}  X_{h_1} \cdot (... ) $ one has
$$
X_{e}\left\{ (Y^{n_0-1}  X_{g_1} \cdot ...  ) \circ_{Sh} {\cal B}\right\} \stackrel{S}{=} 
X_{e}\left\{ (X_{g_1} \cdot ...  ) \circ_{Sh} {\cal B}\right\} 
$$
Similarly 
$X_{e}\left\{ {\cal B}\right\} \stackrel{S}{=} X_{e}\left\{ X_{h_1} \cdot ...\right\} $, thus
$$
X_{e}\left\{ {\cal B}\circ_{Sh} (X_{g_1} \cdot ...  ) \right\} \stackrel{S}{=} 
X_{e}\left\{ (X_{h_1} \cdot ...  )\circ_{Sh} (X_{g_1} \cdot ...  ) \right\} 
$$
The statement, and hence the part b) of the theorem are proved.

\begin{proposition} \label{mfor} The quotient ${\cal D}_{\bullet \bullet}'(G)$ of 
$\widetilde {\cal D}_{\bullet \bullet}(G)$ inherits from $\widetilde {\cal D}$ a cobracket. 
\end{proposition} 

{\bf Proof}. Formulas (\ref{ccc3}) and (\ref{gauss6}) define a cobracket on ${\cal D}^{''}$. 
This cobracket, transported by $i$ to ${\cal D}^{'}$, is induced from $\widetilde {\cal D}$. 
For another proof see in 
the end of s. `6.

The Lie coalgebra 
${\cal D}_{\bullet \bullet}(G)$ is a quotient of ${\cal D}'_{\bullet \bullet}(G)$.

\section{The dihedral Lie algebras and  special equivariant  derivations}

 {\bf 1. The special derivations and cyclic words (after [Dr], [K])}. Let $A$ be the free associative algebra generated by a 
finite set $S$. We are going to define a Lie algebra ${\rm Der}^S{A}$ of special derivations of $A$ and describe it via cyclic words in $S$. 
A derivation ${D}$ of the algebra $A$ is called special if there are elements  
$B_s \in A$ such that 
\begin{equation} \label{g2*}
{D}(s) = [B_s, s] \quad \mbox{and}   
\quad {D}(\sum_{s \in S} s) = 0
\end{equation}
Thus 
\begin{equation} \label{S3.1}
\mbox{a system $B = \{B_s\}$ of elements of $A$ with $\sum[B_s, s]=0$} 
\end{equation} 
defines a 
special derivation $D_B$. One has 
$$
[D_B, D_C] = D_{[B,C]}, \quad \mbox{ where} \quad  
$$
\begin{equation} \label{S3.2}
[B,C]_s:= D_B(C_s) - D_C(B_s) - [B_s, C_s]
\end{equation} 
Formula (\ref{S3.2}) defines a Lie bracket on systems $(B_s)$ obeying (\ref{S3.1}). Indeed, as was pointed out to me by the referee, one can interpret them as infinitesimal automorphisms of the structure consisting of 
$$
\mbox{algebra $A$, for each $s\in S$, left module $A$ (noted $A(s)$), }
$$
$$
\mbox{endomorphism $a \lms as$ of $A(s)$, element $\sum s$ in $A$}
$$
The bracket (\ref{S3.2}) is the bracket coming from this interpretation. 

Define a map of $\Q$-linear spaces 
$$
{\rm Cycl}: {\cal C}(A) \longrightarrow A:\quad s_1 ... s_k \lms \sum_{i=1}^k 
s_{i} ... s_{i-1+k}, \quad 1\lms 0
$$
given on the generators by the sum of the cyclic permutations. Then set
$$
\partial_s: s_1 ... s_k \lms \left\{ \begin{array}{ll}
s_2 ... s_k& \mbox{if} \quad s=s_1\\ 
 0&  \mbox{otherwise} \end{array} \right.
$$
Similarly define the map $\partial_s'$ by $s_1 ... s_k \lms s_1 ... s_{k-1} $ if 
$s_k = s$ and $0$ otherwise. Then the image of the map ${\rm Cycl}$ is precisely 
$\cap {\rm Ker}(\partial_s - \partial_s')$ in the positive degree part $A^+$ of $A$. 
If $x$ is in this image 
\begin{equation} \label{**a}
\sum[\partial_sx, s] = \sum \partial_s x \cdot s - \sum s \cdot \partial_s x = 
\sum \partial'_s x \cdot s - \sum s \cdot \partial_s x = x-x =0
\end{equation}

Set ${\cal D}_s:= \partial_s \circ {\rm Cycl}$.
Let us define a Lie bracket in ${\cal C}({A})$ by 
\begin{equation} \label{**}
 [C_1, C_2]:=  -\sum_{s \in S}{\cal C}\Bigl( [{\cal D}_s C_1, 
{\cal D}_s C_2] \cdot s\Bigr)
\end{equation}
Then ${\cal C}({A})$ is the product of the central $\Q\cdot 1$ and of ${\cal C}^+({A})$. 
The map
$$
\kappa': {\cal C}^+({A})  \lra \mbox{the Lie algebra (\ref{S3.1}) (\ref{S3.2})}; \qquad \kappa': x \lms {\cal D}_s (x)
$$
is an isomorphism of $\Q$-vector spaces. This follows from theorem 4.2 in [K]. 
For the convenience of the reader we reproduce the argument. 
Let $F^1(A):= \oplus_{s \in S} A \otimes d s$. There is a sequence 
\begin{equation} \label{fff}
0 \lra {\cal C}^+(A) = A/([A,A] +\Q\cdot 1) \stackrel{d}{\lra}   F^1(A) \stackrel{t}{\lra}  [A,A] \lra 0
\end{equation}
where $d({\cal C}):= \sum {\cal D}_{s} {\cal C} \otimes d s$ and $t[a \otimes d s ]:= 
[a, s]$. It is a complex by (\ref{**a}), and it is clearly exact from the left and right. There is an isomorphism of graded (by the weight) $\Q$-vector spaces
$
A/\Q\cdot 1 \lra F^1(A), \quad s_{1} ... s_{m} \lra s_{1} ... s_{m-1} \otimes d s_{m}
$.   
Thus the Euler characteristic of the weight $\geq -m$ part of (\ref{fff}) is zero. 
So the complex is exact. 

The map $\kappa'$ respects the Lie brackets. 
In particular this proves that (\ref{**}) satisfies the 
Jacoby identity. Thus $\kappa'$ is a Lie algebra isomorphism.

The centralizer of $s$ in ${A}$ is $\Q[s]$. 
So $Ker \kappa' = \oplus_{s \in S} \Q[s]$. Let 
\begin{equation} \label{M,M}
\widetilde {\cal C}({A}) := {\cal C}({A})/\oplus_{s \in S}\Q[s]
\end{equation}
 The map $\kappa'$ induces a Lie algebra morphism 
$\kappa: \widetilde {\cal C}({A}) \lra 
{\rm Der}^S({A})$.

\begin{proposition} 
\label{ssw} The morphism $\kappa: \widetilde {\cal C}({A}) \lra 
{\rm Der}^S({A})$ is an isomorphism. 
\end{proposition}

 {\bf Proof}. We need only to show that $\kappa$ is surjective. 
This follows from exactness of (\ref{fff}).

{\bf 2.  The   Lie algebra  of   special equivariant derivations}.   Now assume $S = \{0\} \cup G$ where $G$ 
is a finite commutative group. The corresponding algebra $A$ is the algebra $A(G)$ of 2.4, with 
generators $Y$ and $X_g$ ($g \in G$). Then $A(G)$ is graded by the weight and depth, and as before 
${\rm Der}^SA(G)$ inherits a weight grading compatible with a
 depth filtration. On $\widetilde {\cal C}(A(G))$ it corresponds to the weight grading and 
depth filtration induced from $A(G)$, both shifted by one. It is clear for the weight grading. 
For depth: if $x$ in $\widetilde {\cal C}(A)$ has depth $m$, $y := {\rm Cycl}(x)$ (in $A/\oplus_s \Q[s])$ has the same depth. If no 
$\partial_g y$ has depth $m-1$, it follows that $y$ has no word of depth $m$ except possibly 
$X_0^k$ in the case $m=0$, which is zero in $\widetilde {\cal C}(A(G))$ thanks to (\ref{M,M}).  
 
The map $\kappa$ provides a linear map 
$$
{\kappa}_G: \widetilde {\cal C}({A}(G))^G \lra   {\rm Der}^{S}{A}(G)^G = : {\rm Der}^{SE}{A}(G)
$$
 
Proposition \ref{ssw} implies that it is an isomorphism.

We will identify the vector space 
$\widetilde {\cal C}({A}(G))^G$ with its graded for the depth filtration: the graded for the depth filtration of the Lie algebra $\widetilde {\cal C}({A}(G))^G$ is $\widetilde {\cal C}({A}(G))^G$ with the Lie bracket given by the sum (\ref{**}), with $s=0$ omitted:
$$
[C_1, C_2]:=  -\sum_{g \in G}{\cal C}\Bigl( [{\cal D}_g C_1, 
{\cal D}_g C_2] \cdot X_g\Bigr)
$$

{\bf 3. Formulation of the result}. Recall that $Y$ and $X_g$, $g\in G$ are the generators of the algebra $A(G)$. Let $$
{\cal C}(X_{g_0} Y^{n_0-1} \cdot ... \cdot  
X_{g_m} Y^{n_m-1})^G:= \sum_{h \in G} {\cal C}(X_{hg_0} Y^{n_0-1} \cdot ... \cdot  
X_{hg_m} Y^{n_m-1})   
$$
The expressions on the left are parametrized by  
{\it cyclic $G$-equivariant words}, i.e. 
a $G$-orbits on the set of all cyclic words. 
Consider the following formal expression:
\begin{equation} \label{cw2}
\widehat \xi_G :=\sum \frac{1}{|{\rm Aut}{\cal C}|}
    \widetilde I_{n_0, ..., n_m}(g_0: ... : g_m) \otimes {\cal C}(X_{g_0} Y^{n_0-1} \cdot ... \cdot  
X_{g_m} Y^{n_m-1})^G
\end{equation}
where the sum is over all  cyclic $G$-equivariant words ${\cal C}$ in $X_g , Y$. The weight 
$1/|{\rm Aut}{\cal C}|$ of a given cyclic word ${\cal C}$ 
is the order of automorphism group of this 
cyclic word, as was pointed out by the referee.  

Applying the map $Id \otimes  {\rm Gr}({\kappa})$ we get 
a  bidegree $(0,0)$ element 
$$
\widetilde \xi_G \in \widetilde {\cal D}_{\bullet \bullet}(G) \widehat \otimes_{\Q} 
{\rm Gr}{\rm Der}^{SE}_{\bullet \bullet}{A}(G)
$$

  Since $G$ is finite ${\cal D}_{w, m}(G) $ is finite dimensional $\Q$-vector space. 
Let ${D}_{-w,-m}(G)= {\cal D}_{w,m}(G)^{\vee}$. Then ${D}_{\bullet \bullet}(G):= 
\oplus_{w,m \geq 1} {D}_{-w,-m}(G)$ is a bigraded Lie algebra. We similarly define 
$\widetilde {D}_{\bullet \bullet}(G)$ and ${D}'_{\bullet \bullet}(G)$. 

We may view the element $\widetilde \xi_G$ as a map between the bigraded $\Q$-vector spaces: 
 \begin{equation} \label{g3}
\widetilde \xi_G \in Hom_{\Q-Vect}(  \widetilde {D}_{\bullet \bullet}(G) , {\rm Gr}{\rm Der}^{SE}_{\bullet \bullet}{A}(G) 
)
\end{equation}

Notice that ${D}_{\bullet \bullet}(G) \subset {D}'_{\bullet \bullet}(G) \subset \widetilde {D}_{\bullet \bullet}(G)$.

  \begin{theorem} \label{ga}
 a) $\widetilde \xi_G: \widetilde {D}_{\bullet \bullet}(G) \stackrel{=}{\lra} {\rm Gr}{\rm Der}^{SE}_{\bullet \bullet}{A}(G)$ is an isomorphism of bigraded  Lie algebras.

b) Restricting $\widetilde \xi_G$ to ${D}'_{\bullet \bullet}(G) $ we get an isomorphism 
\begin{equation} \label{12.22.99.1}
\xi'_G: {D}'_{\bullet \bullet}(G) \stackrel{=}{\lra} {\rm Gr}{\rm Der}^{SE}_{\bullet \bullet}{L}(G)
  \end{equation}

c) Restricting $\xi_G'$ to ${D}_{\bullet \bullet}(G) $ we get an 
injective Lie algebra morphism
\begin{equation} \label{12.22.99}
\xi_G: {D}_{\bullet \bullet}(G) \quad {\hookrightarrow} \quad {\rm Gr}{\rm Der}^{SE}_{\bullet \bullet}{L}(G)
\end{equation}
\end{theorem}

{\bf 4. Proof}. a) 
 We start with a remark from linear algebra. Let $L_1$, $L_2$ be Lie algebras 
and $\zeta \in Hom_{\Q-Vect}(L_1,L_2) = L_1^* \otimes L_2$. Denote by 
$\delta: L_1^* \lra \Lambda^2L_1^*$ the Lie cobracket on $L_1^*$ and by 
$[,]:\Lambda^2L_2 \lra L_2$ the Lie bracket on $L_2$. 
Define a symmetric product $
(L_1^* \otimes L_2) \otimes (L_1^* \otimes L_2) \lra \Lambda^2L_1^* \otimes \Lambda^2L_2$ by 
$ (a_1 \otimes a_2) \circ (b_1 \otimes b_2):= (a_1 \wedge b_1) \otimes (a_2 \wedge b_2)$.  
Then $\zeta$ is a morphism of Lie algebras if and only if 
$$
(\delta \otimes id)(\zeta) = \frac{1}{2}(id \otimes [,])(\zeta \circ \zeta)
$$
If $\zeta = \sum_i A_i \otimes B_i $ then 
$
 (id \otimes [,])(\zeta \circ \zeta) = \sum_{i,j} A_i\wedge A_j \otimes [B_i, B_j]
$ 
 and so the condition is
\begin{equation} \label{galois5}
\sum_i\delta(A_i) \otimes B_i = \frac{1}{2}\sum_{i,j} A_i\wedge A_j \otimes [B_i, B_j]
\end{equation}
To check that $1/2$ is needed notice that if  $\{e_i\}$ 
be a basis of a Lie algebra $L$, $\{f^i\}$  the dual basis, and 
$[e_i, e_j] = \sum c_{ij}^k e_k$, then 
$\delta f^k = \frac{1}{2} \sum_{i,j=1}^n c_{ij}^k f^i \wedge f^j$. 

Let us suppose in addition that the Lie algebra bracket $[,]$ on $L_2$ is obtained by alternation of 
 a 
(non necessarily associative) product $*$ on the vector space $L_2$, i.e $[x,y] := x*y - y*x$.  Then 
 $\zeta$ is a Lie algebra morphism if and only if 
$$
(\delta \otimes id)(\zeta) = (id \otimes *)(\zeta \circ \zeta)
$$
 If $\zeta = \sum_i A_i \otimes B_i $ then it looks as follows
\begin{equation} \label{galois05}
\sum_i\delta(A_i) \otimes B_i = \sum_{i,j} A_i\wedge A_j \otimes B_i * B_j
\end{equation}

Let us return to our situation. The Lie bracket on  
${\rm Gr}\widetilde {\cal C}({A}(G))$ is given by formula (\ref{**}) 
where the sum is over $X_g$ only. It follows that the Lie algebra structure on 
the subspace ${\rm Gr}\widetilde {\cal C}({A}(G) )^G$ is obtained 
by alternation of the non associative product $*_1$ given by (see the picture below)
$$
{\cal C}(X_{g_0} Y^{n_0-1} \cdot ... \cdot  
X_{g_k} Y^{n_k-1})^G \ast_{1} {\cal C}(X_{h_0} Y^{m_0-1} \cdot ... \cdot  
X_{h_l} Y^{m_l-1})^G : =
$$
$$
\sum_{i,j}{\cal C}(X_e Y^{n_i-1} X_{g_i^{-1}g_{i+1}} Y^{n_{i+1}-1} \cdot ... \cdot  
X_{g_i^{-1}g_{i-1}} Y^{n_i + m_j-2} \cdot 
$$
$$
X_{h_j^{-1}h_{j+1}} Y^{m_{j+1}-1} \cdot ... \cdot  
X_{h_j^{-1}h_{j-1}} Y^{m_{j-1}-1})^G
$$
It can be defined by the same formulae for any abelian group $G$.
\begin{center}
\hspace{4.0cm}
\epsffile{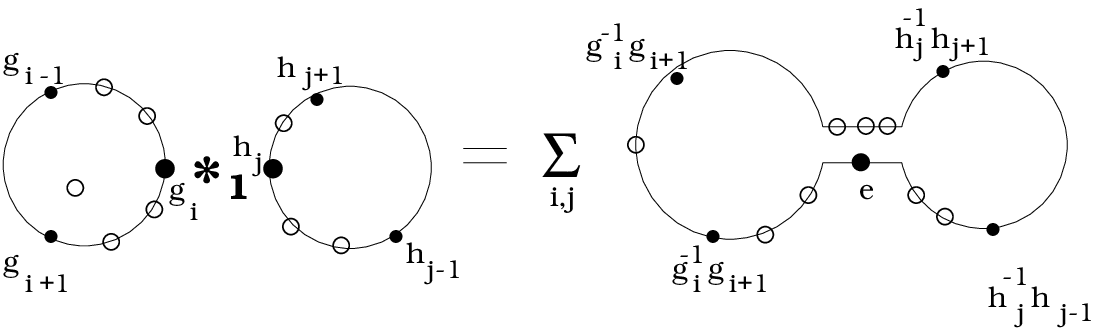}
\end{center}

To visualize a single expression  
\begin{equation} \label{cw1}
\delta 
(  \widetilde I_{n_0,n_1,...,n_m})(g_0:... :g_m) \otimes {\cal C}(X_{g_0} Y^{n_0-1} \cdot ... \cdot  
X_{g_m} Y^{n_m-1})^G 
\end{equation}
in $(\delta \otimes id)(\widehat \xi_G)$ we proceed as follows. 
Take an oriented circle divided into $ m+1 $ arcs by $ m+1 $ black points. Label each 
of the black points by an element of the set $\{X_g\}$, $g \in G$ and call it an $X_g$-point. 
The $i$-th arc is subdivided into $n_i$ little arcs by $ n_i-1 $  points labeled by $Y$ (presented by little circles on the picture and called  $Y$-points). Such a data corresponds to a cyclic word
$
{\cal C}(X_{g_0} Y^{n_0-1} \cdot ... \cdot  
X_{g_m} Y^{n_m-1}) 
$.

\begin{center}
\hspace{4.0cm}
\epsffile{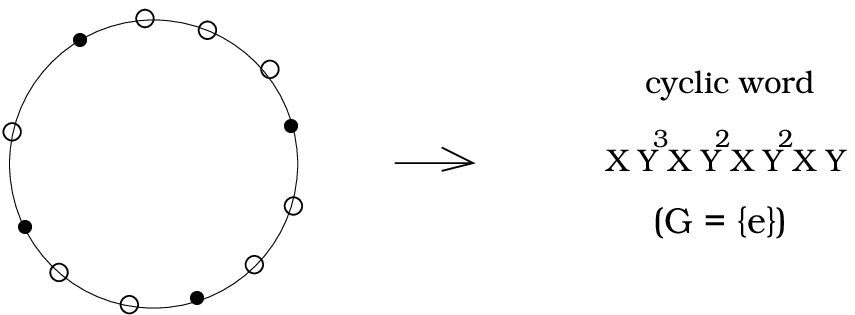}
\end{center}
The group $G$ acts naturally on them: 
an element $h \in G$ transforms  $X_g$-points to $X_{hg}$-points and leaves 
untouched $Y$-points.  
The orbits are called {\it circles  with $(n_0,...,n_m)$-cyclic $G$-structure}. 
They correspond to cyclic $G$-equivariant words, and thus to 
expressions (\ref{cw1}), as well as to the generators 
$\widetilde I_{n_0,n_1,...,n_m}(g_0:... :g_m) $ obeying the cyclic symmetry 
and homogeneity condition (\ref{gauss10}).

 Let us mark such a picture by choosing a little arc and a black point different from the ends of the arc containing the little arc. 
We call it a {\it marked circle  with $(n_0,...,n_m)$-cyclic $G$-structure}.
They correspond to {\it marked cyclic $G$-equivariant words}.

\begin{center}
\hspace{4.0cm}
\epsffile{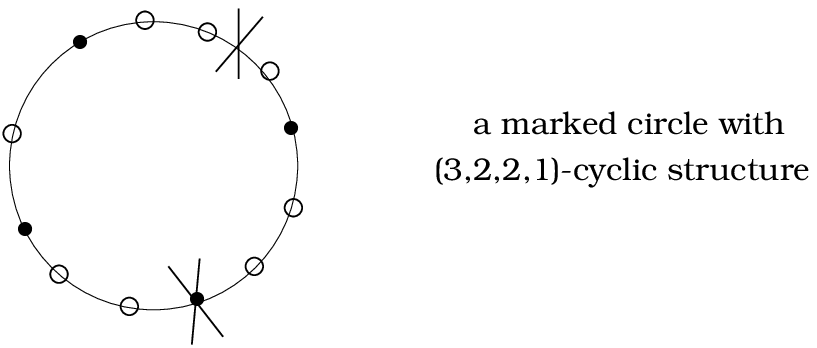}
\end{center}

It follows from formula (\ref{gauss6}) that  expression  (\ref{cw1}) is a sum of the terms 
which are in bijective correspondence with markings of the particular   
circle which 
corresponds 
to ${\cal C}(X_{g_0}Y^{n_0-1} ... )^G$.  
 For 
example the marked  circle with $(3,2,2,1)$-cyclic 
structure 
on the picture corresponds to 
$
I_{2,1 } \wedge I_{2,1,2} \otimes {\cal C}(XY^2XYXY^3XY^2)
$.     
In general we use  the marks  to cut the circle on $2$ oriented semicircles and   
make $2$ new circles by gluing the endpoints of the semicircles,  
adding a new black point on each of 
the new circles instead of the marked black point on the initial circle, and using 
the rest of the points on each of the new circles. The new circles 
are getting natural cyclic structures.

Going from a single expression (\ref{cw1}) to the sum, weighted by $1/{\rm Aut}{\cal C}$, of such expressions over  all isomorphism classes of {\it cyclic $G$-equivariant words} we get the  sum over all 
{\it marked cyclic $G$-equivariant words}. Notice that marking a particular 
cyclic word ${\cal C}$ we get a sum of marked cyclic words, each of them appearing precisely 
$|{\rm Aut}{\cal C}|$ times. The weight makes the coefficient equal to $1$.

Now let us investigate how  the right hand side of (\ref{galois05}) looks for our element 
$\widehat \zeta_G$. 
Consider the sum 
\begin{equation} \label{UP}
\sum_{} \frac{1}{|{\rm Aut} {\cal C}_P|}\frac{1}{|{\rm Aut} {\cal C}_Q|} 
 \widetilde I_{p_0,...,p_k}(g_0: ... : g_k) \wedge \widetilde I_{q_0,...,q_l}(h_0: ... : h_l) \otimes 
\end{equation}
\begin{equation} \label{UP1}
{\cal C}(X_{g_0} Y^{p_0-1}
\cdot ... \cdot X_{g_k} Y^{p_k-1})^G *_1 {\cal C}(X_{h_0} Y^{q_0-1}
\cdot ... \cdot X_{h_l} Y^{q_l-1})^G  
\end{equation}
where the summation is over all (ordered) pairs of cyclic $G$-equivariant 
words ${\cal C}_P= {\cal C}(X_{g_0} Y^{p_0-1}
\cdot ... \cdot X_{g_k} Y^{p_k-1})^G $ and ${\cal C}_Q = {\cal C}(X_{h_0} Y^{q_0-1}
\cdot ... \cdot X_{h_l} Y^{q_l-1})^G $.  

The particular product (\ref{UP1}) is a sum of $(k+1)(l+1)$ cyclic $G$-equivariant words.
 Each of them corresponds to a pair 
$$
\{ \mbox{\it a circle with $P:= (p_0, ..., p_k)$-cyclic $G$-structure $+$ a black point on it},
$$
 $$
\mbox{\it a circle with $Q:= (q_0, ..., q_l)$-cyclic $G$-structure $+$ a black point on it}\} 
$$ 
We call a circle with cyclic $G$-structure and a choice of a black point on it a {\it labelled circle with cyclic $G$-structure}. It corresponds to a {\it labelled cyclic $G$-equivariant 
word}.  

Now going from a single expression (\ref{UP1}) to the weighted 
sum over all pairs of cyclic $G$-equivariant 
words we get the terms corresponding to pairs of labelled cyclic $G$-equivariant 
words, each with the coefficient $1$ thanks to the weighting. 

The pairs of {\it labelled} cyclic $G$-equivariant 
words are in bijective correspondence with the {\it marked} cyclic $G$-equivariant 
words, as demonstrated on the picture using labelled or 
marked circles instead of the corresponding words:
\begin{center}
\hspace{4.0cm}
\epsffile{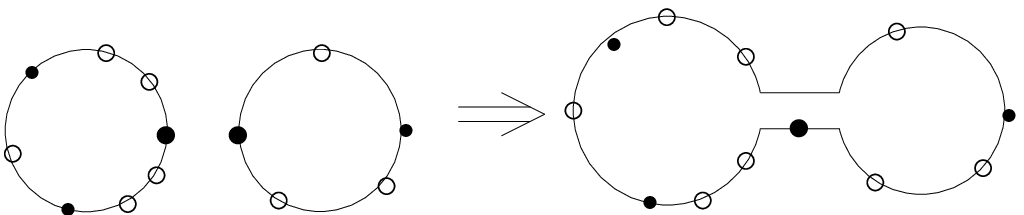}
\end{center}

\begin{center}
\hspace{4.0cm}
\epsffile{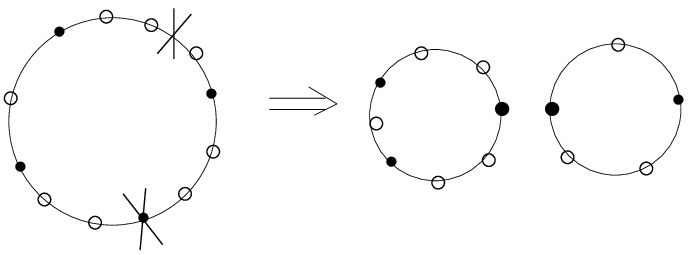}
\end{center}
Namely, to get a marked circle 
with 
cyclic equivariant $G$-structure we proceed as follows. Make a connected sum of the $P$ and $Q$-circles 
by connecting the chosen black points   by a bridge. The orientation of the initial circles induces 
an orientation of their connected sum. Instead of the chosen two black points on the 
initial circles we put a single black point on the bottom part of the bridge. We keep 
the rest of the 
points, getting a cyclic structure on the new circle. The black point on the bottom bridge  
together with the top part of the bridge provide the marks on the circle with the cyclic structure 
we get. The distinguished point on each of the three circles is chosen
 to be the $X_e$-point. 
The procedure is reversible; the inverse map is shown on the bottom of the picture.

We proved that $\widetilde \xi_G$ is a morphism of bigraded Lie algebras. By lemma \ref{ill} it is an isomorphism. The part a) of the 
theorem is proved. 

{\it Proof of the part b)}. 
\begin{lemma} \label{DD}
Let $D \in {\rm Der}^{SE}{A}(G)$. Then 
$D(X_e) \in {L}(G) <=> D \in {\rm Der}^{SE}{L}(G) $. 
\end{lemma}

{\bf Proof}. 
We need the following fact. Let the free associative algebra $A(S)$ contain the free Lie 
algebra $L(S)$. Then, if $x \in A(S)$ is such that $[x,s] \in L(S)$, one has 
$x\in L(S) + \Q[s]$. 
Indeed, by Poincar\'e-Birkhof-Witt, it suffices to check that the kernel of ${\rm ad} s$ 
acting on ${\rm Sym}^n(L(S))$ is reduced to $s^n$. Indeed, as ${\cal U}(\Q\cdot s)$ module, 
$L(S)$ is the sum 
of $\Q s$ and of a free module, and $L(S)^{\otimes n}$ hence the sum of 
$s^{\otimes n}$ and a free module. 

We are now ready to prove the only nontrivial implication, $=>$, of the lemma. 
Since $D$ is equivariant $D(X_e) \in {L}(G)$ implies that 
 $D(X_g) \in {L}(G)$ for any $g \in G$. 
Since $D(\sum X_g +Y) = 0$ we get also $D(Y) \in {L}(G)$. Thus thanks to 
the statement we just 
proved $D \in {\rm Der}^{SE}{L}(G)$.  The lemma is proved. 

Recall that $\widetilde {\cal C}(G)$ is a bigraded vector space and there is the isomorphism
\begin{equation} \label{DD+DD}
\kappa: \widetilde {\cal C}(G) \lra {\rm Der}^{S}A(G)
\end{equation}
The {\it vector space} ${\rm Der}^{S}A(G)$ admits, besides the weight grading, 
 a depth grading, with $D$ of degree $d$ if all 
$D(X_g)$ are of homogeneous degree $d+1$. Indeed, if $[B_0, Y] + \sum [B_g, X_g] =0$, the same 
holds with the $B_g$ (resp. $B_0$) replaced by their part of degree $d$ (resp. $d+1$). 
The isomorphism (\ref{DD+DD}) is compatible with the weight and depth grading (after a shift by one in 
$\widetilde {\cal C}(G)$). Further, ${\rm Der}^{S}L(G) \hookrightarrow {\rm Der}^{S}A(G)$, defined by 
``the $D(X_g)$ are primitive'' is a bigraded subspace. The action of the group $G$ preserves 
the gradings. Thus ${\rm Der}^{SE}L(G) \hookrightarrow {\rm Der}^{SE}A(G)$ is also a bigraded subspace. 
Notice that the Lie algebra structure on either sides of (\ref{DD+DD}) 
does not respect the depth grading, only the depth filtration. 

Let $X \in \widetilde {D}_{-w, -m}(G)$. Then we claim that 
\begin{equation} \label{DDDD}
\widehat \xi_G(X) \in {\rm Der}^{SE}L(G) \quad <=> \quad X \in \widehat  {D}'_{-w, -m}(G)
\end{equation} 

Recall that there is a coproduct $\Delta$ on ${A}(G)$ which 
dualizes to the shuffle product 
$\circ_{Sh}$, see (\ref{SH}). 
Let $\overline \Delta(Z) := \Delta(Z) - 1 \otimes Z  - Z \otimes 1$. 
Then $L(G) = {\rm Ker}\overline \Delta$. 
To prove (\ref{DDDD}) consider the expression
\begin{equation} \label{DDD1}
(id \otimes \overline \Delta \circ {\cal D}_{X_{e}})
(\widehat \xi_G) 
\in \widetilde {\cal D}_{-w, -m}(G) \otimes
\Bigl({A}(G) \otimes {A}(G)\Bigr) 
\end{equation}
Choose two elements $A$, $B$ of the natural basis in $A(G)$. 
Then $A \otimes B \in {A}(G) \otimes {A}(G)$ appears in $(id \otimes 
\overline \Delta \circ {\cal D}_{X_{e}})(\widehat \xi_G)$
with coefficient $X_e (A \circ_{Sh} B)$. This implies (\ref{DDDD}). 
So map (\ref{DD+DD}) leads to  an isomorphism of the bigraded {\it vector spaces}
$$
\widehat \xi_G: D_{\bullet \bullet}'(G) \lra {\rm Der}^{SE}L(G)
$$
After taking the associated graded for the depth filtration it becomes  an isomorphism 
of bigraded {\it Lie algebras} denoted $\xi'_G$.

Summarizing we see that  shuffle relations (\ref{ill1}) are equivalent to the condition 
$\widetilde \xi_G({D}'_{\bullet \bullet}(G)) \subset {\rm Gr}{\rm Der}_{\bullet \bullet}^{SE}
{ L}(G) $. The part b) of the 
theorem is proved.

{\bf Another proof of proposition \ref{mfor}}. We just proved that $\xi'_G$ is an isomorphism of bigraded 
$\Q$-vector spaces. 
Using a) (and an obvious  fact that ${\rm Der }^{SE}{ L}(G)$ is a Lie algebra!) 
we conclude that 
${D}'_{\bullet \bullet}(G)$ is a Lie subalgebra of $\widetilde {D}_{\bullet \bullet}(G)$.

c) It follows from a), b) and propositions \ref{mfor1} and \ref{mfor}.
The theorem is proved.

{\bf 5. The Lie algebra of outer semi-special derivations}. 
Let 
$$
{\cal X}_N:= \{0\} \cup \{\infty\} \cup \mu_N
$$
 For every point $x \in {\cal X}_N$ there is a 
well defined conjugacy class in $\pi_1(X_N, z)$ provided by a little lop around $x$. 

For any point $y \in {\cal X}_N$ there is a tangent vector $v_y$ at $y$ provided 
by the canonical coordinate $t$ on $X_N$ if $y \not = \infty$ and $t^{-1}$ at $\infty$.
Similarly for the pro-$l$ group 
$\pi^{(l)}_1(X_N, v_y)$ there is a well defined 
up to a conjugation map:
\begin{equation} \label{8.20.00.1}
i_x: \Z_l(1) \lra \pi^{(l)}_1(X_N, v_y)
\end{equation}
Therefore the  Galois group ${\rm Gal}(\overline \Q/ \Q(\zeta_{l^{\infty}N}))$ 
acts on $\pi^{(l)}_1(X_N, v_y)$ preserving the 
the conjugacy classes of the maps $i_x$, $x \in {\cal X}_N$. 
More precisely, define a {\it semi-special 
automorphism} of the group  $\pi^{(l)}_1(X_N, v_y)$ as the automorphism  
preserving all the conjugacy classes of 
the maps $i_x$ for $x \in {\cal X}_N$.
Denote by ${\rm Aut}^{SS}\pi^{(l)}_1(X_N, v_y)$ the group of all semi-special automorphisms 
of $\pi^{(l)}_1(X_N, v_y)$. We define the group 
 of {\it outer} semi-special automorphisms  ${\rm Out}^{SS}\pi^{(l)}_1(X_N)$ 
as the quotient of the group ${\rm Aut}^{SS}\pi^{(l)}_1(X_N, v_y)$ modulo the 
inner automorphisms. The corresponding groups defined 
for different base points $v_y$ are canonically isomorphic, so we can 
drop the base vector from the notations.  We get a canonical homomorphism
\begin{equation} \label{8.20.00.2}
\Psi_N: {\rm Gal}(\overline \Q/ \Q(\zeta_{l^{\infty}N})) \quad \lra \quad {\rm Out}^{SS}\pi^{(l)}_1(X_N)
\end{equation}

Similarly we define the Lie algebra ${\rm ODer}^{SS}{\Bbb L}(X_N)$ of outer semi-special 
derivations 
of the fundamental Lie algebra of $X_N$:
\begin{equation} \label{8.21.00.1}
{\rm ODer}^{SS}{\Bbb L}(X_N):= \quad \frac{{\rm Der}^{SS}{\Bbb L}(X_N, v_y)}
{{\rm InDer}{\Bbb L}(X_N, v_y)}
\end{equation}
Here the denominator is the Lie algebra of all inner derivations, and the numerator is the Lie algebra of semi-special derivations, i.e. derivations 
preserving the conjugacy classes of loops around $x \in {\cal X}_N$.

For any point $x_0 \in {\cal X}_N$ there is a Lie subalgebra 
\begin{equation} \label{8.20.00.3}
{\rm Der}^{S}{\Bbb L}(X_N, 
v_{x_0}) \quad \hookrightarrow \quad {\rm Der}^{SS}{\Bbb L}(X_N, 
v_{x_0})
\end{equation} 
of all derivations {\it special} with respect to the point $x_0$. 
It consists of the semi-special derivations killing  $i_{x_0}(\Q(1))$. 
Thus there  is a natural Lie algebra homomorphism
$$
p_{{x_0}}: {\rm Der}^{S}{\Bbb L}(X_N, v_{x_0}) \quad \lra \quad {\rm ODer}^{SS}{\Bbb L}(X_N)
$$
It is obviously surjective. Indeed, 
if ${\cal D}$ is a semi-special derivation then there exists $C \in {\Bbb L}(X_N, v_{x_0})$ 
such that 
${\cal D}(i_{x_0}(\Q(1))) = [C, i_{x_0}(\Q(1))]$, so subtracting the inner derivation 
$X \lms [C,X]$ we get a derivation special with respect to $x_0$. 
Since the pro-nilpotent 
Lie algebra ${\Bbb L}(X_N, v_{x_0})$ is free, ${\rm Ker}p_{x_0}$ is one dimensional, and 
consists of inner derivations provided by commutator with $i_{x_0}(\Q(1))$. 

There is also the $\Q_l$-version of this story. We denote by 
${\rm ODer}^{SS}{\Bbb L}^{(l)}(X_N)$ the $l$-adic analog of  (\ref{8.21.00.1}). 
%Linearization of the map (\ref{8.20.00.2}) provides a Lie subalgebra
%$$
%{\cal G}^{(l)}_{{\rm out}}(\mu_N) \quad \hookrightarrow \quad {\rm ODer}^{SS}{\Bbb L%}^{(l)}(X_N)
%$$
%Similarly linearization of the action of Galois group ${\rm Gal}(\overline \Q/ \Q(\z%eta_{l^{\infty}N}))$ on ${\Bbb L}^{(l)}(X_N, v_{x_0})$ provides a Lie subalgebra 
%$$
%{\cal G}^{(l)}_{{\rm x_0}}(\mu_N) \quad \hookrightarrow \quad 
%{\rm Der}^{S}{\Bbb L}^{(l)}(X_N, v_{x_0})
%$$
%The composition 
%$$
%{\cal G}^{(l)}_{{\rm x_0}}(\mu_N) 
%\quad \hookrightarrow \quad {\rm Der}^{S}{\Bbb L}^{(l)}(X_N, x_0) \quad \stackrel{p_%{x_0} \otimes \Q_l}{\lra} \quad {\rm ODer}^{SS}{\Bbb L}^{(l)}(X_N)
%$$
%provides a surjective Lie algebra morphism
%$$
%{\cal G}^{(l)}_{{\rm x_0}}(\mu_N) \lra {\cal G}^{(l)}_{{\rm out}}(\mu_N)
%$$
% Its kernel is at most one dimensional. We will see later on that it is one dimensio%nal 
%if $N\not =1$, 
%and zero for $N=1$. 

{\bf 6. A symmetric construction of the Lie algebra of special derivations}. 
Below we present a variation on the theme discussed in the previous subsection, 
in a bit different (and general) setup. 

Denote by $A_{{\cal X}}$ the free associative algebra generated by the set 
${\cal X}$. Denote by $X_x$ the generator corresponding to $x \in {\cal X}$. 
Let ${\cal X}$ be a finite set. Choose an element $x_0 \in {\cal X}$. Let 
$S:= {\cal X} - \{x_0\}$. Then there is the Lie algebra ${\rm Der}^S(A_{S})$ of 
special derivations of the algebra 
$A_{S}$. It contains the one dimensional subspace generated by 
inner derivation $\ast \lms [X_{x_0}, \ast]$. 
Our goal is to describe the structure of the quotient of  
the Lie algebra ${\rm Der}^S(A_{S})$ by this subspace 
in terms of the set ${\cal X}$ independently  
of the choice of 
an element $x_0$. 
We will use this for ${\cal X}:= {\cal X}_N$ to describe
 the corresponding quotient of the Lie algebras 
${\rm Der}^S{\Bbb L}(X_N, v_y)$, where  $y \in  {\cal X}_N$. As a result we get 
 natural {\it canonical isomorphisms} between these quotients.

Let $\overline A_{{\cal X}}$ be the quotient of 
$A_{{\cal X}}$ by the ideal $I({\cal X})$ generated by the 
element $ \sum_{x \in {\cal X}} X_x$. 
Then for any element $x \in {\cal X}$ there is a canonical isomorphism 
\begin{equation} \label{8.17.00.3} 
%{\rm iso}_x: 
\overline A_{{\cal X}} \quad \stackrel{=}{\lra} \quad A_{S}
\end{equation}

A semi-special derivation of the algebra $\overline A_{{\cal X}}$ is a derivation ${\cal D}$ 
preserving the conjugacy classes of the generators $X_x$, i.e. ${\cal D}(X_x) = [C_x, X_x]$
 for certain $C_x \in \overline A_{{\cal X}}$. 
The semi-special derivations form a Lie algebra  ${\rm Der}^{SS}(A_{{\cal X}})$. 
Let ${\rm In Der} (\overline A_{{\cal X}})$ be the Lie subalgebra of inner derivations. The quotient 
$$
{\rm ODer}^{SS}(\overline A_{{\cal X}}) := \quad \frac{{\rm Der}^{SS}
(\overline A_{{\cal X}}) }{{\rm In Der} (\overline A_{{\cal X}})}
$$
is called the Lie algebra of outer semi-special derivations.   There is canonical surjective homomorphism
$$
p_{{x_0}}: {\rm Der}^{S}A_{S} \quad \lra \quad {\rm ODer}^{SS}(\overline A_{{\cal X}})
$$
whose kernel is the subspace generated by the inner derivations 
corresponding to elements $(\sum_{x \in {\cal X}} X_x)^n$, $n \geq 1$. 

There is a similar story for the  Lie algebras. Namely, let $L_{{\cal X}}$ be the free Lie algebra generated by the set ${\cal X}$. Denote by $\overline L_{{\cal X}}$ its quotient by the ideal generated by the element $ \sum_{x \in {\cal X}} X_x$. 
Then for any $x \in {\cal X}$ 
there is canonical isomorphism $\overline L_{{\cal X}} \lra L_S$. We 
define the Lie algebra of outer semi-special derivations 
${\rm ODer}^{SS}(\overline L_{{\cal X}})$. There is 
a canonical surjective morphism 
$$
p_{x_0}: {\rm Der}^{S}L_{S} \quad \lra \quad {\rm ODer}^{SS}(\overline L_{{\cal X}})
$$
with one dimensional kernel generated by the inner derivation given by 
commutator with $\sum_{x \in {\cal X}} X_x$.

Denote by $\overline {\cal C}(A_{S}) $ the quotient of 
 the Lie algebra 
$\widetilde {\cal C}(A_{S})$ of 
cyclic words in 
$A_{{\cal X}}$ by the 
subspace generated by 
the elements $(\sum_{s \in S}X_s)^n$, $n\geq 1$.

\begin{theorem} \label{8.17.00.2} a) The space 
$\widetilde {\cal C}(I_{{\cal X}})$ is a Lie algebra 
 ideal in $\widetilde {\cal C}(A_{{\cal X}})$.

b) The element $(\sum_{s \in S}X_s)^n$ 
is in the center of the Lie algebra 
$\widetilde {\cal C}(A_{S})$. 

c) Let us choose an element $x_0 \in {\cal X}$. Then there is a canonical isomorphism of Lie algebras 
\begin{equation} \label{8.23.00.5} 
\frac{\widetilde {\cal C}(A_{{\cal X}})}{ \widetilde {\cal C}(I_{{\cal X}})} \quad 
\stackrel{=}{\lra} \quad \overline {\cal C}(A_{S}) \quad \stackrel{\overline \kappa }{= }
\quad {\rm ODer}^{SS}(\overline A_{{\cal X}})
\end{equation}
\end{theorem}

{\bf Proof}. Take cyclic words $A = 
{\cal C}\Bigl((\sum_{x \in {\cal X}}X_x) 
\cdot A_1\Bigr) \in \widetilde {\cal C}(I_{{\cal X}})$ and $B \in 
\widetilde {\cal C}(A_{{\cal X}})$. Then $[A,B]$, 
by the very definition of the commutator 
in the Lie algebra of cyclic words, 
is a sum of the terms 
corresponding to pairs (a generator in $A$, a generator in $B$). 
The pairs (a generator in $A_1$, a generator in $B$) clearly 
produce an element of the subspace $I_{{\cal X}}$. It is easy 
to see that the sum of the terms corresponding to the remaining pairs 
($\sum_{x \in {\cal X}} X_x$, a generator in $B$) is zero: this is very similar 
to the proof of the fact that a cyclic word provide a derivation killing $\sum_{s \in S}X_s$, 
 and the same kind of argument proves b). 
So $\widetilde {\cal C}(I_{{\cal X}})$ 
is a Lie ideal in $\widetilde {\cal C}(A_{{\cal X}})$, and $\overline {\cal C}(A_{{S}})$ 
is a Lie algebra. 

The isomorphism (\ref{8.17.00.3} ) 
shows that one obviously has an isomorphism of vector spaces
$$
\frac{{\cal C}(A_{{\cal X}})}{ {\cal C}(I_{{\cal X}})} \quad 
\stackrel{=}{\lra} \quad {\cal C}(A_{S}) 
$$
Passing to the $\widetilde {\cal C}$-quotients on the left we get an isomorphism 
of vector spaces
$$
\frac{\widetilde {\cal C}(A_{{\cal X}})}{ \widetilde {\cal C}(I_{{\cal X}})} \quad 
\stackrel{=}{\lra} \quad \overline {\cal C}(A_{S}) 
$$
Notice that the element $X^n_{x_0}$ is zero in $\widetilde {\cal C}(A_{{\cal X}})$ but it is not zero in $\widetilde {\cal C}(A_{{S}})$. So we had to define the quotient 
$\overline {\cal C}(A_{{S}})$. 
Moreover, since any element of $\overline A_{{\cal X}}$ can be written modulo 
$I({\cal X} )$ 
as element of 
$A_{{S}}$ (i.e. we use again (\ref{8.17.00.3})) 
this isomorphism commutes with the Lie brackets.

The second equality in (\ref{8.23.00.5}) follows from 
the description of ${\rm Ker}p_{x_0}$ and proposition \ref{ssw}.  The theorem is proved. 

We define the group ${\rm Out}^{SSE}\pi_1^{(l)}(X_{N})$ 
of outer semi-special {\it equivariant} automorphisms 
of $\pi_1^{(l)}(X_{N})$ as the invariants of the action of $\mu_N$ on 
${\rm Out}^{SS}\pi_1^{(l)}(X_{N})$, and proceed similarly in the other cases. 

If ${\cal X}_G:= \{0\} \cup G \cup \{\infty\}$, $x_0 = \{\infty\}$, 
theorem \ref{8.17.00.2} provides a description of the Lie algebra ${\rm ODer}^{SSE}A(G)$ 
of outer semi-special {\it equivariant} derivations of $A(G)$:
\begin{equation} \label{8.24.00.1} 
\Bigl(\frac{\widetilde {\cal C}(A_{{\cal X}_G})}
{ \widetilde {\cal C}(I_{{\cal X}_G})}\Bigr)^G \quad 
\stackrel{=}{\lra} \quad {\rm ODer}^{SSE}(\overline A_{{\cal X}_G})
\end{equation}

Combining this with theorem \ref{ga}b) we get a description of the Lie algebra 
${\rm ODer}^{SSE}(\overline L_{{\cal X}_G}) \subset 
{\rm ODer}^{SSE}(\overline A_{{\cal X}_G})$. 

{\bf 7. The distribution relations.} 
Recall that  $X_N = {\Bbb P}^1 \backslash \{0, \mu_N, \infty\}$ and  there is canonical isomorphism ${\rm Gr}^W{\Bbb L}(X_{M}, v_{\infty}) = L(\mu_N)$. 
The maps 
$$
i_N: X_{NM} \hookrightarrow  X_M, \quad z \lms z, \qquad \qquad m_N: X_{NM} \lra X_M, \quad z \lms z^N
$$
induce  the surjective Lie algebra homomorphisms 
$$
i_{N*}: { L}(\mu_{MN}) \lra {L}(\mu_{M}); \qquad m_{N*}: {L}(\mu_{MN}) \lra {L}(\mu_{M}) 
$$
given on  the generators by
$$
i_{N*}: Y \lra Y, \quad X_{\zeta} 
\lra \left\{ \begin{array}{ll}
0 &  \zeta \not \in \mu_M \\ 
 X_{\zeta}  &   \zeta \in \mu_M \end{array} \right.
$$
$$
m_{N*}: Y \lra NY,  \quad X_{\zeta} \lra X_{\zeta^N} 
$$

\begin{lemma} Elements of ${\rm ODer}^{SSE}{L}(\mu_{MN})$ preserve
${\rm Ker} (i_{N*})$ and  ${\rm Ker} (m_{N*})$. 
\end{lemma}

{\bf Proof}.  ${\rm Ker} (i_{N*})$ 
is generated by the elements $ X_{\zeta}$ where $\zeta \not \in \mu_N$. Since 
 any semi-special 
derivation preserves conjugacy class of $X_{\zeta}$, it in particular preserves 
the ideal generated by such $X_{\zeta}$. 

The statement about $m_N$ follows 
from the fact that equivariant derivations, by their very definition, commute 
with the natural action of the group $G$ on ${L}(\mu_{G})$. 
Namely,  let ${\cal D} \in {\rm Der}^E{L}(\mu_{MN})$. 
Then $\xi_*({\cal D}(X_{\zeta})) = {\cal D}(X_{\xi\zeta})$, and so $m_{N*}$ sends 
${\cal D}(X_{\zeta})$ and ${\cal D}(X_{\xi\zeta})$ to the same element. The lemma is proved.

It follows that  there are  well defined homomorphisms 
\begin{equation} \label{8.26.00.31}
{\widetilde i}_{N}, {\widetilde m}_{N}: 
{\rm ODer}^{SSE}_{}{L}(\mu_{MN}) \quad \lra \quad 
{\rm ODer}^{SSE}_{}{L}(\mu_{M})
\end{equation}
uniquely characterized by the property that 
for any element 
$ l \in {L}(\mu_{MN})$ and a derivation 
$ D \in {\rm ODer}^{SSE}_{}{ L}(\mu_{MN})$ one has 
\begin{equation} \label{8.26.00.3}
f_* \circ D(l) \quad = \quad {\widetilde f} (D)\circ f (l); \qquad f = i_N, m_N
\end{equation}

The same arguments provide us with  homomorphisms 
\begin{equation} \label{8.26.00.1}
{\overline i}_{N}, {\overline m}_{N}: 
{\rm Out}^{SSE}{L}^{(l)}(\mu_{MN}) \quad \lra \quad
{\rm Out}^{SSE}{L}^{(l)}(\mu_{M})
\end{equation}
Recall the maps 
$$
\psi^{(l)}_N: {\rm Gal}\Bigl(\overline \Q/ \Q(\zeta_{l^{\infty}N}) \Bigl)
\quad \lra \quad {\rm Out}^{SSE}{L}^{(l)}(\mu_N)
$$
They together with homomorphisms (\ref{8.26.00.1}) provide  the following two 
 commutative diagrams (similar to (\ref{8.26.00.3})), where $j_M$ is the natural inclusion:
$$
\begin{array}{ccc}
{\rm Gal}\Bigl(\overline \Q/ \Q(\zeta_{l^{\infty}MN}) \Bigr)& 
\stackrel{\psi^{(l)}_{MN}}{\lra} & 
{\rm Out}^{SSE} {L}^{(l)}(\mu_{MN})\\
\downarrow j_M& &\downarrow  \overline i_N\\
{\rm Gal}\Bigl(\overline \Q/ \Q(\zeta_{l^{\infty}M})\Bigr) & 
\stackrel{\psi^{(l)}_{M}}{\lra} & {\rm Out}^{SSE} {L}^{(l)}(\mu_{M})
\end{array}
$$
and 
$$
\begin{array}{ccc}
{\rm Gal}\Bigl(\overline \Q/ \Q(\zeta_{l^{\infty}MN}) \Bigr)& 
\stackrel{\psi^{(l)}_{MN}}{\lra} & 
{\rm Out}^{SSE} {L}^{(l)}(\mu_{MN})\\
\downarrow j_M& &\downarrow  \overline m_N\\
{\rm Gal}\Bigl(\overline \Q/ \Q(\zeta_{l^{\infty}M})\Bigr) & 
\stackrel{\psi^{(l)}_{M}}{\lra} & {\rm Out}^{SSE}{L}^{(l)}(\mu_{M})
\end{array}
$$
The composition $\psi^{(l)}_{M} \circ j_M$ does not depend 
on the choice of the maps 
$\overline i_N, \overline m_N$. 
Passing to the Lie 
algebras we get
\begin{equation} \label{9.9.00.111}
{\rm Gr}^W{\cal G}_{NM}^{(l)} \quad \subset \quad
{\rm Ker}({\widetilde i}_{N} - {\widetilde m}_{N}) \otimes \Q_l 
\end{equation}
The maps $i_{N*}, m_{N*}$  provide the  linear maps 
$$
i'_{N*}, m'_{N*}: {\cal C}(A(\mu_{MN}))^{\mu_{MN}} \quad \lra \quad 
{\cal C}(A(\mu_{M}))^{\mu_{N}}
$$
We set ${i}^{''}_{N}:= {i}'_{N*}$ and ${m}^{''}_{N}:= N^{-1}{m}_{N*}$.

\begin{lemma}
The maps ${i}^{''}_{N}, {m}^{''}_{N}$ preserve the Lie brackets. 
\end{lemma}

{\bf Proof}. Direct check. 

{\it Warning}. The map ${m}^{''}_{N}$ is a  Lie algebra morphism 
only when restricted to the subspace of $G$-invariant cyclic words. 

Since  
\begin{equation} \label{9.8.00.111}
{i}^{''}_{N} - {m}^{''}_{N}: \quad Y^k \lms (1-N^{k-1}) Y^k, \quad 
\sum_{g \in \mu_{MN}}X_g^k \lms 0
\end{equation}
the identification ${\rm Der}^{SE}{A}(\mu_{G}) = 
\widetilde {\cal C}(A(\mu_{G}))^G $ leads to the Lie algebra  maps 
$$%\begin{equation} \label{9.8.00.1}
{\rm Der}^{SE}{A}(\mu_{MN}) \quad \lra \quad 
{\rm Der}^{SE}{A}(\mu_{M})
$$%\end{equation}
Then they restrict  to  the  maps
\begin{equation} \label{9.8.00.1}
{i}'_{N}, {m}'_{N}: 
{\rm Der}^{SE}{L}(\mu_{MN}) \quad \lra \quad 
{\rm Der}^{SE}{L}(\mu_{M})
\end{equation}

\begin{proposition} \label{9.9.00.1}
The distribution relations are equivalent to the conditions
\begin{equation} \label{8.23.00.1}
{\rm Gr}^W{\cal G}_{N}^{(l)} \quad \subset \quad
{\rm Ker}({i}'_{L} - {m}'_{L}) \qquad \mbox{for each integer $L|N$}
\end{equation}
\end{proposition}

{\bf Proof}. 
An element of ${\rm Der}^{SE}{ L}(\mu_{N}) $ is given by expression 
\begin{equation} \label{9.7.00.1}
\sum a(g_0:g_1: ... : g_m)_{n_0, .., n_m} \cdot 
\kappa_{\mu_{N}}{\cal C}(X_{g_0} Y^{n_0-1} ... X_{g_m} Y^{n_m-1}) 
\end{equation}
 The sum is over cyclic words. The $a$'s are coefficients in $\Q$. They  satisfy 
the relation   
$a(g_0: ... : g_m)_{*} = a(hg_0: ... : hg_m)_{*}$ expressing the fact that the cyclic word (\ref{9.7.00.1}) is G-invariant.  
It is in $Ker({i}'_{L} - {m}'_{L})$ if and only if 
\begin{equation} \label{9.7.00.2}
a(g_0: g_1: ... : g_m)_{n_0, .., n_m} = L^{-1}L^{w-m}\sum_{h_i^L = g_i}a(h_0: h_1: ... : h_m)_{n_0, .., n_m}
\end{equation}
except 
the relation 
\begin{equation} \label{9.9.00.10}
a(g:g) = L^{-1}\sum_{h_i^L = g}a(h_0: h_1)
\end{equation}
 for  $m=1$ and 
$g_0 = g_1$.  
This is precisely the distribution relations.  
Here is the origin of the exception. 
An element 
$$
\sum a(g_0:g_1: ... : g_m)_{n_0, .., n_m} \cdot 
{\cal C}(X_{g_0} Y^{n_0-1} ... X_{g_m} Y^{n_m-1}) 
$$
 is  in 
$Ker({i}^{''}_{L} - {m}^{''}_{L})$ if and only if 
 equations (\ref{9.7.00.2}) are valid. Then to get 
the maps (\ref{9.8.00.1}) we first restrict to the subspace of those 
cyclic words which 
provide derivations of the 
Lie algebra $L(\mu_N)$, and then kill the elements 
by $Y^2$ and $\sum_{g \in G}X_g^2$. Therefore, thanks to (\ref{9.8.00.111}), 
the kernel is increased by a one 
dimensional subspace:  we add the element projected by 
${i}^{''}_{L} - {m}^{''}_{L}$ onto 
$\sum_{g \in \mu_{N/L}}X_g^2$. The equation we thus removed 
is exactly (\ref{9.9.00.10}). 
Indeed, this equation just means that the $\sum_{g \in \mu_{N/L}}X_g^2$ 
component of $({i}^{''}_{L} - {m}^{''}_{L})(\ref{9.7.00.1})$ is zero, which we no longer require.  
The proposition is proved. 

One has 
$$
{i}'_{L}- { m}'_{L}: \quad (Y + \sum_{g \in \mu_N} X_g)^2 \quad \lms \quad 
(1-L)(Y + \sum _{g \in \mu_{N/L}}X_g)^2
$$
So  killing the elements $(Y + \sum X_g)^2$ we get well defined maps
$$
{\widetilde i}'_{L}, {\widetilde  m}'_{L}: 
{\rm ODer}^{SSE}{L}(\mu_{N}) \quad \lra \quad 
{\rm ODer}^{SSE}{ L}(\mu_{N/L})
$$
and moreover, since $1 - L \not = 0$, we have  ${\rm Ker}({i}'_{L}- { m}'_{L}) = {\rm Ker}
({\widetilde i}'_{L}- {\widetilde m}'_{L})$. 

\begin{lemma} \label{9.9.00.11}
One has ${\widetilde i}_{L} = {\widetilde i}'_{L}$ and ${\widetilde m}_{L} = 
{\widetilde m}'_{L}$
\end{lemma}

{\bf Proof}. One checks that each of the maps satisfy 
condition 
(\ref{8.26.00.3}). Since this condition characterizes these maps 
the lemma follows. 

Now the distribution relations 
follow from(\ref{9.9.00.111}),  proposition \ref{9.9.00.1}  and 
lemma \ref{9.9.00.11}.

{\bf Proof of theorem \ref{8.17.00.1}}. It follows immediately from the inversion relation, which just has been proved as $N=-1$ case of 
the distribution relations. The argument is identical with
 the one given in the proof of corollary \ref{8.26.00.100}.

{\bf 8. Conjecture \ref{ramierz}}. 
To complete the proof of conjecture \ref{ramierz} it remains to 
prove that ${\rm Gr}{\cal G}_{\bullet \bullet}^{(l)}(\mu_N)$ lies in the subspace of 
${\rm Gr}{\rm Der}^{SE}_{\bullet \bullet}{L}(\mu_{N}) $ 
defined by the power shuffle relations (\ref{gw4}). 
So these relations  provide the most nontrivial constraints on 
the image of the Galois group. This  is rather amazing since they are 
the simplest relations on the level of functions. 
Unfortunately I do not know a good motivic proof of them. 
There exists, however, an indirect argument: relations (\ref{gw4}) 
hold for functions $=>$ for the corresponding framed mixed Hodge structures $=>$ for the corresponding mixed Tate motives over $\Q(\zeta_N)$ $=>$ valid for their $l$-adic realization. The 
detailed exposition will appear in [G4]. In s. 7.3-7.4 I  prove  
 certain special cases of conjecture \ref{ramierz} relevant to our situation 
which can be obtained by some ad hoc methods.

\section{Cohomology of some discrete subgroups of $GL_2(\Z)$ and $GL_3(\Z)$}

{\bf 1. The general scheme.} Let $G$ be a reductive group over $\Q$. Denote by ${\Bbb A}_f$ the ring of 
finite adels for $\Q$. Choose a finite index subgroup $K_f \subset G(\widehat \Z)$. 
Let $Z\subset G$ be the maximal split torus in the center of $G$, and $Z^0(\R)$ the 
connected component of identity of its $\R$-valued points. Denote by $K_{\infty}'$ 
the connected component of the maximal compact subgroup of the derived group 
$G'(\R)$. Set
$$
K_{\infty}^{{\rm min}}:= K'_{\infty} \times Z^0(\R)
$$
$$
K_{\infty}^{{\rm max}}:= \mbox{maximal compact subgroup of $G(\R) \times Z^0(\R)$}
$$
Choose a subgroup $K_{\infty}$  sitting in between:
$$
K_{\infty}^{{\rm min}} \subset K_{\infty} \subset K_{\infty}^{{\rm max}}
$$
One defines the modular variety 
corresponding 
to a given choice of the subgroups $K_f$ and $K_{\infty}$ as follows:
\begin{equation} \label{6.6.00.1}
S^G_{K_{\infty} \times K_f} = G(\Q) \backslash \Bigl(G(\R)/K_{\infty} 
\times G({\Bbb A}_f)/K_f \Bigr)
\end{equation}
In general it is, of course, not a variety in the algebraic geometry sense, 
and not even a manifold, only an orbifold.  
A rational representation $
\rho: G\otimes_{\Q}\overline \Q \lra GL(V)
$ provides a local system on the orbifold 
(\ref{6.6.00.1}): 
$$
{\cal L}_V:= \quad  \Bigl(G(\R)/K_{\infty} 
\times G({\Bbb A}_f)/K_f \Bigr)\times_{G(\Q)} V 
$$

{\bf Example}. $G = GL_m$, $K^{{\rm max}}_{\infty} = O_m \cdot \R^*_+$, 
$K^{{\rm min}}_{\infty} = SO_m \cdot \R^*_+$. The space 
$GL_m(\R)/K^{{\rm max}}_{\infty}$ is identified with positive definite 
symmetric $m \times m$ matrices. The group $GL_m(\R)$ acts on it by $X \lms AXA^t$.  
The symmetric space is
$$
{\Bbb H}_m  = GL_m(\R)/K^{{\rm max}}_{\infty}\cdot \R^*_{>0} \quad = \quad 
\{\mbox{$>0$ definite symmetric $m \times m$ matrices}\}/\R^*_+
$$
One defines 
$$
\Gamma := \quad GL_m(\Q) \cap K_f \quad \subset \quad GL_m(\Z)
$$
If $\Gamma$ is torsion free then 
$$
H^*(\Gamma, V) \quad = \quad H^*(S^G_{K_{\infty} \times K_f}, {\cal L}_V)
$$
If $\Gamma$  contains a finite index torsion free subgroup 
$\widetilde \Gamma$ then
$
H^*(\Gamma, V)  = H^*(\widetilde \Gamma, V)^{\Gamma/\widetilde \Gamma } 
$. 

We will use the shorthand $S_{\Gamma}$ for the modular variety corresponding to the subgroup
$\Gamma$. 
The Borel-Serre compactification $\overline S_{\Gamma}$ 
is a compact manifold with corners. The boundary $\partial 
\overline S_{\Gamma}$ is a topological manifold. It is a disjoint union of faces which 
correspond bijectively to the $\Gamma$-conjugacy classes of proper parabolic 
subgroups of $G$ defined over $\Q$. The closure of the face corresponding to a parabolic subgroup $P$ is  called a 
strata and denoted $\partial_{P}S$. 
There is a natural restriction map
$$
H^*(S_{\Gamma}, {\cal L}_V) \stackrel{{\rm Res}}{\lra} H^*(\partial \overline S_{\Gamma}
, {\cal L}_V) 
$$
The cohomology at infinity
$H_{{\rm inf}}^*(\partial \overline S_{\Gamma}, {\cal L}_V)$ are, by definition,
  the image of this map.   
The computation of these groups is an important step in 
understanding of the cohomology of ${\cal L}_V$. To carry it out one should determine first the cohomology 
of the restriction of ${\cal L}_V$ to the Borel-Serre boundary $\partial 
\overline S_{\Gamma}$.

Pick a rational parabolic subgroup $P$. Let $U_P$ be its unipotent radical and $M_P:= P/U_P$ 
the Levi quotient.  Denote by ${\cal N}_P$ the Lie algebra of $U_P$. Take the subgroup $P(\R) \cap K_{\infty}$ 
and project it 
to $K_{\infty}^M \subset M(\R)$. 
In our case 
$
P({\Bbb A}_f) \cdot K_f = G({\Bbb A}_f)
$.  
Then the cohomology $H^*(\partial_PS, {\cal L}_V)$ of the 
restriction of the local 
system ${\cal L}_V$ to the boundary strata
 $\partial_PS$  are calculated using the Leray spectral sequence whose $E_2$-term looks as follows:
\begin{equation} \label{6.7.00.1}
H^p(S^M_{K^M_{\infty}}, H^q({N}_P(\Z), V)) \quad => \quad H^{p+q}(\partial_PS, {\cal L}_V)
\end{equation} 
One has 
 $$
H^{q}(N_P(\Z), V) = H^{q}({\cal N}_P, V)
$$
The first step of the calculation of the cohomology at infinity of a subgroup $\Gamma$ is 
the calculation of  these groups for all of the strata.

{\it The Kostant theorem}. Let $P$ be a rational parabolic subgroup of a reductive group $G$.  
Fixing Cartan and Borel subgroups  
$H(\C) \subset B(\C)$ with $H(\C)\subset M_P(\C)$ we get a system of positive roots and the Weyl group $W$. Denote by $W_P$ the Weyl group for $M_P$ with $W_P \subset W$. 
It is known that in each coset class $W_P \backslash W$ there is a unique
 element of minimal possible length. Denote by $W_P^1$ the set of 
the minimal length representatives for $W_P \backslash W$. 
For any dominant weight $\mu$ denote by $L_{\mu}$ the irreducible representation of the group $M_P(\C)$ with the highest weight $\mu$. 
Let $\rho$ be the half of the sum of the 
positive roots for $G$. 

\begin{theorem} \label{8-29.1/99} {\rm ([K], Theorem 5.14).} The $M_P(\C)$-representation $H^*({\cal N}_P, V)$ is algebraic and is given by
$$
H^*({\cal N}_P, V) \quad = \quad \oplus_{\omega \in W^1_P}
L_{\omega(\lambda+\rho) - \rho}[-l(\omega)]
$$
where $l(\omega)$ is the length of $\omega$, and $[-l(\omega)]$ indicates that the module appears in degree $l(\omega)$. If $V$ is defined over $\Q$, then so is $H^*({\cal N}_P, V)$. 
\end{theorem}

Let $H $ be the
 maximal Cartan subgroup of $GL_m$ given by the diagonal matrices.  
Denote by $[n_1, ..., n_m]$ the character of $H$ given by 
$(t_1, ..., t_m) \lra t_1^{n_1} \cdot ... \cdot t_m^{n_m}$. 
Let $P$ be a parabolic 
subgroup $P$ containing $H$ and $M_P$ the Levi quotient of $P$. 
 Then $H$ is a maximal torus for $M_P$. Denote by $V^P_{[n_1, ..., n_m]}$ or simply 
$V_{[n_1, ..., n_m]}$ the representation with the highest weight 
$[n_1, ..., n_m]$ of $M_P$.

{\bf Example}. The weights $\rho$ for $GL_2$ and $GL_3$ are given by
$$
\rho = [1/2, -1/2]; \qquad \rho = [1, 0, -1]; 
$$

{\bf 2. Cohomology of $GL_2(\Z)$}. 
The cohomology $H^{*}(GL_2(\Z), S^{w-2}V_2\otimes \varepsilon_2)$ vanish for odd $w$. 
Indeed, if $w$ is odd the central element ${\rm diag}(-1,-1)$ acts by $-1$ in $S^{w-2}V_2$. 
So we will assume that $w-2$ is even. Clearly 
$$
H^{0}(GL_2(\Z), S^{w-2}V_2\otimes \varepsilon_2) = 0 
$$
For any $GL_2$-module $V$ over $\Q$  the Hochshild-Serre
 spectral sequence corresponding to the exact sequence 
$$
0 \lra SL_2(\Z) \lra GL_2(\Z) \stackrel{{\rm det}}{\lra} \{\pm 1\} \lra 0
$$
gives 
$$
H^*(GL_2(\Z), V) = H^*(SL_2(\Z), V)^+ 
$$
Here $+$ is the invariants under the involution provided  by conjugation 
by the matrix $\left (\matrix{-1& 0\cr 0& 1\cr}\right )$.  
There is an exact sequence ([Sh])
$$
0 \lra H_{{\rm cusp}}^{1}(SL_2(\Z), S^{w-2}V_2) \lra H^{1}(SL_2(\Z), S^{w-2}V_2) 
\stackrel{{\rm Res}}{\lra} 
$$
$$
H_{{\rm inf}}^{1}(SL_2(\Z), S^{w-2}V_2) \lra 0
$$
It admits a natural splitting as a module over the Hecke algebra. 

Let us compute the cohomology of the restriction of the local system 
${\cal L}_{S^{w-2}V_2}$ 
 to the boundary strata $\partial_BS$ corresponding to the 
upper triangular Borel subgroup $B \subset GL_2$. In this case $M_B = T$ is the diagonal torus and $K^M_{\infty} = \left (\matrix{\xi_1& 0\cr 0& \xi_2\cr}\right )$, 
$\xi_i = \pm 1$. 
We need to compute 
$$
H^p(K^M_{\infty}, H^q({\cal N}_B, S^{w-2}V_2 \otimes \varepsilon_2)
$$
We have $S^{w-2}V_2\otimes \varepsilon_2 = L_{[w-1,1]}$ and 
$$
H^q({\cal N}_B, S^{w-2}V_2\otimes \varepsilon_2 ) = \quad \left\{ \begin{array}{llll}
L_{[w-1, 1]} &  q=0 \\ 
 L_{[0, w]}  &  q=1\end{array} \right.
$$
Therefore $H^0 =0$ (since $w-1$ is odd), and 
$$
H^1(GL_2(\Z), S^{w-2}V_2 \otimes \varepsilon_2) = H_{{\rm cusp}}^1(SL_2(\Z), S^{w-2}V_2)^-  \oplus 
H_{{\rm inf}}^1(SL_2(\Z), S^{w-2}V_2)
$$
It is known that the $+$ and $-$ parts of $H_{{\rm cusp}}^{1}(SL_2(\Z), S^{w-2}V_2)$ 
are of the same dimension which, thanks to the Eichler-Shimura isomorphism,  
coincides with the dimension of the space of the holomorphic weight $w$ cusp 
forms on $SL_2(\Z)$.

{\bf 2. Cohomology of $GL_3(\Z)$ with coefficients in $S^{w-3}V_3$}. 
The main result is 
\begin{theorem} \label{8-28/99} 
\begin{equation} \label{8-28.1/99}
 H^{i}(GL_3(\Z), 
S^{w-3}V_3)  \quad = \quad \left\{ \begin{array}{lll}
\Q & i=0, w=3\\
 H_{{\rm cusp}}^{1}(GL_2(\Z),S^{w-2}V_3\otimes \varepsilon_2) &  i=3 \\
0& {\rm otherwise}  \end{array} \right.
\end{equation}
\end{theorem}

{\bf Proof}. Since the central element $-{\rm Id} \in GL_3$ acts by $-1$ on $V_3$, and hence on $S^{w-3}V_3$ if $w-3$ is odd, we get  
$
H^{i}(GL_3(\Z), 
S^{w-3}V_3)  =0$ if $w-3$ is odd. We will assume from now on that $w-3$ is even. 

For the definition of  cuspidal cohomology 
$H^*_{{\rm cusp }}(\Gamma, V)$ which we use below see s. 1.3 of [LS].  

\begin{proposition} \label{6.8.00.1}
a) The following  sequence is exact
$$
0 \lra H^i_{{\rm cusp}}(GL_3(\Z), 
S^{w-3}V_3) \lra H^i(GL_3(\Z), 
S^{w-3}V_3) \lra 
$$
$$
\lra H_{{\rm inf}}^i(GL_3(\Z), 
S^{w-3}V_3) \lra 0
$$

b) $H^*_{{\rm cusp}}(GL_3(\Z), S^{w-3}V_3) = 0$.
\end{proposition}

{\bf Proof}. Let $sl_3$ be the Lie algebra of $SL_3(\R)$. Let $\pi$ be a unitary irreducible representation  of $SL_3(\R)$.  Denote by $H^{\infty}_{\pi}$ 
the space of $C^{\infty}$-vectors in $\pi$. Since for $w>3$ $S^{w-3}V_3$ is not self dual 
(i.e. the Cartan involution $\theta: g \lms (g^{-1})^t$ transforms it to a non isomorphic 
representation)   one has according to proposition 6.12 in chapter II of [BW] 
$$
H^{*}(sl_3, SO_3; S^{w-3}V_3\otimes H^{\infty}_{\pi}) = 0
$$
The proposition follows from the argumentation given in chapter  1 of [LS].

The boundary of the Borel-Serre compactification of 
our modular variety is a disjoint union of three faces, $e(P_1), e(P_2)$ and $e(B)$ which 
correspond to the $GL_3(\Q)$-equivalence classes of the following proper 
parabolic $\Q$-subgroups:
$$
P_1:= \left (\matrix{*& *&*\cr 
0&  *&*\cr
0&  *&*\cr}\right ); \qquad P_2:= \left (\matrix{*& *&*\cr 
* &  *&*\cr
0&  0&*\cr}\right ) \qquad B:= \left (\matrix{*& *&*\cr 
0 &  *&*\cr
0&  0 &*\cr}\right )
$$
Their nilradicals are denoted by $N_1$, $N_2$ and $N$; the corresponding Lie algebras are ${\cal N}_1$, ${\cal N}_2$ and ${\cal N}$. We are going to calculate the groups 
(\ref{6.7.00.1}). 

Applying Kostant's theorem we get the following: 
$$
H^q({\cal N}_1, S^{w-3}V_3) = \quad \left\{ \begin{array}{lll}
L_{[w-3, 0, 0]} &  q=0 \\ 
 L_{[-1, w-2, 0]} &  q=1\\
 L_{[-2, w-2, 1]} &  q=2\end{array} \right.
$$

$$
H^q({\cal N}_2, S^{w-3}V_3) = \quad \left\{ \begin{array}{lll}
L_{[w-3, 0, 0]} &  q=0 \\ 
 L_{[w-3, -1, 1]} &  q=1\\
 L_{[-1, -1, w-1]} &  q=2\end{array} \right.
$$

$$
H^q({\cal N}, S^{w-3}V_3) = \quad \left\{ \begin{array}{llll}
L_{[w-3, 0, 0]} &  q=0 \\ 
 L_{[-1, w-2, 0]} \oplus L_{[w-3, -1, 1]} &  q=1\\
 L_{[-1, -1, w-1]} \oplus L_{[-2, w-2, 1]} &  q=2\\
L_{[-2, 0, w-1]} &  q=3\end{array} \right.
$$
Let 
$$
S_1:= \left (\matrix{\pm 1& 0\cr 0& GL_2(\Z)\cr }\right ); \quad 
S_2:= \left (\matrix{GL_2(\Z)& 0\cr 0& \pm 1\cr }\right ); \quad
S:= \left (\matrix{\pm 1& 0&0\cr 0& \pm 1&0 \cr 0&0&\pm 1}\right )
$$ 

If there is a $\pm 1$ subgroup in one of these groups, $S_{?}$, 
 which acts by $-1$ 
on a module $V$, then $H^i(S_{?}, V)=0$ for all $i$. 
This simple remark shows that we get 
 a nonzero contribution only the following cases:

i) For the group $S_1$: 
${L}_{[w-3, 0, 0]} $ and ${ L}_{[-2, w-2, 1]} [-2]$.

ii) For the group $S_2$: ${L}_{[w-3, 0, 0]} $ and ${L}_{[-1, -1, w-1]} [-2]$.

iii) For the group $S$: 
${L}_{[w-3, 0, 0]} $ and ${ L}_{[-2, 0, w-1]} [-3]$. 

We consider them now case by case, indicating the value $(p,q)$ for the 
$E^{p, q}_2$-term of the spectral sequence where they appear. 
 Recall that  $w-3$ is even. 

i) For the strata $\partial_{P_1}S$ we have the following:
\begin{equation} \label{6.7.00.2}
H^0(S_1, L_{[w-3, 0, 0]}) = L_{[w-3, 0, 0]}; \qquad (p,q)=(0,0)
\end{equation}
and 
$$
H^1(S_1, L_{[-2,w-2,1]}) \quad = \quad 
H^1_{{}}(GL_2(\Z), S^{w-3}V_3 \otimes \varepsilon_2) \quad = 
$$
\begin{equation} \label{6.7.00.3}
H^1_{{\rm cusp}}(SL_2(\Z), S^{w-3}V_2)^- \quad \oplus \quad 
H^1_{{\rm inf}}(SL_2(\Z), S^{w-3}V_2);\qquad (p,q) = (1,2)
\end{equation}
The other cohomology groups are zero.

ii) Further, for the strata $\partial_{P_2}S$
$$
H^0(S_2, L_{[w-3, 0, 0]}) \quad = \left\{ \begin{array}{lll}
\Q &  w=3 \\ 
 0 &  w>3 \end{array} \right. ; \qquad (p,q) = (0,0)
$$
$$
H^1(S_2, L_{[w-3, 0, 0]}) \quad = \quad H^1(GL_2(\Z), S^{w-3}V_2) \quad = 
$$
$$
H^1_{{\rm cusp}}(SL_2(\Z), S^{w-3}V_2)^+;
\qquad (p,q) = (1,0)
$$
The other cohomology groups vanish. 

iii) Finally, for the strata $\partial_{P_2}S$ the only nonzero cohomology groups are  
$$
H^0(S, H^q({\cal N}, S^{w-3}V_3)) \quad = \left\{ \begin{array}{lll}
\Q &  q=0 \\ 
 \Q &  q=3\\
 0 &  \mbox{otherwise}\end{array} \right.
$$

Since 
$$
\partial \overline S = \partial_{P_1} S \cup \partial_{P_2} S; \qquad 
\partial_{B} S = \partial_{P_1} S \cap \partial_{P_2} S
$$
 we have the Mayer-Vietoris long exact sequence
$$
... \lra  H^*(\partial \overline S, {\cal L}_V) \lra H^*(\partial_{P_1} S, {\cal L}_V) \oplus 
H^*(\partial_{P_2} S, {\cal L}_V) \lra  
H^*(\partial_{B} S, {\cal L}_V) \lra  ...
$$
It is easy to see that the differential in this complex maps isomorphically the group 
(\ref{6.7.00.2}) to  the group from the step iii) with $q=0$. 

Further it follows from the definitions that the differential
 maps isomorphically the group $H^1_{{\rm inf}}$ from (\ref{6.7.00.3}) onto the 
group from the step iii) with $q=3$. 

Therefore the boundary cohomology are
$$
H^i(\partial \overline S, {\cal L}_{S^{w-3}V_3}) =  \left\{ \begin{array}{lll}
 H^1_{{\rm cusp}}(SL_2(\Z), S^{w-3}V_2)^+&  i=1 \\ 
  H^1_{{\rm cusp}}(SL_2(\Z), S^{w-3}V_2)^-&  i=3
\end{array} \right.
$$
To determine the subgroup $H_{{\rm inf}}^i(\partial \overline S, {\cal L}_{S^{w-3}V_3})$ one usually uses the theory of Eisenstein cohomology classes originated by G. Harder. 
See [H6] for details in our case. However just 
in our case there is another approach. Namely, we will show in the next chapter that 
$
H^3(GL_3(\Z), {\cal L}_{S^{w-3}V_3}) 
$
contains the subgroup $H^1_{{\rm cusp}}(SL_2(\Z), S^{w-3}V_2)^-$, 
(see (\ref{dvas})). 
The proof uses theorem 1.3 in [G3] and some simple arguments presented in the section 7.1 which 
are independent of the rest of the paper. Therefore 
\begin{equation} \label{6.8.00.11}
H^3_{{\rm inf}}(\partial \overline S, {\cal L}_{S^{w-3}V_3}) = 
H^1_{{\rm cusp}}(SL_2(\Z), S^{w-3}V_2)^-
\end{equation}
It remains to show that $H_{{\rm inf}}^1$ is zero. 
The Poincar\'e duality provides a nondegenerate pairing 
\begin{equation} \label{6.8.00.10}
H^1(\partial \overline S, {\cal L}_{S^{w-3}V_3}) \otimes 
H^3(\partial \overline S, {\cal L}_{S^{w-3}V^{\vee}_3}) \lra \Q
\end{equation}
Employing 
the Cartan involution $ \theta$ 
of the group $GL_3$ which transforms $V_3$ into $V^{\vee}_3$ we get
\begin{equation} \label{6.9.00.11}
H_{{\rm inf}}^3(\partial \overline S, {\cal L}_{S^{w-3}V^{\vee}_3}) = 
H^1_{{\rm cusp}}(SL_2(\Z), S^{w-3}V_2)^-
\end{equation}
If $x \in H^1(S, {\cal L}_{S^{w-3}V_3})$ and $y \in 
H^3( S, {\cal L}_{S^{w-3}V^{\vee}_3})$ 
then ${\rm Res} (x) \cup {\rm Res} (y) =0$. 
Therefore (\ref{6.8.00.11}) together with the 
Poincar\'e pairing 
(\ref{6.8.00.10}) imply that 
${\rm Res} (x) =0$. The theorem is proved.

{\bf 3. Cohomology of $\Gamma_1(3;p)$}. 
Recall that if $\Gamma$ is a torsion free subgroup of $SL_3(\Z)$ then 
$$
H^*(\Gamma, \Q) \quad = \quad H^*(S_{\Gamma}, \Q)
\quad = \quad \quad H^*(\overline S_{\Gamma}, \Q) 
$$ 
Lee and Schwermer derived in the section 1.8 of [LS] an exact sequence
$$
0 \lra H^*_{{\rm cusp }}(\Gamma, \C) \lra H^*
(\overline S_{\Gamma}, \C)\lra 
H^*_{{\rm inf}}(\overline S_{\Gamma}, \C)\lra 0
$$
One has $H^q(S_{\Gamma}, \Q) =0$ for $q>3$. Further, $H^{1}(S_{\Gamma}, \Q) =0$ 
by ch. 16 in [BMS]. 
The cuspidal cohomology satisfy the Poincar\'e duality. Therefore
\begin{equation} \label{6.9.00.5}
 {\rm dim }H_{{\rm cusp}}^{2}(\Gamma_1(3;p),\Q) \quad = \quad 
{\rm dim }H_{{\rm cusp}}^{3}(\Gamma_1(3;p),\Q)
\end{equation}

\begin{theorem} \label{6.19.00} Let $p$ be a prime number. Then 
\begin{equation} \label{2.6.00.1} 
H^i_{{\rm inf}}(\Gamma_1 (3; p), \Q)\quad  = \left\{ \begin{array}{ll}
\Q &  q=0  \\ 
 0 &    q=1, q \geq 4\\
H^1_{{\rm cusp}}(\Gamma (2; p), V_2) & q=2 \\
H^1_{{\rm cusp}}(\Gamma (2; p), \varepsilon_2) \otimes \Q^2 
 \quad \oplus \quad \Q^{\frac{p-3}{2}} & q=3 \end{array} \right.
 \end{equation}
\end{theorem}

{\bf Proof}. 
Let $\Gamma(m;p)$ be the full congruence subgroup of level $p$ of $GL_m(\Z)$. 
The cohomology of the subgroup $\Gamma(3;p)$ has been computed by Lee and Schwermer [LS]. 
Let me recall (and slightly correct) their result. 

Set $SL^{\pm}_m(F_p):= GL_m(\Z)/\Gamma(m;p)$. Here 
$F_p$ is the finite field of order  $p$.  Then 
$H^*_{{\rm inf}}(\overline S_{\Gamma(3m;p)}, \Q)$ is a 
$SL^{\pm}_m(F_p)$-module,

Consider the following subgroups in $SL^{\pm}_3(F_p)$:
$$
\widetilde P_1:= \quad \left (\matrix{*&*&*\cr
*&*&*\cr
0&0&\pm 1\cr}\right ) \quad 
\widetilde M_1:= \quad \left (\matrix{*&*&0\cr
*&*&0\cr
0&0&\pm 1\cr}\right ) 
\quad 
$$
$$
\widetilde P_2:= \quad \left (\matrix{\pm 1&*&*\cr
0&*&*\cr
0&*&*\cr}\right ) \qquad \widetilde B:= \quad \left (\matrix{\pm 1&*&*\cr
0&\pm 1&*\cr
0&0&\pm 1\cr}\right ) 
$$
Denote by $ P_1'$, $M'_1$ and $P'_2$ 
the  similar subgroups of $SL_3(F_p)$. 
Let $\widetilde \Gamma_1(m; p)$ be the 
projection of  $\Gamma_1 (m; p)$ onto $SL^{\pm}_m(F_p)$. 
Denote by $\xi'_i$ (resp. $\widetilde \xi_i$) 
the nontrivial one dimensional representation of $P'_i$ (resp. $\widetilde P_i$)  
pulled from  the one of the $GL_2$-part of the Levi quotient given by the determinant. 
Denote   by [V] the isomorphism class of a representation $V$. 
Recall that  $V_2$ is the standard two dimensional $GL_2$-module. 

Let  $\widetilde {\rm St_3}$ be the generalized Steinberg representation of 
$SL^{\pm}_3(F_p)$; it sits in the exact sequence
$$
0 \lra \Q \lra \oplus_{i=1}^2{\rm Ind}_{\widetilde P_i}^{SL_3^{\pm}(F_p)} \Q \lra 
{\rm Ind}_{\widetilde B}^{SL_3^{\pm}(F_p)} \Q \lra \widetilde {\rm St_3}\lra 0
$$

It was proved in [LS], see s.2.5 there, that,  considered as a $SL_3(F_p)$-module, 
$$
[H_{{\rm inf}}^q(\overline S_{\Gamma(3;p)}, \Q)] = 
$$
\begin{equation} \label{2.6.00.1} 
\quad  = \left\{ \begin{array}{ll}
\Q &  q=0  \\ 
 0 &    q=1, q \geq 4\\
2[{\rm Ind}_{P'_1}^{SL_3(F_p)}\xi'_1] \quad \oplus \quad [{\rm Ind}_{P'_1}^{SL_3(F_p)}
H^1_{{\rm cusp}}(\Gamma (2; p), V_2)] & q=2 \\
\oplus_{i=1}^2[{\rm Ind}_{P'_i}^{SL_3(F_p)}H^1_{{\rm cusp}}(\Gamma (2; p), \varepsilon_2)]
 \quad \oplus \quad [\widetilde {\rm St}_3] & q=3 \end{array} \right.
 \end{equation}

In  fact for $q=3$ in [LS] appears  the module
$
\oplus_{i=1}^2[{\rm Ind}_{P'_i}^{SL_3(F_p)}H^1_{{\rm cusp}}(\Gamma (2; p), \Q)]
$, 
but the $P'_i$-module $[H^1_{{\rm cusp}}(\Gamma (2; p), \Q)]$ should be changed to 
$[H^1_{{\rm cusp}}(\Gamma (2; p), \varepsilon_2)]$. As abelian groups they are isomorphic. 
See also an argument in s. 7.7 below which shows that $\Q$ has to be changed to 
$\varepsilon_2$. 

One can check, following the arguments in [LS], that, as a $SL^{\pm}_3(F_p)$-module, 
$H^q_{{\rm inf}}(\overline S_{\Gamma(3;p)}, \Q)$ in the nontrivial cases $q=2,3$ 
looks as follows: 
$$
[H_{{\rm inf}}^q(\overline S_{\Gamma(3;p)}, \Q)] = 
$$
\begin{equation} \label{6.9.00.1} 
\quad  = \left\{ \begin{array}{ll}
2[{\rm Ind}_{\widetilde P_1}^{SL^{\pm}_3(F_p)}\xi_1] \quad \oplus 
\quad [{\rm Ind}_{\widetilde P_1}^{SL^{\pm}_3(F_p)}
H^1_{{\rm cusp}}(\Gamma (2; p), V_2)]; & q=2 \\
\oplus_{i=1}^2[{\rm Ind}_{\widetilde P_i}^{SL^{\pm}_3(F_p)}H^1_{{\rm cusp}}(\Gamma (2; p), \varepsilon_2)]
 \quad \oplus \quad [\widetilde {\rm St}_3] & q=3 \end{array} \right.
 \end{equation}
(We replaced $SL_3$ by $SL^{\pm}_3$ and put $'$ by  instead 
of tilde over $P_i$ and $\xi_i$). 

To  compute $H^*_{{\rm inf}}(\Gamma_1(3;p), \Q))$ 
we need to find the invariants of the action of 
$\widetilde \Gamma_1(3;p) $ on these $SL^{\pm}_3(F_p)$-modules. 
We will use the following general statement. 

\begin{lemma} \label{KIR}
Let 
 $\rho: \widetilde P_1 \lra Aut(W)$ be a representation trivial on the unipotent radical and $\pm 1$ subgroup, and such that $H^0(SL_2^{\pm}(F_p), W) =0$. Then 
$$
\Bigl({\rm Ind}_{\widetilde P_1}^{SL^{\pm}_3(F_p)}
W\Bigr)^{\widetilde \Gamma_1(3; p) } = W^{\widetilde \Gamma_1(2;p) }
$$
\end{lemma}

{\bf Proof}. The left hand side  is given by the space of $W$-valued 
functions $f(g)$ on $SL^{\pm}_3(F_p)$ satisfying the condition
$$
f(\gamma g y) = \rho(y^{-1})\cdot f(g) , \qquad \gamma \in \widetilde \Gamma_1(3; p), \quad 
y \in \widetilde P_1. 
$$ 
The coset $\widetilde \Gamma_1(3; p) \backslash SL^{\pm}_3(F_p)$ is identified with the 
set of 
rows $X = (x_1, x_2, x_3) \in F^3_p - 0$. The group $\widetilde P_1$ acts on it from the right. 
So we need to determine 
the space of the $W$-valued 
functions $f(X)$ 
such that $f(Xy) = \rho(y^{-1})f(X)$. 
 Any non zero vector $(x_1, x_2, x_3)$ is $\widetilde  \Gamma_1(3; p)$-equivalent to 
$(0,1,0)$ or $(0,0,x)$. The stabilizer of $(0,1,0)$ is $\widetilde 
\Gamma_1(2;p)\times \{\pm 1\}$. So 
 $f(0,1,0) \in W^{\widetilde \Gamma_1(2;p)}$. The stabilizer of $(0,0,x)$ 
is $SL_2^{\pm}(F_p)\times \{\pm 1\}$. Since $H^0(SL_2^{\pm}(F_p), W) =0$ we have $f(0,0,x) =0$.
The lemma is proved. 

Applying the lemma we get the following spaces of $\widetilde \Gamma_1 (3; p)$-invariants: 

i) 
For $\rm Ind_{\widetilde P_1}^{SL^{\pm}_3(F_p)}\widetilde \xi_1$ it is zero.

ii) For 
${\rm Ind}_{\widetilde P_1}^{SL^{\pm}_3(F_p)}
H^1_{{\rm cusp}}(\Gamma (2; p), V_2)$ it is $H^1_{{\rm cusp}}(\Gamma_1 (2; p), V_2)$.

iii) For 
${\rm Ind}_{\widetilde P_i}^{SL^{\pm}_3(F_p)}
H^1_{{\rm cusp}}(\Gamma (2; p), \varepsilon_2)$ it is $H^1_{{\rm cusp}}(\Gamma_1 (2; p), \varepsilon_2)$.

iv) The dimension of the space of $\widetilde \Gamma_1 (3; p)$-invariants 
on  $\widetilde {\rm St_3}$ is $\frac{p-3}{2}$. 

{\it Argumentation}. i) Indeed,  $\xi(g) = -1$ for $ g:= {\rm diag}(-1, 1, 1) \in \widetilde \Gamma_1 (3; p)
$. 

ii) and iii) are clear.

iv) It is a corollary of the following elementary statements, 
which are checked similarly to the proof of the lemma above:
$$
{\rm dim}\Bigl( {\rm Ind}_{\widetilde B}^{SL_3^{\pm}(F_p)} \Q 
\Bigr)^{\widetilde \Gamma_1 (3; p)} = 3 \cdot \frac{p-1}{2}
$$
$$
{\rm dim}\Bigl( {\rm Ind}_{\widetilde P_i}^{SL_3^{\pm}(F_p)} \Q 
\Bigr)^{\widetilde \Gamma_1 (3; p)} = \frac{p-1}{2}+1
$$
 and the trivial module $\Q$ is, of course, $\widetilde \Gamma_1 (3; p)$-invariant.
So we get
$$
3 \cdot \frac{p-1}{2} - 2 (\frac{p-1}{2} + 1) -1  = \frac{p-3}{2}
$$
The theorem is proved.

\section{Proofs of the theorems from section 2}

{\bf 1. The Soul\'e elements}. Since ${\rm Gr}{\cal G}^{(l)}_{\bullet, -1}$ is abelian we may identify it with the $\Q_l$-points of the corresponding algebraic group. So it makes sense to talk about projection of $\varphi^{(l)} (\sigma)$ on 
 ${\rm Gr}{\cal G}^{(l)}_{\bullet, -1}$. Let $\varphi_{-w, -1}^{(l)} (\sigma)$ be the component of this projection in 
${\rm Gr}{\cal G}^{(l)}_{-w, -1}$.

In [So] Soul\'e constructed for each integer $m >1$ a ${\rm Gal}(\Q(\zeta_{l^{\infty}})/\Q)$-homomorphism
$$
\chi_m : {\rm Gal}(\overline \Q/\Q(\zeta_{l^{\infty}}))^{ab} \lra \Q_l(m)
$$
 He proved that it is zero if and only if $m$ is even. 
Let $I_{2m-1}(e:e)^{\vee}$ be the generator of 
${D}_{-2m+1, -1}$ dual to $I_{2m-1}(e:e )$. It follows from the Key Lemma B in  [Ih2] that 
$$
\varphi^{(l)} (\sigma) = \frac{(1-l^{m-1})^{-1}}{(m-1)!}\chi_m(\sigma) \xi \Bigl(I_{2m-1}(e:e)^{\vee} \Bigr)  \qquad (m: {\rm odd} \geq 3)
$$

{\bf 2. The depth 1 case}. The distribution relations, proved in s. 4.7, in the depth one case just mean that 
$$
{\rm Gr} 
{\cal G}^{(l)}_{-w, -1} (\mu_N) \quad \subset \quad \xi_{\mu_N}\Bigl(D_{-w, -1 }(\mu_N)\Bigr) \otimes \Q_l
$$
Easy classical arguments combined with Borel's theorem show that 
$$
{\rm dim} D_{-w, -1 }(\mu_N) \quad = \quad 
{\rm dim} {\rm Hom}(K_{2w-1}(S_N), \Q)
$$
For $w=1$ it is a reformulation of the Bass theorem, 
and the general case is completely similar.  
From the motivic theory of classical 
polylogarithms at $N$-th roots of unity we get the following (see theorem \ref{4-12.60090}):
$$
{\rm dim} {\rm Gr} 
{\cal G}^{(l)}_{-w, -1} (\mu_N) \quad = \quad {\rm dim} {\rm Hom}(K_{2w-1}(S_N), \Q)
$$
Indeed, it has been proved in [BD], see also  [HW], that 
motivic classical polylogarithms at $N$-th roots of unity provide 
classes generating $K_{2w-1}(S_N)\otimes \Q$; in particular the $l$-adic 
regulators of these classes has been computed. 
Denote by $H_1^{(-w)}$ the part of degree $-w$ in $H_1$, and similarly 
$H_1^{(-w, -m)}$ the part of degree $-w$, depth $-m$ in $H_1$. 
We conclude that 
\begin{equation} \label{9.8.00.2}
H_1^{(-w)}\Bigl({\cal G}^{(l)}_{N} /
{\cal F}_{ -2}{\cal G}^{(l)}_{N}  \Bigr)   \quad = \quad   {\rm Gr} 
{\cal G}^{(l)}_{-w, -1} (\mu_N) \quad = 
\end{equation}
\begin{equation} \label{9.10.00.3}
{\rm Hom}(K_{2w-1}(S_N), \Q_l)\quad = \quad 
\xi_{\mu_N}\Bigl(D_{-w, -1 }(\mu_N)\Bigr) \otimes \Q_l
\end{equation}

{\bf 3. Description of the image of the Galois group: the depth two case}. 
\begin{theorem} \label{tbav} a) Conjecture \ref{ramierz} is valid in the depth $2$ case, i.e. 
$$
{\rm Gr}{\cal G}^{(l)}_{\bullet, \geq -2}(\mu_N) \quad \hookrightarrow \quad 
\xi_{\mu_N}\Bigl(D_{\bullet, \geq -2}(\mu_N)\Bigr) \otimes \Q_l
$$

b) Moreover  
$$
{\rm Gr}{\cal G}^{(l)}_{\bullet, -2 }(\mu_N) \quad = \quad 
\xi_{\mu_N}\left(\frac{{\cal D}_{\bullet, -2 }(\mu_N)}{{\rm Ker}\Bigl(\delta : {\cal D}_{\bullet, -2 }(\mu_N) \lra \Lambda^2 {\cal D}_{\bullet, -1 }(\mu_N)\Bigr)}\right)^{\vee}
$$
\end{theorem}

{\bf Proof}. The part b) implies the part a). 
Since ${\rm Gr}^W{\cal G}^{(l)}_{N}$ 
is a quotient of the fundamental Lie algebra $L_{{\cal T}\Q_l}(S_N)$ (see s. 2.7) one has 
\begin{equation} \label{9.8.00.1}
H_1^{(-w)}\Bigl({\cal G}^{(l)}_{N} /
{\cal F}_{-3}{\cal G}^{(l)}_{N}  \Bigr) \quad \subset 
\quad {\rm Hom}(K_{2w-1}(S_N), \Q_l)
\end{equation}
Therefore 
\begin{equation} \label{8.10.00.2}
H_1^{(-w, -2)}\Bigl({\cal G}^{(l)}_N
/{\cal F}_{-3}{\cal G}^{(l)}_{N}  \Bigr)  = 0
\end{equation}
Indeed, $[{\cal G}^{(l)}_N, {\cal G}^{(l)}_N] \subset 
{\cal F}_{-2}{\cal G}^{(l)}_N$, so a nontrivial depth $-2$ 
part of $H_1$ plus (\ref{9.8.00.2}) 
will make the left hand side of (\ref{9.8.00.1}) bigger then the right hand side. 

If ${\cal G}$ is a Lie algebra with a filtration ${\cal F}_{\bullet}$ indexed 
by integers $n = -1, -2, ...$, such that 
${\cal F}_{-1}{\cal G} = {\cal G}$, then the Lie algebra
$$
{\cal G}/ {\cal F}_{ -3}{\cal G} \quad 
\mbox{ is isomorphic to the Lie algebra} \quad 
{\rm Gr}^{{\cal F}}_{\geq -2}{\cal G}:= 
\quad {\rm Gr}^{{\cal F}}_{-1}{\cal G} \oplus {\rm Gr}^{{\cal F}}_{-2}{\cal G}
$$  
In particular 
${\cal G}^{(l)}_N /
{\cal F}_{-3}{\cal G}^{(l)}_N $ is isomorphic to  ${\rm Gr}{\cal G}^{(l)}_{\bullet, \geq -2}(\mu_N) $. Therefore thanks to (\ref{8.10.00.2}) we have 
\begin{equation} \label{9.10.00.4}
H_1^{(-w, -2)}
\Bigl({\rm Gr}{\cal G}^{(l)}_{\bullet, \geq -2}(\mu_N) \Bigr)  = 0
\end{equation}
This means that ${\rm Gr}{\cal G}^{(l)}_{\bullet, \leq -2 }(\mu_N)$ 
is generated by ${\rm Gr}{\cal G}^{(l)}_{\bullet, -1 }(\mu_N)$.

 Therefore
$$
{\rm Gr}{\cal G}^{(l)}_{\bullet, -2 }(\mu_N) \quad \stackrel{(\ref{9.10.00.4})}{= } \quad
[{\rm Gr}{\cal G}^{(l)}_{\bullet, -1 }(\mu_N), {\rm Gr}{\cal G}^{(l)}_{\bullet, -1 }(\mu_N)]  \quad \stackrel{(\ref{9.10.00.3})}{= }
$$
$$
\xi_{\mu_N}([D_{\bullet, -1 }(\mu_N), D_{\bullet, -1 }(\mu_N)]) \quad \subset \quad
\xi_{\mu_N}(D_{\bullet, -2 }(\mu_N))
$$
This is equivalent to the statement of the theorem. 
The theorem is proved. 

{\bf Remark}. If $N=1$, or $N =p$ is a prime and $\bullet = -2$,  
we prove below that the space
$$
{\rm Ker}\Bigl(\delta : {\cal D}_{\bullet, -2 }(\mu_N) \lra \Lambda^2 {\cal D}_{\bullet, -1 }(\mu_N)\Bigr) 
$$
is zero. It may not be zero otherwise. It would be 
interesting to construct explicitly the elements in the kernel.

Using theorem \ref{tbav} one can show that if $N$ is not prime then 
 $\xi_{\mu_N}({D}^{\Delta}_{\bullet}(\mu_N))$ could be bigger then ${\rm Gr}{\cal G}^{(l)}_{\bullet}(\mu_N)$.

{\bf 4. The image of the Galois group for $m=3, N=1$ case}

\begin{theorem} \label{tbav10} 
$$
{\rm Gr}{\cal G}^{(l)}_{-w, -3}\quad \hookrightarrow \quad 
\xi(D_{-w, -3}) \otimes \Q_l
$$
\end{theorem}

{\bf Proof}. When $w$ is odd this follows from theorem \ref{8.17.00.1}. So we may assume 
 that $w$ is even. 
The same argumentation as in the proof of theorem \ref{tbav} gives 
\begin{equation} \label{8.10.00.12}
H_1^{(-w, -3)}\Bigl({\cal G}^{(l)}_N
/{\cal F}_{-4}{\cal G}^{(l)}_{N}  \Bigr)  = 0
\end{equation}
However {\it a priori} the Lie algebra ${\cal G}^{(l)}_N
/{\cal F}_{-4}{\cal G}^{(l)}_{N}$ may  be different from 
${\rm Gr}{\cal G}_{\bullet, \geq -3}^{(l)}(\mu_N)$.  Recall that 
${\rm Gr}^W{\cal G}_{N}^{(l)}$ is non canonically isomorphic to 
${\cal G}_{N}^{(l)}$. 
We are going to show that if $N=1$ and  $w$ is odd the weight $w$, depth $\leq 3$ 
parts 
of the standard cochain complexes of 
${\rm Gr}^W{\cal G}^{(l)}$ and ${\rm Gr}{\cal G}_{\bullet \bullet}^{(l)}$ are isomorphic. For the depth $\leq 2$ parts this has been proved above. 
So the descrepency between them can appear only if there are 
elements $x\in {\cal F}_3{\rm Gr}^W({\cal G}^{(l)})^{\vee}$ such that 
$$
\delta (x) \in  {\cal F}_1{\rm Gr}^W({\cal G}^{(l)})^{\vee} \wedge 
{\cal F}_1{\rm Gr}^W({\cal G}^{(l)})^{\vee} 
$$ 
Since  ${\cal F}_1{\rm Gr}_p^W({\cal G}^{(l)})^{\vee}$ 
can be non zero only if $p$ is odd, the weight of $\delta(x)$ 
must be even. This contradicts to the assumption that $w$ is odd. 
Therefore 
\begin{equation} \label{9.10.00.5}
H_1^{(-w, -3)}
\Bigl({\rm Gr}{\cal G}^{(l)}_{\bullet, \geq -3} \Bigr)  \quad = \quad 
H_1^{(-w, -3)}\Bigl({\cal G}^{(l)}_N
/{\cal F}_{-4}{\cal G}^{(l)}_{N}  \Bigr)  \quad = \quad 0
\end{equation}
Therefore this together with the previous theorem implies that when $w$ is odd 
${\rm Gr}{\cal G}^{(l)}_{-w, \geq -3}$ is generated by 
iterated 
commutators of triples of elements of ${\rm Gr}{\cal G}^{(l)}_{\bullet, -1}$. 
Since  ${\rm Gr}{\cal G}^{(l)}_{\bullet, -1} = 
\xi (D_{\bullet, -1})\otimes \Q_l$ we conclude that 
${\rm Gr}{\cal G}^{(l)}_{-w, \geq -3}$ coincides with the 
subspace generated by iterated commutators of triples of elements of 
$\xi (D_{\bullet, -1})\otimes \Q_l$, and so belongs to 
$\xi (D_{\bullet, -3})\otimes \Q_l$. The theorem is proved.

Let ${\cal L}^{(l)}_{\bullet, -k} \subset 
{\rm Gr}{\cal G}^{(l)}_{\bullet, -k}$ be the subspace generated by iterated 
commutators of $k$ elements of 
${\rm Gr}{\cal G}^{(l)}_{\bullet, -1}$. Then  
$
{\cal L}^{(l)}_{\bullet, \geq -m}:=  \oplus_{k=1}^m {\cal L}^{(l)}_{\bullet, -k}
$
is a bigraded Lie subalgebra of ${\rm Gr}{\cal G}^{(l)}_{\bullet, \geq -m}$. 

\begin{corollary} \label{2/04/00.4}
${\cal L}^{(l)}_{\bullet, -m} \quad\stackrel{}{=}\quad
{\rm Gr}{\cal G}^{(l)}_{\bullet, -m} \quad\quad \mbox{for $m=2,3$}
$.  
\end{corollary}

{\bf Remark}. 
Denote by $\varepsilon_2$ the one dimensional $GL_2$-module 
given by the determinant. 
We proved recently that, assuming conjecture \ref{ramierz} below, one has 
$$
{\rm Gr}{\cal G}^{(l)}_{-w, -4} /{\cal L}^{(l)}_{-w, -4}  = 
H^1_{{\rm cusp}}(GL_2(\Z), S^{w-2}V_2 \otimes \varepsilon_2) \otimes \Q_l
$$

{\bf 5. The cohomology of the Lie algebras $D_{\bullet \bullet}$ and $\widehat 
D_{\bullet \bullet}$ of depths $2, 3$}. 
Recall the  isomorphism of bigraded Lie algebras 
$$\widehat {D}_{\bullet \bullet}(G) = 
{D}_{\bullet \bullet}(G)  \oplus \Q(-1,-1)
$$
where $\Q(-1,-1)$ is a  one dimensional Lie algebra 
of the bidegree $(-1,-1)$. 
So there is the decomposition 
of the standard cochain complex of $\widehat {\cal D}_{\bullet \bullet}(G)$:
\begin{equation} \label{4-16.10}
\Lambda^* \widehat {\cal D}_{\bullet \bullet}(G) \qquad =\qquad 
\Lambda^* {\cal D}_{\bullet \bullet}(G)  \quad \oplus \quad
\Lambda^* {\cal D}_{\bullet \bullet}(G) \otimes \Q_{(1,1)}
\end{equation}
It provides a canonical morphism of complexes
$$
\partial_m: \Lambda_{(m)}^* \widehat {\cal D}_{\bullet \bullet}(G) \lra 
\Lambda_{(m-1)}^*  {\cal D}_{\bullet \bullet}(G)[-1]
$$
Surprisingly it is easier to describe the structure of the Lie algebra 
 $\widehat {D}_{\bullet \bullet}(\mu_N)$, although the Galois-theoretic or motivic meaning of 
its $\Q_{(-1,-1)}$-component is unclear: it should correspond to $\zeta(1)$, or, better, 
 Euler's $\gamma$-constant,  
whatever it means.

{\bf Remark}. i) Notice the classical formula
$$
d \log \Gamma (1-z) = \sum_{m\geq 1} \zeta(m) z^{m-1} dz; \qquad \mbox{where we set} \quad 
\zeta(1):= \gamma
$$
ii) One has 
$$
\gamma = - \int_0^{\infty} e^{-t}\log t dt = - \int_0^{\infty} \frac{dt}{t} \circ (e^{-t}  dt) 
$$  
(The iterated integral on the right has to be regularized). 
So the Euler $\gamma$-constant 
should be thought of as an ``irregular'' period of a motive.  

It was proved in theorem 1.3 in [G3], see also lemma \ref{modcoho},  that 
\begin{equation} \label{d1f.}
H^i_{(w, 2)}(\widehat {D}_{\bullet \bullet}) =  H^{i-1}(GL_2(\Z), 
S^{w-2}V_2 \otimes \varepsilon_2), \quad i=1,2
\end{equation}
\begin{equation} \label{d2f.}
H^i_{(w, 3)}(\widehat {D}_{\bullet,\bullet}) = H^{i}(GL_3(\Z), 
S^{w-3}V_3)  = 0, \quad i=1,2,3
 \end{equation}

\begin{theorem} \label{C} Denote by  $H_{(w, m)}$ the  weight $w$, 
depth $m$ part of $H$. Then 
  \begin{equation} \label{d1f}
H^i_{(w, 2)}({D}_{\bullet \bullet}) =  H_{{\rm cusp}}^{i-1}(GL_2(\Z), 
S^{w-2}V_2\otimes \varepsilon_2), \quad i=1,2
\end{equation}
\begin{equation} \label{d2f}
H^i_{(w, 3)}({D}_{\bullet,\bullet}) = 0, \quad i=1,2,3
 \end{equation}
\end{theorem}

{\bf Proof}.  Consider the inclusion of complexes 
provided by the depth $2$ part of (\ref{4-16.10}):
\begin{equation} \label{dwa}
\Bigl( {\cal D}_{\bullet, 2} \lra \Lambda^2{\cal D}_{\bullet, 1} \Bigr) 
\quad \hookrightarrow \quad 
{\cal D}_{\bullet, 2} \lra \Lambda^2 \widehat {\cal D}_{\bullet, 1} 
\stackrel{\widehat \partial_2}{\lra} \widehat 
{\cal D}_{\bullet, 1} 
\end{equation} 
Here  $\widehat \partial_2$ is the composition of the map
 $\partial_2$ followed by the natural inclusion $
{\cal D}_{\bullet, 1} \hookrightarrow \widehat 
{\cal D}_{\bullet, 1} $. 
The isomorphism ${\mu}$ from theorem 1.2 in [G3] is 
easily extended to an isomorphism between the complex 
${\Bbb M}_{(2)}^{\ast} \otimes_{GL_2(\Z)} S^{\bullet -2}V_2 $ and the complex on the right of 
 (\ref{dwa}),  
which takes $\tau_{[1,2]}{\Bbb M}_{(2)}^{\ast} \otimes_{GL_2(\Z)} S^{\bullet -2}V_2 $ 
just to the subcomplex on the left. This together with lemma \ref{modcoho} prove  (\ref{d1f}). 

Notice that the second summand 
${\cal D}_{\bullet, 1}[-1]$ in the decomposition (\ref{4-16.10}) of the depth $2$ part 
of $\Lambda^*\widehat {\cal D}_{\bullet \bullet}$  provides the 
Eisenstein part of the cohomology.

The depth $3$ part of the decomposition  (\ref{4-16.10})
  is 
\begin{equation} \label{raz}
\Bigl({\cal D}_{\bullet, 3} \stackrel{\delta}{\lra}  
{\cal D}_{\bullet, 2}\otimes {\cal D}_{\bullet, 1}   \stackrel{\delta}{\lra}  
\Lambda^3{\cal D}_{\bullet, 1}\Bigr) \quad \oplus \quad 
\Bigl({\cal D}_{\bullet, 2}   \stackrel{\delta}{\lra}  
\Lambda^2{\cal D}_{\bullet, 1} \Bigr)
\end{equation}
Therefore we see that 
\begin{equation} \label{dvas}
H^3_{(w, 3)}(\widehat {D}_{\bullet \bullet}) \quad \mbox{contains as a direct summand} \quad 
H_{{\rm cusp}}^{1}(GL_2(\Z), 
S^{w-2}V_2\otimes \varepsilon_2)
\end{equation}
This was the last bit needed to prove theorem \ref{8-28/99}. 

Now we can start using theorem \ref{8-28/99}.
By theorem 1.3 in [G3] and theorem \ref{8-28/99}  we have
  $$
 H^i_{(w, 3)}(\widehat {D}_{\bullet \bullet})\quad   \stackrel{1.3 }{= } \quad H^{i}(GL_3(\Z),S^{w-3}V_3)\quad 
\stackrel{\ref{8-28/99}}{= }
$$
\begin{equation} \label{d2f+}
\left\{ \begin{array}{ll}
0 &  i=1,2 \\ 
 H_{{\rm cusp}}^{1}(GL_2(\Z),S^{w-2}V_2\otimes \varepsilon_2) &  i=3 \end{array} \right.
\end{equation} 
Therefore using (\ref{d1f}) we see that the complex on the left of (\ref{raz}) is acyclic. 
This proves formula 
(\ref{d2f}). 
Theorem \ref{C} is proved.

{\bf 6.  Proofs of theorems  \ref{4-13.1}, \ref{mth1t}, and \ref{dwhn1}}. 
Set 
${D}^{(l)}_{\bullet \bullet} := {D}_{\bullet \bullet} \otimes \Q_l$. 
Applying to $\{e:e|t_0:t_1\} = \sum_{n>0} I_n(e:e) (t_0 - t_1)^{n-1} $ 
the distribution relation (\ref{inversion}) we 
get ${\cal D}_{-2w, -1} =0$. 
Therefore 
$$
{\cal G}^{(l)}_{\bullet, -1} = {D}^{(l)}_{\bullet, -1} 
$$  
Since 
$
{D}^{(l)}_{\bullet \bullet}$ and 
${\cal L}^{(l)}_{\bullet \bullet}$ are Lie subalgebras of 
${\rm Gr}{\rm Der }^{SE}_{\bullet \bullet}{\Bbb L}^{(l)}$ this implies that  
\begin{equation} \label{ini}
{\cal L}^{(l)}_{\bullet \bullet} \hookrightarrow {D}^{(l)}_{\bullet \bullet}
\end{equation}
 Thus 
there are commutative diagrams where $[,]$ are the commutator maps:
$$
\begin{array}{cccccccc}
\Lambda^2 {\cal L}^{(l)}_{\bullet, -1} & \stackrel{=}{\lra } & 
\Lambda^2 {D}^{(l)}_{\bullet, -1}  
&&&{\cal L}^{(l)}_{\bullet, -2} \otimes {\cal L}^{(l)}_{\bullet, -1} 
& \stackrel{}{\lra } & {D}^{(l)}_{\bullet, -2} \otimes 
{D}^{(l)}_{\bullet, -1} \\
\downarrow [,] & & \downarrow [,]   &&\mbox{and} &  \downarrow [,] && \downarrow [,] \\
{\cal L}^{(l)}_{\bullet, -2} & \hookrightarrow &{D}^{(l)}_{\bullet, -2} &&&{\cal L}^{(l)}_{\bullet, -3} & \hookrightarrow &{D}^{(l)}_{\bullet, -3} 
\end{array}
$$

 \begin{lemma}
$
H^1_{(\bullet, m)}(D_{\bullet \bullet}) =0 \quad \mbox{for} \quad m=2, 3
$. 
\end{lemma}

{\bf Proof}. 
Since $H^{0}(GL_2(\Z), 
S^{w-2}V_2\otimes \varepsilon_2) = 0
$ 
the lemma follows from formulas (\ref{d1f}) and (\ref{d2f}). 

The lemma means that the right vertical arrows in the  diagrams above are surjective. 
Therefore the left diagram shows that the map
$$
{\cal L}^{(l)}_{\bullet, -2} \hookrightarrow {D}^{(l)}_{\bullet, -2} \quad 
$$ 
is surjective, and hence is an isomorphism. 
Thus the top arrow in the right diagram is an isomorphism.  
Therefore the second diagram shows that the map
$$
{\cal L}^{(l)}_{\bullet, -3} \hookrightarrow {D}^{(l)}_{\bullet, -3} 
$$
is also an isomorphism. Summarizing we have
\begin{equation} \label{ini1}
{\cal L}^{(l)}_{\bullet, -2} =  {D}^{(l)}_{\bullet, -2} \quad 
\mbox{and } \quad {\cal L}^{(l)}_{\bullet, -3} =  {D}^{(l)}_{\bullet, -3} 
\end{equation}
Combining  this with corollary \ref{2/04/00.4} we get 
\begin{equation} \label{ini1}
{\cal G}^{(l)}_{\bullet, -2} =  {D}^{(l)}_{\bullet, -2} \quad 
\mbox{and } \quad {\cal G}^{(l)}_{-w, -3} =  {D}^{(l)}_{-w, -3}, \quad 
\mbox{$w$ is odd}
\end{equation}
%On the other hand one has 
%$$
%{\cal L}^{(l)}_{\bullet \bullet} \quad \hookrightarrow \quad
% {\cal G}^{(l)}_{\bullet \bullet} 
%\stackrel{{\rm conj. \ref{ramierz} }}{\hookrightarrow} 
%{D}^{(l)}_{\bullet \bullet}
%$$
%where the second inclusion is provided by conjecture \ref{ramierz}. 

After this we reduced the study of the depth $-2$ and $-3$ quotients of the Lie algebras 
${\cal G}^{(l)}_{\bullet \bullet}$ and  
 $\widehat {\cal G}^{(l)}_{\bullet \bullet}$ 
to the study of the corresponding quotients of the Lie algebras $D_{\bullet \bullet}$ and 
$\widehat D_{\bullet \bullet}$, which were investigated in [G3]. Therefore 
theorems \ref{4-13.1}a), \ref{mth1t}a), b) and \ref{dwhn1} 
follow from theorems 1.2 in [G3], lemma \ref{modcoho}, proposition \ref{6.8.00.1}b) and theorem 6.2 in [G3]. It remains to prove the formula from theorem 
\ref{mth1t}c). Thanks to theorem \ref{8.17.00.1} we may  assume that $w$ is odd. 

{\it Computation of the dimensions}. i) Consider the generating series
$$
a_m(x):= \sum_{w>0} {\rm dim}\Lambda^m_{(w)}({\cal D}_{\bullet, 1})x^w
$$
We know that $a_1(x) = \frac{x^{3 }}{1-x^{2}}$. Thus one has 
$$
a_m(x) = \frac{x^{3 + 5+ ... + (2m+1)}}{\prod_{i=1}^m(1-x^{2i})}
$$
In particular
$$
a_2(x) = \frac{x^{8 }}{(1-x^{2})(1-x^{4})}, 
\qquad a_3(x) = \frac{x^{15}}{(1-x^{2})(1-x^{4})(1-x^{6})}
$$

ii) The ring of modular forms for $SL_2(\Z)$ is generated by the Eisenstein 
series $E_4(z)$ and $E_6(z)$. By the Eichler-Shimura theorem the dimension of 
the space of modular forms for $SL_2(\Z)$ coincides with 
${\rm dim}H^1(GL_2(\Z), S^{w-2}V_2)$. Therefore
$$
\sum_{w\geq 2} {\rm dim}H^1(GL_2(\Z), S^{w-2}V_2) x^w \quad = \quad \frac{1}
{(1-x^{4})(1-x^{6})}-1
$$

iii) Let 
$$
d_m(x) := \sum_{w>0} {\rm dim}{\cal D}_{w,m} x^w
$$
Computing the Euler characteristic of the depth $2$ part of the standard cochain complex of the Lie coalgebra  ${\cal D}_{\bullet \bullet}$ and using theorem \ref{C} we get
$$
d_2(x) - a_2(x) \quad = \quad  - \frac{1}{(1-x^{4})(1-x^{6})} +1
$$
Therefore using the formulas above we get
$$
d_2(x) \quad = \quad \frac{x^8}{(1-x^{2})(1-x^{6})}
$$
This formula is equivalent to theorems \ref{4-13.1}b). 

iv) Computing the Euler characteristic of the depth $3$ part of the standard cochain complex of the Lie coalgebra  ${\cal D}_{\bullet \bullet}$ and using theorem \ref{C} we get
$$
d_3(x) - d_2(x) \cdot a_1(x) + a_3(x)\quad = \quad  0
$$
Therefore using the formulas above we get
$$
d_3(x) \quad = \quad \frac{x^{11}(1 + x^{2} - x^{4}) }{(1-x^{2})(1-x^{4})(1-x^{6})}
$$
This formula is equivalent to theorems \ref{mth1t}c). 

These results also provide the computation of the 
right hand sides of the formulas appearing in 
theorems 1.4 and 1.5 in [G3].

These theorems, in particular, imply that the double shuffle relations  provide 
a {\it complete} list 
of constraints on the Lie algebra of the image of the Galois group in 
${\rm Aut}\pi^{(l)}_1({\Bbb P}^1 - \{0, 1, \infty\}, v_{\infty})$ in the depths 
$-2$ and $-3$.

{\bf 7.  Proofs of theorems \ref{4-15.10}  and \ref{4-15.11}}. They are similar 
to the proofs in the subsection 3 above. Recall that $p$ is a prime number. 
We will use  shorthands like  ${D}^{(l)}_{\bullet}(\mu_p)$ for the diagonal Lie algebra  
${D}^{\Delta}_{\bullet}(\mu_p)^{(l)}$.
By theorem \ref{4-12.6} 
$$
{\rm Gr}{\cal L}^{(l)}_{-1}(\mu_p) = {D}^{(l)}_{ -1}(\mu_p) 
$$  
This implies that  
\begin{equation} \label{ini+}
{\cal L}^{(l)}_{\bullet}(\mu_p) \hookrightarrow {D}^{(l)}_{\bullet}(\mu_p)
\end{equation}
 Thus 
there are commutative diagrams
$$
\begin{array}{cccccccc}
\Lambda^2 {\cal L}^{(l)}_{ -1}(\mu_p) & \stackrel{=}{\lra } & 
\Lambda^2 {D}^{(l)}_{-1}(\mu_p)  \\
&&\\
\downarrow [,] & & \downarrow [,]   \\
&&\\
{\cal L}^{(l)}_{ -2}(\mu_p) & \hookrightarrow &{D}^{(l)}_{-2}(\mu_p)
\end{array}
$$
and
$$
\begin{array}{ccc}
{\cal L}^{(l)}_{ -2}(\mu_p) \otimes {\cal L}^{(l)}_{ -1}(\mu_p) 
& \stackrel{}{\lra } & {D}^{(l)}_{-2}(\mu_p)\otimes {D}^{(l)}_{-1}(\mu_p)\\
&&\\
  \downarrow [,] && \downarrow [,] \\
&&\\
{\cal L}^{(l)}_{ -3}(\mu_p)& \hookrightarrow &{D}^{(l)}_{-3}(\mu_p)
\end{array}
$$  
\begin{lemma} \label{4-17.1}
$
H^1_{(m)}(D_{\bullet}(\mu_p)) =0 \quad \mbox{for} \quad m=2, 3.
$
\end{lemma}

 {\bf Proof}.  
Thanks to decomposition (\ref{4-16.10}) one has 
$
H^1(\widehat D_{\bullet \bullet}(G))  = H^1(D_{\bullet \bullet}(G))
$. This implies   
$
H^1(\widehat D_{ \bullet}(G))  = H^1(D_{\bullet}(G))
$. 
Theorems 1.2, 6.1, 6.2  in [G3] provide  the following crucial result:
\begin{equation} \label{d1f_}
H^i_{(w, 2)}(\widehat {D}_{\bullet}(\mu_p)) \quad =  
\quad H^{i-1}(\Gamma_1(2;p), \varepsilon_2), \quad =  
\quad H^{i-1}(\Gamma_1(p), \Q)^-
\end{equation}
\begin{equation} \label{d1f_qeq}
H^i_{(w, 3)}(\widehat {D}_{\bullet}(\mu_p)) \quad =  
\quad H^{i}(\Gamma_1(3;p), \Q); \qquad i = 1,2,3
\end{equation}
In particular 
\begin{equation} \label{KRYAK}
H^1_{(2)}(\widehat D_{\bullet}(\mu_p)) = H^{0}(\Gamma_1(2;p), \varepsilon_2) 
\qquad H^1_{(3)}(\widehat D_{\bullet}(\mu_p)) = H^{1}(\Gamma_1(3; p), \Q)
\end{equation}
Both groups are zero: for the first one it is clear, and for the second see 
ch. 16 in [BMS] or use Kazhdan's theorem. The lemma is proved.

So the right vertical arrows in the  diagrams above are surjective. 
Using the same arguments as in s. 6.3 we get Lie algebra isomorphisms 
\begin{equation} \label{4-15.20}
{\cal L}^{(l)}_{ \geq -m}(\mu_p) \quad  = \quad {D}_{\geq  -m}(\mu_p)\otimes \Q_l, 
\qquad m = 2,3
\end{equation}

 After this theorems \ref{4-15.10} follows  
from  theorem \ref{tbav} the results of [G3]. 
Similarly theorem  \ref{4-15.11} follows  from (\ref{4-15.20}) and 
the results of [G3].

These results imply that the double shuffle relations  provide 
a {\it complete} list 
of constraints on the weight = depth part of the Lie algebra of the image of the Galois group in 
${\rm Aut}\pi^{(l)}_1(X_p, v_{ \infty})$ in the depths 
$-2$ and $-3$.

{\bf 8. The Lie coalgebra ${\cal D}^{{\rm un}}_{\bullet}(\mu_p)$}. Recall that $p$ is a prime number. 
Let us define a linear map $v_p: {\cal D}_{\bullet}(\mu_p) \lra \Q$ by 
$$
 v_p: 
{\cal D}_{m}(\mu_p) \lms 0 \quad \mbox{for $m>1$}, \quad v_p: I_{1,1}(1: \zeta_p^a) \lms 1, 
 \quad (\zeta_p^a \not = 1)
$$ 
Since there is  no distribution relations when $p$ is prime the map $v_p$ is   well defined. 
Let ${\cal D}^{{\rm un}}_{\bullet}(\mu_p) \hookrightarrow {\cal D}_{\bullet}(\mu_p)$ 
be the codimension $1$ 
subspace ${\rm Ker} (v_p)$ 
%$$
%{\cal D}^{{\rm un}}_{1}(\mu_p):= {\rm Ker}\Bigl(v_p: {\cal D}_{1}(\mu_p) \lra \Q\Bigr); \qquad 
%{\cal D}^{{\rm un}}_{m}(\mu_p) = {\cal D}_{m}(\mu_p) \quad \mbox{for $m>1$}
%$$

\begin{proposition} \label{6.1.00.100}
${\cal D}^{{\rm un}}_{\bullet}(\mu_p)$ is a sub Lie coalgebra 
of ${\cal D}_{\bullet}(\mu_p)$. 
\end{proposition}

{\bf Proof}. 
Let $V$ be a vector space and $f \in V^*$. There is a map 
$$
\partial_f: 
\Lambda^nV \lra \Lambda^{n-1}V; \qquad \partial^2 =0
$$
$$
v_{1} \wedge ... \wedge v_n \lms \sum_{i=1}^n(-1)^{i-1}f(v_i) v_{1} \wedge ... 
\wedge \widehat v_i \wedge ... \wedge v_n
$$
In particular the map $v_p$  provides a degree $-1$ maps  $$
 \partial_{v_p}: \Lambda^n  
{\cal D}_{\bullet}(\mu_p) \lra \Lambda^{n-1} {\cal D}_{\bullet}(\mu_p)
$$ 
There is an exact sequence
$$
0 \quad \lra \quad \Lambda^2_{(m)}{\cal D}^{{\rm un}}_{\bullet}(\mu_p)   \quad \hookrightarrow
 \quad \Lambda^2_{(m)}{\cal D}_{\bullet}(\mu_p) \quad \stackrel{\partial_{v_p}}{\lra} \quad 
{\cal D}^{{\rm un}}_{m-1}(\mu_p) \quad \lra \quad 0
$$
So to prove the proposition we need to check that the composition
$$
{\cal D}_{m}(\mu_p) \quad \stackrel{\delta}{\lra} \quad 
\Lambda^2_{(m)}{\cal D}_{\bullet}(\mu_p) \quad \stackrel{\partial_{v_p}}{\lra} \quad 
{\cal D}_{m-1}(\mu_p)
$$
is zero. To simplify notations we set $\{g_0: ... :g_m\}:= I_{1, ..., 1}(g_0: ... :g_m)$.
 Then 
\begin{equation} \label{6.12.00.11}
\partial_{v_p}\circ \delta \{g_0: ... :g_m\} = 
\end{equation}
$$
\partial_{v_p}\Bigl( \sum_{i=0}^{m+1} \{g_{i+1}: g_{i+2}: ... : g_{i-1}\} 
\wedge \{g_{i-1}: g_{i}\} + 
\{g_{i-1}: g_{i}\} \wedge \{g_{i}: g_{i+1}: ... : g_{i-2}\}\Bigr)
$$
Since $\{g_{i-1}: g_{i}\} =0$ if $g_{i-1}= g_i$ we calculate 
(\ref{6.12.00.11}) as follows. 
We locate $g_0, ..., g_{m+1}$ cyclically on the circle. Say that a string in this cyclic word 
is a sequence 
of letters 
$g_i, g_{i+1}, ..., g_j$ such that $g_i= g_{i+1}= ...= g_j$ and $g_{i-1} \not = g_i$, 
$g_{j} \not = g_{j+1}$. The cyclic word splits 
into a union of  strings. For instance if all $g_i$'s are distinct, 
we get $m+1$ one element strings. The sum (\ref{6.12.00.11}) equals to  
$$
\sum_{{\rm strings}}\partial_{v_p} \Bigl(\{g_{i+1}: g_{i+2}: ... : g_{i-1}\} \wedge \{g_{i-1}: g_{i}\} \Bigr)+ 
$$
$$
\sum_{{\rm strings}}\partial_{v_p} \Bigl(\{g_{j}: g_{j+1}\} \wedge 
\{g_{j+1}: ... : g_{i-1}: g_{i}\} \Bigr)
$$
The sum of the two terms corresponding to a given string is obviously zero, 
so (\ref{6.12.00.11}) is zero. The proposition is proved. 

It follows from proposition \ref{6.1.00.100} that the map 
$$
\partial_{v_p}: \Lambda^*_{(m)}{\cal D}_{\bullet}(\mu_p)  \lra
\Lambda^*_{(m-1)}{\cal D}_{\bullet}(\mu_p) [-1]  
$$
is a homomorphism of complexes.

{\bf 9. The depth $2$ part of the cochain complex of ${\cal D}^{{\rm un}}_{\bullet}(\mu_p)$ 
and the cuspidal cohomology of $\Gamma_1(2; p)$}.  
\begin{theorem}\label{6.12.00.3}
There is a canonical isomorphism of complexes
$$
{\cal D}_{ 2}(\mu_p) \lra 
\Lambda^2 {\cal D}_{ 1}^{{\rm un}}(\mu_p) \quad = \quad \tau_{[1,2]}\Bigl({\Bbb M}_{(2)}^* \otimes_{\Gamma_1(2;p)} \Q\Bigr)
$$
\end{theorem}

{\bf Proof}. This theorem is a version of the results proved in s. 3.5-3.6 in [G2]. 
Let us extend $v_p$ to $\widehat {\cal D}_{ 1}(\mu_p)$ by putting it zero on $\Q_{(1,1)}$. 
We define a map 
$$
\Lambda^2 \widehat {\cal D}_{ 1}(\mu_p) \quad \stackrel{\delta'}{\lra } 
\quad 
{\cal D}_{ 1}(\mu_p) \oplus {\cal D}_{ 1}(\mu_p)
$$
by setting $\delta':= (-\partial, \partial + v_p)$, i.e. 
$$
\{1: \zeta_p^{\alpha}\} \wedge \{1: \zeta_p^{\beta}\} \lms  
\left\{ \begin{array}{ll}
0 \quad \oplus  \quad  \{1: \zeta_p^{\beta}\}  - \{1: \zeta_p^{\alpha}\} 
& \alpha, \beta \not = 0 \\ 
 -\{1: \zeta_p^{\beta}\}\quad \oplus  \quad  \{1: \zeta_p^{\beta}  \}&  \alpha = 0 , \beta \not = 0  \end{array} \right.
$$
Thus we get  the following complex, placed in degrees $[1,3]$:
\begin{equation} \label{6.10.00.4}
{\cal D}_{ 2}(\mu_p) \quad \stackrel{\delta}{\lra } \quad 
\Lambda^2 \widehat {\cal D}_{ 1}(\mu_p) \quad \stackrel{\delta'}{\lra } 
\quad 
{\cal D}_{ 1}(\mu_p) \oplus {\cal D}_{ 1}(\mu_p)
\end{equation}
\begin{theorem} \label{6.12.00.2}
The complex (\ref{6.10.00.4}) is canonically identified with the complex 
${\Bbb M}_{(2)}^* \otimes_{\Gamma_1(2;p)} \Q$. 
\end{theorem}

{\bf Proof}. We start with a lemma describing more explicitly 
the rank $2$ modular complex. 

\begin{lemma} \label{6.10.00.5} The double shuffle relations for $m=2$ are {\rm equivalent} 
to the  dihedral relations:
 $$
<v_0, v_1,v_2> = - <v_1,v_0,v_2>  = - <v_0, v_2,v_1> = <-v_0,-v_1,-v_2>
$$
\end{lemma}

{\bf Proof}. The double shuffle relations in $<\cdot ,\cdot ,\cdot >$ 
generators look as follows:
$$
<v_0, v_1,v_2> + <v_1,v_0,v_2> =0, \quad <v_1, v_2-v_1,-v_2> + <v_2, v_1-v_2, -v_1> =0
$$
Changing variables $u_0 = v_1, u_1 = v_2-v_1, u_2 = -v_2$ we write the second of them as $<u_0, u_1, u_2> + <-u_2, -u_1, -u_0> =0$. The lemma follows.

It follows from this lemma that the modular complex for $GL_2(\Z)$ 
is canonically isomorphic to the following complex of left $GL_2(\Z)$-modules:
$$
\Z[GL_2(\Z)]\otimes_{D_2} \xi_2 \stackrel{\partial}{\lra} 
\Z[GL_2(\Z)]\otimes_{D_{1,1}}\xi_{1,1} \stackrel{\partial'}{\lra} \Z[GL_2(\Z)/\widetilde B]
$$
where $D_2 \subset GL_2(\Z)$ is the order $12$ dihedral subgroup, $$
D_{1,1} = 
\{\left (\matrix{\pm 1& 0\cr 
0& \pm 1 \cr}\right ), \left (\matrix{0& \pm 1\cr 
 \pm 1 &0\cr}\right ) \},\qquad 
\widetilde B:= \left (\matrix{\pm 1&* \cr 
0& \pm 1 \cr}\right )
$$ 
$\xi_2$  is the character of $D_2$ given by the determinant, and  
$\xi_{1,1}$ is a nontrivial character of $D_{1,1}$ killing $\left (\matrix{\pm 1& 0\cr 
0& \pm 1 \cr}\right )$. 
The differentials commute with the left action of the group $GL_2(\Z)$ 
and so are determined by their action on the unit in $GL_2$. They are  given by 
$$
 \partial: \left (\matrix{1&0 \cr 
0&1\cr}\right )\lms \quad  -\left(\matrix{1&0 \cr 
0&1\cr}\right ) -  
\left(\matrix{0&-1 \cr 1&-1\cr}\right ) 
- \left(\matrix{-1&1 \cr -1&0\cr}\right )
$$
$$
 \partial': \left (\matrix{1&0 \cr 
0&1\cr}\right )\lms \quad  \left(\matrix{0&1 \cr 1&0\cr}\right )  - \left(\matrix{1&0 \cr 
0&1\cr}\right )  
$$
The coset 
$
\Gamma_1(2;p)\backslash GL_2(\Z)$ is identified with the set of rows 
$\{(\alpha, \beta) \in F_p^2- 0\}$. One can identify it with the set 
of non zero rows $\{(\alpha, \beta, \gamma) \}$ with $\alpha + \beta + \gamma = 0$. 
Then 
$$
M_{(2)}^1 =  \Z[\Gamma_1(2;p)\backslash GL_2(\Z)]\otimes_{D_2} \xi_2  =
\frac{\Z[(\alpha, \beta, \gamma)\in F_p^3- 0 \quad | \quad \alpha + \beta + \gamma = 0]}{ \mbox{the dihedral relations}}
$$

$$
M_{(2)}^2 \quad = \quad \Z[\Gamma_1(2;p)\backslash GL_2(\Z)]\otimes_{D_{1,1}} \xi_{1,1} \quad  =
\quad \frac{\Z[ (\alpha, \beta)\in F_p^2-0]}{ (\alpha, \beta) =  
-(\beta, \alpha) = (\pm \alpha, \pm \beta) }
$$

$
M_{(2)}^3 \quad = \quad \Z[\Gamma_1(2;p)\backslash GL_2(\Z)/ \widetilde B] \quad  =\quad
$ the abelian group with the generators $[\beta, 0]$ and $[0, \beta]$,  where 
$\beta \not = 0$, and the only relation is symmetry under $\beta \lms -\beta$. 

The complement $X_1(p) - Y_1(p)$ consists of $p-1$ cusps. 
The natural covering $X_1(p) \lra X_0(p)$ is 
unramified of degree $(p-1)/2$. There are two cusps on $X_0(p)$: 
the $0$ and $\infty$ cusps. So there are $(p-1)/2$ cusps over $0$ and  over $\infty$. 
Under the identification with 
${\cal D}_{ 1}(\mu_p) \oplus {\cal D}_{ 1}(\mu_p)$ they correspond to the summands 
${\cal D}_{ 1}(\mu_p)$.

The desired isomorphism of complexes is defined by 
$$
(\alpha, \beta, \gamma) \lms  
\{\zeta_p^{\alpha}, \zeta_p^{\beta},\zeta_p^{\gamma}\}
$$
$$
(\alpha, \beta) \lms  
\{\zeta_p^{\alpha}, \zeta_p^{-\alpha}\} \wedge \{\zeta_p^{\beta}, \zeta_p^{-\beta}\}
$$
$$
[\beta, 0] \oplus [0, \beta'] \lms 
\{\zeta_p^{\beta}, \zeta_p^{-\beta}\} \oplus \{\zeta_p^{\beta'}, \zeta_p^{-\beta'}\} 
$$

To check that it is indeed an isomorphism of the vector spaces notice that, 
similarly to lemma \ref{6.10.00.5}, the only relations among the elements
$\{\zeta_p^{\alpha}, \zeta_p^{\beta},\zeta_p^{\gamma}\}$ are the dihedral symmetry. 
Notice also that $\{1,1,1\}=0$ and the only distribution relation
 $
0= \{1, 1, 1\} = 
\sum_{\alpha, \beta}\{\zeta_p^{\alpha}, \zeta_p^{\beta},\zeta_p^{-\alpha - \beta}\}
$ follows from the skew symmetry. 

Theorem \ref{6.12.00.2} is proved. 
Theorem \ref{6.12.00.3} follows immediately from theorem \ref{6.12.00.2}.

{\bf 10. The depth $3$ part of $\Lambda^*{\cal D}^{{\rm un}}_{\bullet}(\mu_p)$}. 

\begin{theorem} \label{6.11.00.1}
\begin{equation} \label{6.10.00.1}
H^i_{(w, 3)}({D}^{{\rm un}}_{\bullet}(\mu_p)) = \quad \left\{ \begin{array}{llll}
0&  i=1 \\ 
H_{{\rm cusp}}^{2}(\Gamma_1(3;p), \Q) \oplus H_{{\rm cusp}}^{1}(\Gamma_1(2;p), V_2)&  i=2 \\ 
  H_{{\rm cusp}}^{3}(\Gamma_1(3;p), \Q)  &  i=3\end{array} \right.
 \end{equation}
\end{theorem}

{\bf Proof}. 
There  is an exact sequence of complexes
$$
\begin{array}{ccccc}
{\cal D}_{3}(\mu_p) &\stackrel{\delta}{\lra} & 
{\cal D}_{2}(\mu_p)\otimes {\cal D}^{{\rm un}}_{1}(\mu_p)   
&\stackrel{\delta}{\lra} & 
\Lambda^3{\cal D}^{{\rm un}}_{ 1}(\mu_p) \\
\downarrow =&& \downarrow   && \downarrow \\
{\cal D}_{3}(\mu_p) &\stackrel{\delta}{\lra}
&{\cal D}_{2}(\mu_p) \otimes {\cal D}_{1}(\mu_p)& \stackrel{\delta}{\lra} &\Lambda^3
{\cal D}_{1}(\mu_p) \\
&& \downarrow \partial_{v_p}&& \downarrow \partial_{v_p}  \\
&&{\cal D}_{2}(\mu_p)&\stackrel{\delta}{\lra} & \Lambda^2 {\cal D}^{{\rm un}}_{1}(\mu_p) 
\end{array}
$$
A choice of splitting $\Q \lra {\cal D}_{1}(\mu_p) $ of the map $v_p$ provides 
a splitting of the middle complex into a direct sum of the top and bottom complexes. 
Taking into account canonical decomposition (\ref{4-16.10}) we end up with 
$$
\Lambda^*_{(3)}(\widehat {\cal D}_{\bullet}(\mu_p)) \quad = \quad
\Lambda^*_{(3)}({\cal D}^{{\rm un}}_{\bullet}(\mu_p)) \quad \oplus \quad 
\Lambda^*_{(2)}({\cal D}^{{\rm un}}_{\bullet}(\mu_p))[-1] \quad \oplus \quad 
\Lambda^*_{(2)}({\cal D}_{\bullet}(\mu_p))[-1]
$$
Notice also a noncanonical splitting 
$$
\Lambda^*_{(2)}({\cal D}_{\bullet}(\mu_p)) \quad = \quad
\Lambda^*_{(2)}({\cal D}^{{\rm un}}_{\bullet}(\mu_p)) \quad \oplus \quad 
{\cal D}^{{\rm un}}_{1}(\mu_p)[-1]
$$
Therefore we have
$$
\Lambda^*_{(3)}(\widehat {\cal D}_{\bullet}(\mu_p)) \quad = \quad 
\Lambda^*_{(3)}({\cal D}^{{\rm un}}_{\bullet}(\mu_p)) \quad \oplus \quad 
\Lambda^*_{(2)}({\cal D}^{{\rm un}}_{\bullet}(\mu_p))\otimes \Q^2 [-1]
\quad \oplus  \quad
{\cal D}^{{\rm un}}_{\bullet}(\mu_p)[-2]
$$

This implies that 
$$
H^1_{{\rm cusp}}(\Gamma_1(2; p), \varepsilon_2) \otimes \Q^2 \quad \oplus \quad 
\Q^{\frac{p-3}{2}}
$$
is a direct summand of $H_{{\rm inf}}^3(\Gamma_1(3; p), \Q)$. 
This has been known to us through 
theorem \ref{6.19.00}. On the other hand this gives an additional 
confirmation that the result of [LS] concerning 
$H^3_{{\rm inf}}(\Gamma(3; p), \Q)$ has to be corrected, as was explained in s. 6.3.

Using the relation (\ref{d1f_})-(\ref{d1f_qeq}) between the 
depth $2$ and $3$ pieces of the cohomology of the Lie coalgebra 
$\widehat {\cal D}_{\bullet}(\mu_p)$ and the cohomology of 
groups $\Gamma_1(2; p)$ and $\Gamma_1(3; p)$,  
and combining this with  the description of the cohomology of the group $\Gamma_1(3; p)$ 
given in the previous chapter, in particular theorem \ref{6.19.00},  we arrive to the proof of theorem \ref{6.11.00.1}.

{\bf Remark}. Here are some cycles in $H^2(\widehat {\cal D}_{\bullet}(\mu_p))$:
$$
\{\zeta_p^{\alpha}, \zeta_p^{-\alpha}\} \otimes \{\zeta_p^{\alpha}, \zeta_p^{\beta}, \zeta_p^{\gamma} \} + 
\{\zeta_p^{\beta}, \zeta_p^{-\beta}\} \otimes \{\zeta_p^{\beta}, \zeta_p^{\alpha},  \zeta_p^{\gamma}\}
$$
They probably generate $H^2$, and one should be able to prove 
directly using these cycles that $H^2 (\widehat {\cal D}_{\bullet}(\mu_p)) = 
H^1_{{\rm cusp }}(\Gamma_1(2; p), V_2) $ 

{\bf 11.  Proof of corollary \ref{4-19.10}}. 
Let us introduce the following notation. Let ${\cal D}_{\bullet}$ be a $\Z_{+}$-graded 
Lie algebra. Denote by $\chi_{(m)}({\cal D}_{\bullet})$ the Euler characteristic of the degree $m$ part of the standard cochain complex of ${\cal D}_{\bullet}$ 

\begin{corollary} \label{CORZ}
For $m=2, 3$ one has 
\begin{equation} \label{6.10.00.2}
\chi_{(m)}({\cal D}^{{\rm un}}_{\bullet}(\mu_p)) = \left\{ \begin{array}{llll}
{\rm dim}H_{{\rm cusp}}^{1}(\Gamma_1(2;p),  \varepsilon_2) &  m=2 \\ 
{\rm dim}H_{{\rm cusp}}^{1}(\Gamma_1(2;p), V_2)&  m=3 \end{array} \right.
\end{equation}
\end{corollary}

{\bf Proof}. For $m=2$ this was done before. For $m=3$ this follows from the above theorem and the Poincar\'e duality for cuspidal cohomology for $\Gamma_1(3; p)$.

{\it Corollary \ref{4-19.10}: the $m=2$ case}. Thanks to lemma \ref{modcoho} and theorem \ref{6.12.00.3}
$$
{\rm dim}{\cal D}_{2}(\mu_p)   - {\rm dim}\Lambda^2 
{\cal D}^{{\rm un}}_{1}(\mu_p) \quad = \quad 0
-  {\rm dim} H^1_{{\rm cusp}}(\Gamma_1(2;p), \varepsilon (2)))
$$
One has 
$$
{\rm dim}{\cal D}^{{\rm un}}_{1}(\mu_p) = \frac{p-3}{2};\quad => \quad  
{\rm dim}\Lambda^2{\cal D}^{{\rm un}}_{1}(\mu_p) =  \frac{(p-3)(p-5)}{8}
$$
Using proposition 1.40 in [Sh] we have, for $p \geq 5$, 
$$
{\rm dim}H^1(X_1(p), \Q)^- = {\rm dim}H_{{\rm cusp}}^1(\Gamma_1(2;p), \varepsilon (2))
 = 1 + \frac{p^2-1}{24}-\frac{p-1}{2}
$$
(In our case, using the notations of [Sh], $\mu = \frac{p^2-1}{2}$, 
$\nu_2 = \nu_3 =0, \nu_{\infty} = p-1$ and $\overline \Gamma_1(p):= \Gamma_1(p)/\pm {\rm Id}$). 
Therefore
$$
{\rm dim}{\rm Gr}{\cal L}^{(l)}_{-2}(\mu_p)  \quad \stackrel{(\ref{4-15.20})}{=} \quad 
{\rm dim}{\cal D}^{}_{2}(\mu_p)\quad = \quad \frac{(p-5)(p-1)}{12}
$$

{\it Corollary \ref{4-19.10}: the $m=3$ case}. We use corollary \ref{CORZ} and the following result ([Sh], theorem 2.25)
$$
{\rm dim}H_{{\rm cusp}}^{1}(\Gamma_1(2;p), V_2\otimes \varepsilon_2) \quad =\quad
2 \Bigl(\frac{p^2-1}{24} - \frac{p-1}{2} \Bigr) + \frac{p-1}{2}
$$
Using this and (\ref{4-15.20}) computation of the Euler characteristic gives 
$$
{\rm dim}{\rm Gr}{\cal L}^{(l)}_{-3}(\mu_p) \quad \stackrel{(\ref{4-15.20})}{=}\quad  {\rm dim}{\cal D}_3(\mu_p) \quad 
$$
$$
 \frac{(p-1)(p-5)(p-3)}{24} 
-\frac{(p-3)(p-5)(p-7)}{48} - {\rm dim}H_{{\rm cusp}}^{1}(\Gamma_1(2;p), V_2
\otimes \varepsilon_2) = 
$$
$$ 
\frac{(p-5)(p^2-2p-11)}{48} 
$$
I believe that there exists an explicit degree $m$ polynomial in $p$  giving 
for ${\rm dim }{\cal D}_{m}(\mu_p)$ for all $m$. There should be a formula for 
$\chi_m(\widehat {\cal D}_{\bullet}(\mu_p))$ similar to (\ref{6.10.00.2}). 
Notice that there is no closed formula 
for ${\rm dim }H_{{\rm cusp}}^{2}(\Gamma_1(3;p))$. 
 
The computation of $\chi_4(\widehat {\cal D}_{\bullet}(\mu_p))$ 
seams to be a  feasible problem for the following reasons. 
First, according to our recent results ([G4]) the relation between the depth $m$ part of the 
cohomology of the Lie coalgebra $\widehat {\cal D}_{\bullet}(\mu_p)$ and 
cohomology of $\Gamma_1(m;p)$ is still in case for $m=4$: 
$$
H_{(4)}^i(\widehat {\cal D}_{\bullet}(\mu_p)) \quad = \quad H^{i+2}(\Gamma_1(4; p), 
\varepsilon_4); \qquad 
1 \leq i \leq 4
$$
Here $\varepsilon_4$ is the character of $\Gamma_1(4; p)$ given by the determinant. Therefore
$$
\chi_4(\widehat {\cal D}_{\bullet}(\mu_p)) = 
\sum_{i=3}^6 (-1)^i {\rm dim} H^i(\Gamma_1(4; p), \varepsilon_4)
$$
Further, the cuspidal cohomology of $\Gamma_1(4; p)$ are in the degrees 
$4$ and $5$, and the rest of the cohomology is given by the cohomology at infinity 
$H^*_{{\rm inf}}$, i.e. 
$$
0 \lra H^*_{{\rm cusp}}(\Gamma_1(4; p), \varepsilon_4)\lra H^*(\Gamma_1(4; p), \varepsilon_4)\lra H^*_{{\rm inf}}(\Gamma_1(4; p), \varepsilon_4) \lra 0
$$
 is an exact sequence.
 This can be proved using the following facts: 
besides the trivial representation 
there is only one unitary cohomological representation of $SL_4(\R)$ 
(the so called special representation) and 
this representation is tempered (see theorems 5.6 and 6.16 
in [VZ] for the general results). Then we  
employ 
 the arguments similar to the one given in chapter 1 of [LS]. 
Since the cuspidal cohomology satisfy the Poincar\'e duality we get
$$
\chi_4(\widehat {\cal D}_{\bullet}(\mu_p)) = 
\sum_{i=3}^6 (-1)^i {\rm dim} H_{{\rm inf}}^i(\Gamma_1(4; p), \varepsilon_4)
$$

{\bf Problem}. Compute of the right hand side of this formula. 

I would expect that the cuspidal cohomology of $\Gamma_1(3; p)$ will not appear in the answer, 
so we should have an explicit formula for it.

{\bf 12. The Galois groups unramified at $p$. } Let us assume 
that $p \not = l$. 
Let ${\rm Gr}{\cal G}^{(l)}_{-m}(\mu_p)^{{\rm un}}$ be the maximal quotient of 
${\rm Gr}{\cal G}^{(l)}_{-m}(\mu_p)$ 
 unramified at the place $1-\zeta_p$. 
It is easy to see that 
$$
{\rm Gr}{{\cal G}^{(l)}_{-1}}(\mu_p)^{{\rm un}} = 
{\rm Hom}\Bigl(\mbox{the group of the cyclotomic units in $
\Z[\zeta_p]$}, \quad \Q_l\Bigr)
$$
It follows from theorem \ref{6.12.00.3} that 
$$
{\rm Gr}{{\cal G}^{(l)}_{ -2}}(\mu_p)^{{\rm un}} = {\rm Gr}{{\cal G}^{(l)}_{ -2}}(\mu_p)\qquad {\rm Gr}{{\cal L}^{(l)}_{ -3}}(\mu_p)^{{\rm un}} = {\rm Gr}{{\cal L}^{(l)}_{ -3}}(\mu_p)
$$

{\bf Acknowledgement}. 
  I am  grateful to M. Kontsevich and G. Harder for very useful discussions. 
The initial  draft of this paper was submitted in 
1998 to Mathematical Research Letters. 
I am very much indebted to the referee
for correcting  some gaps and errors 
and making  many suggestions, 
 incorporated into the text, which greatly improved and clarified
 the exposition. 

The initial draft of this paper was written 
during the Summer and Fall 1998 in MPI (Bonn) and  reported 
at the MPI Conference on Polylogarithms at Schloss Ringberg (August 1998). 
I am  grateful to the Max-Planck-Institute  
  for hospitality and
 support. 
The  support by the NSF grants  DMS-9500010 and DMS-9800998  is gratefully
acknowledged. Finally I would like to thank IHES for hospitality and support during 
the final stage of work on the paper.

\vskip 5mm \noindent
{\bf Corrections to the paper [G3]}

1. In formulas (10), (35), (36) and formula one line before (35) 
one needs to replace $S^{w-2}V_2$ by 
$S^{w-2}V_2 \otimes \varepsilon_2$, where $\varepsilon_2$ is as above. 

2. In theorem 6.1 $\Gamma$ is  a finite index subgroup of $GL_2(\Z)$; 
replace $H^{i-1}(\Gamma, V)$ by $H^{i-1}(\Gamma, V\otimes \varepsilon_2)$.

3. Section 4.1, first two lines:  use $\Z[X]$ instead of $\Z[[X]]$; line four:  
$\Z[{\cal P}_m]$ instead of $\Z[[{\cal P}_m]]$.

4. In theorems 6.2 $\Gamma$ is  a finite index subgroup of $GL_3(\Z)$.

\vskip 5mm \noindent
{\bf REFERENCES}
\begin{itemize}
\item[{[BMS]}] Bass H., Milnor J., Serre J.-P.: {\it Solution of the congruence 
subgroup problem in $SL_n (n\geq 3)$ and $Sp_{2n} (n\geq 2)$}, Publ. Math. 
I.H.E.S. 33 (1967), 59-137. 
\item[{[Be]}] Beilinson A.A.: {\it Polylogarithms and cyclotomic elements} 
MIT Preprint, 1989. 
\item[{[BD]}] Beilinson A.A., Deligne P.: {\it Motivic polylogarithms and Zagier's 
conjecture} Manuscript (version of 1992). 
\item[{[B]}] : Borel A.: {\it Stable real cohomology of arithmetic groups}. 
Ann. Sci. École Norm. Sup. (4) 7 (1974), 235--272. 
\item[{[BW]}] : Borel A., Wallach N.: {\it Continuous cohomology, discrete 
subgroups, and representations of reductive groups} (Second edition), AMS 
Math. Surveys 
and Monographs, vol. 67, 
2000. 
 \item[{[BS]}] : Borel A., Serre J. P.: {\it Corners and arithmetic groups}. 
Comm. Math. Helv. 48 (1973), 436-491. 
\item[{[D]}] Deligne P.: {\it Le group fondamental de la droite projective moins trois points}, In: Galois groups over $\Bbb Q$. Publ. MSRI, no. 16 (1989) 79-298.  
\item[{[D1]}] Deligne P.: {\it Letter to D. Broadhurst}. (1997).
\item[{[D2]}] Deligne P.: {\it Letter to the author}. (2000).
\item[{[Dr]}] Drinfeld V.G.: {\it On quasi-triangular quasi-Hopf algebras and some group related to  Gal$(\overline{\Bbb Q}/\Bbb Q)$}. Leningrad Math. J., 1991. 
 \item[{[G1]}] Goncharov A.B.: {\it Multiple $\zeta$-numbers,
    hyperlogarithms and mixed Tate motives}, Preprint MSRI 058-93, June 1993.
\item[{[G2]}] Goncharov A.B.: {\it The double logarithm and Manin's
complex for modular curves}. Math. Res. Letters,
 vol. 4. N 5 (1997), pp. 617-636. 
\item[{[G3]}] Goncharov A.B.: {\it Multiple polylogarithms, cyclotomy and modular complexes
},     Math. Res. Letters,
 vol. 5. (1998), pp. 497-516.
\item[{[G4]}] Goncharov A.B.: {\it Multiple polylogarithms at roots of unity, 
and geometry of modular varieties} In preparation. 
\item[{[G5]}] Goncharov A.B.: {\it Galois groups, geometry of modular varieties and graphs},   
Proc. of Arbeitstagung, June 1999, in the MPI preprint 
(http//www.mpim-bonn.mpg.de/) 13 pages. 
\item[{[G6]}] Goncharov A.B.: {\it Multiple $\zeta$-values, Galois groups and geometry of modular varieties}. Proc. of the Third European Congress of Mathematicians, Barcelona, 2000. Posted on alg-geom web, AG/0005069. 
\item[{[G7]}] Goncharov A.B.: 
{\it Modular complexes and geometry of symmetric spaces}. To appear. 
\item[{[G8]}] Goncharov A.B.: 
{\it Mixed elliptic motives}, in London Math. Soc. Lect. Note Series, 254, 
Cambridge Univ. Press, Cambridge, 1998, 147-221. 
\item[{[G9]}] Goncharov A.B.: 
{\it Polylogarithms in arithmetic and geometry}, 
Proc. of the Int. Congress of
Mathematicians, Vol. 1, 2 (Zürich, 1994), 374--387, Birkhäuser, Basel, 1995.
\item[{[Gr]}] Grothendieck A.: {\it Esquisse d'un programme}. Mimeographed note (1984).
\item[{[Ha]}] Harder G.: {\it A letter to the author}, (1999).
\item[{[HM]}] Hain R., Matsumoto M.: {\it Weighted Completion of Galois Groups and Some Conjectures of Deligne}. Preprint June 2000. 
\item[{[HW]}] Huber A, Wildeshaus J.: {\it Classical motivic polylogarithms according to Beilinson and Deligne} Doc. Math. 3 (1998), 27-133.
\item[{[Ih]}] Ihara Y.: {\it Braids, Galois groups, and some arithmetic functions}. Proc. Int. Congress of Mathematicians in Kyoto, (1990). 
\item[{[Ih1]}] Ihara Y.: {\it Galois representation arising from
${\Bbb P}^1- \{0,1,\infty\}$ and Tate twists of even degree},
In: Galois groups over $\Bbb Q$. Publ. MSRI, no. 16 (1989).   
\item[{[Ih2]}] Ihara Y.: {\it Profinite braid groups, Galois representations and complex multiplication   },
 Ann. of Math. 123, (1986)  43-106. 
\item[{[Ih3]}] Ihara Y.: {\it Some arithmetical aspects of Galois 
action on the pro-$p$ fundamental group of 
$\widehat \pi_1({\Bbb P}^1 - \{0,1,\infty\})$}, Preprint RIMS-1229, 1999. 
\item[{[K]}] Kontsevich M.: {\it Formal (non)commutative symplectic 
geometry} The Gelfand mathematical seminars, Birkhauser, 1993, p. 173-187.
\item[{[LS]}] Lee R., Schwermer J.: {\it Cohomology of arithmetic subgroups of ${\rm SL}\sb{3}$ at 
infinity} J. Reine Angew. Math. 330 (1982), 100--131. 
\item[{[So]}] Soul\'e C.: {\it On higher $p$-adic regulators}, Springer LNM 854, (1981), 372-401.
 \item[{[Sh]}] Shimura G.: {\it Introduction to the 
arithmetic theory of automorphic functions}, Publ. Math. Soc. Japan, 11, 
Princeton, 1971. 
\item[{[VZ]}] Vogan D., Zuckerman G,: {\it Unitary representations with non zero cohomology},
 Compositio Math 53 (1984), 51-90.
 \item[{[Z1]}] Zagier D.: {\it Values of zeta functions and their
    applications}, Proceedings of the First European Congress of
  Mathematicians, Paris. 1994, vol 1. 

\end{itemize}

\end{document}